\documentclass[11pt]{article}
\usepackage{amssymb}
\usepackage{amsmath}
\usepackage{bm}
\usepackage{multirow}	
\usepackage{enumerate}
\usepackage{graphicx}
\usepackage{color}
\usepackage{hyperref}
\usepackage{booktabs}
\usepackage{caption}
\usepackage{placeins}

\usepackage{xcolor}
\hypersetup{
    colorlinks,
    linkcolor=blue,
    citecolor=blue,
    urlcolor=blue
}
\usepackage{comment} 
 
\usepackage{epstopdf}
\numberwithin{equation}{section}

\usepackage{amsthm}
\newtheorem{theorem}{Theorem}[section]

\newtheorem{lemma}[theorem]{Lemma}

\newcommand\R{\mathrm{R}}
\newcommand{\tr}{\mathrm{tr}}

\renewcommand{\O}{\Omega}

\newcommand\qan{\quad\hbox{and}\quad}
\renewcommand{\div}{\mathrm{div}}

\newcommand\bsi{{\text{\boldmath $\sigma$}}}
\newcommand\bta{{\text{\boldmath $\tau$}}}
\newcommand\bze{{\text{\boldmath $\zeta$}}}

\newcommand\bt{{\text{\boldmath $t$}}}
\newcommand\bs{{\text{\boldmath $s$}}}

\newcommand\dist{\mathrm{dist}}

\newcommand\calA{\mathcal{A}}
\newcommand\calB{\mathcal{B}}
\newcommand\calC{\mathcal{C}}
\newcommand\calD{\mathcal{D}}
\newcommand\calE{\mathcal{E}}

\newcommand\bw{\boldsymbol{w}}

\newcommand\bfa{\mathbf{a}}
\newcommand\bfb{\mathbf{b}}

\newcommand\bfe{\mathbf{e}}
\newcommand\bfd{\mathbf{d}}

\newcommand\ba{{\text{\boldmath $a$}}}

\newcommand{\bPi}{{\boldsymbol{\Pi}}}

\newcommand{\bfA}{\mathbf{A}}

\newcommand\bx{{\text{\boldmath $x$}}}

\newcommand\bxi{{\text{\boldmath $\xi$}}}
\newcommand\bnu{{\text{\boldmath $n$}}}
\newcommand\bvphi{{\text{\boldmath $\varphi$}}}
\newcommand\bpsi{{\text{\boldmath $\psi$}}}
\newcommand\bphi{{\text{\boldmath $\phi$}}}
\newcommand\bro{{\text{\boldmath $\rho$}}}
\newcommand\bet{{\text{\boldmath $\eta$}}}
\newcommand\balph{{\text{\boldmath $\alpha$}}}

\newcommand\bep{\boldsymbol{\varepsilon}}
\renewcommand\u{{\text{\boldmath $u$}}}
\newcommand\bv{{\text{\boldmath $v$}}}

\newcommand\br{{\text{\boldmath $r$}}}

\newcommand\z{{\text{\boldmath $z$}}}
\newcommand{\dd}{{\mathrm d}}
\newcommand{\I}{\mathbb{I}}

\newcommand{\qon}{{\quad\hbox{on}\quad}}
\newcommand\qin{{\quad\hbox{in}\quad}}

\renewcommand\bet{{\text{\boldmath $\eta$}}}
\newcommand{\bdiv}{\mathbf{div}}
\newcommand{\rL}{\mathrm{L}}
\newcommand{\bL}{\mathbf{L}}

\newcommand{\bH}{\mathbf{H}}
\newcommand{\bbH}{\mathbb{H}}
\newcommand{\Pl}{\mathrm{P}_\ell}
\newcommand{\bPl}{\mathbf{P}_\ell}
\newcommand{\RTl}{\mathbf{RT}_\ell}

\newcommand{\Th}{\mathcal{T}_h}

\def\DOF{\mathtt{DoF}}

\newcommand{\bXi}{\boldsymbol{\Xi}}

\newcommand{\rH}{\mathrm{H}}

\newcommand{\rW}{\mathrm{W}}

\newcommand{\calM}{{\mathcal{M}}}
\newcommand{\bcalM}{{\boldsymbol{\mathcal M}}}
\newcommand{\calH}{\mathcal H}
\newcommand{\bcalH}{{\boldsymbol{\mathcal H}}}

\newcommand{\rX}{\mathrm{X}}
\newcommand{\rM}{\mathrm{M}}

\newcommand{\bQ}{\mathbf{Q}}

\newcommand{\bF}{\mathbf{F}}
\newcommand{\bG}{\mathbf{G}}

\newcommand{\bW}{\mathbf{W}}

\newcommand{\bX}{\mathbf{X}}

\newcommand{\bS}{\mathbf{S}}
\newcommand{\bT}{\mathbf{T}}

\newcommand{\tbS}{{\widetilde{\bS}}}

\newcommand{\ds}{\displaystyle}
\newcommand\rP{{\mathrm P}}
\newcommand\bP{{\mathbf P}}

\theoremstyle{definition}

\theoremstyle{remark}
\newtheorem{remark}[theorem]{Remark}
\numberwithin{equation}{section}
\numberwithin{table}{section}
\numberwithin{figure}{section}
\allowdisplaybreaks

\definecolor{lightblue}{rgb}{0.22,0.45,0.70}
\definecolor{lightgreen}{rgb}{0.22,0.50,0.25}

\newcommand\ff{\boldsymbol{f}}
\newcommand{\cE}{\mathcal{E}}
\newcommand{\jump}[1]{[\![#1]\!]}

\newcommand{\bzeta}{\boldsymbol{\zeta}}
\usepackage{relsize}
\begin{document}

\title{A pseudostress-based fully mixed finite element method
for the stationary Biot--Poisson--Nernst--Planck system
\thanks{This research was partially supported by ANID-Chile
through Fondecyt Postdoctoral project 3260398.}}

\author{{\sc Cristian Inzunza}
\thanks{Departamento de Matem\'aticas, Facultad de Ciencias,
Universidad de La Serena, Juan Cisternas 1200, La Serena, Chile,
email: {\tt cr.inzunza.d@gmail.com}.}}

%\date{Updated: \today}

\date{ }

\maketitle

\begin{abstract}
\noindent
We develop a pseudostress-based fully mixed finite element method
for the stationary Biot--Poisson--Nernst--Planck system describing
ionic transport in a charged, saturated, and deformable porous
medium. The mechanical equations are reformulated in terms of a
non-symmetric pseudostress and the solid displacement, while mixed
formulations are retained for the Darcy, electrostatic, and ionic
equations. The nonlinear electrostatic and advective couplings
are treated in a Banach space setting. The formulation provides
direct approximations of the pseudostress, fluxes, and concentration
gradients, and permits algebraic reconstruction of the displacement
gradient and Cauchy stress. Under block-stability and small-data
assumptions, we establish continuous and discrete well-posedness,
with poroelastic stability bounds uniform with respect to the first
Lam\'e parameter. A conforming Raviart--Thomas discretization,
completed by discontinuous polynomials, yields a priori error
estimates and convergence rates. For the coupled scheme in two
dimensions, we construct a residual-based a posteriori error
estimator and establish its reliability under an additional
smallness condition on the prescribed data. Local efficiency is
proved for polynomial boundary data, with boundary-data oscillation
included for general traces. We also describe the estimator
modifications for nonhomogeneous displacement and normal-flux data
under corresponding stability and solvability hypotheses.
Numerical experiments examine convergence of the mixed and
reconstructed fields, robustness in the nearly incompressible
regime, adaptive refinement, and ionic polarization in a charged
porous medium.
\end{abstract}



\section{Introduction}%\label{sec1}

Ionic transport in deformable porous media couples fluid flow,
mechanical deformation, and electrochemical processes.
Concentrations influence the electric potential and mechanical
stress, while diffusion, electric migration, and fluid motion
determine ionic fluxes.
Homogenized descriptions have been derived for electrolyte-saturated
porous structures \cite{moyne06,allaire17,turja19}, with related
models for charged hydrated tissues and ionic polymer actuators
\cite{levenston98,leronni21,whiteley20}.
These interactions motivate numerical methods that approximate
the mechanical and electrochemical fields together with their
associated stresses, gradients, and fluxes.

Mixed finite element methods provide a natural framework for
these approximations.
In poroelasticity, independent stress and fluid-flux variables
have been used to study nearly incompressible behavior
\cite{ambar20,baerland17,lee18,yi14}.
Mixed formulations in Banach spaces also accommodate
Poisson--Nernst--Planck transport coupled to Stokes and
Navier--Stokes flow \cite{cgrb-M2AN-2023,correa23b}.
The stationary Biot--Poisson--Nernst--Planck system considered
here was analyzed in \cite{gir-CAMWA-2025} using a primal--mixed
poroelastic formulation in terms of displacement, total pressure,
Darcy flux, and pore pressure.
The electrostatic and ionic equations were written in mixed
form, including the electric displacement, concentration
gradients, and total ionic fluxes as unknowns.

The present work introduces a pseudostress-based fully mixed
formulation of this coupled system.
We replace the primal mechanical block of
\cite{gir-CAMWA-2025} with a formulation in terms of a non-symmetric
pseudostress and the solid displacement, while retaining its
Darcy, electrostatic, and ionic formulations.
Since the electrochemical blocks are inherited from
\cite{gir-CAMWA-2025}, we state precisely what is new here:
\begin{itemize}
\item[(i)] a fully mixed mechanical block, based on a non-symmetric
pseudostress and a discontinuous displacement, and its coupling
with the Darcy, Poisson, and ionic blocks, including the
modified storage and poroelastic coupling terms
(Section~\ref{section3}). The model also includes a fixed-charge
density $f_\chi$ and, for $c_0>0$, a mass source $g$ of arbitrary mean, so
that the pressure needs no normalization (Sections~\ref{section2}
and~\ref{section3});
\item[(ii)] continuous and discrete well-posedness of the coupled
problem, in which the stability of the poroelastic block is uniform
in the first Lam\'e parameter (Sections~\ref{section4}
and~\ref{section5});
\item[(iii)] a residual-based a posteriori error estimator for the
fully coupled scheme in two dimensions, with reliability,
local efficiency, and its extension to nonhomogeneous displacement
and normal-flux data (Section~\ref{sec:posteriori-2D});
\item[(iv)] numerical evidence for convergence of the mixed and
algebraically reconstructed fields, robustness in the nearly
incompressible regime, estimator-driven adaptivity, and ionic
polarization in a charged porous medium
(Section~\ref{sec:numerical-results}).
\end{itemize}
The pseudostress avoids a symmetry constraint on the tensor
finite element space, following the approach used for linear
elasticity and stress-assisted diffusion
\cite{ggs-CAMWA-2016,gis-JSC-2022}.
The resulting formulation directly approximates the pseudostress
and allows discontinuous polynomial approximations of the
displacement.
The displacement gradient and Cauchy stress are recovered
algebraically from the mixed fields, without differentiating
the discrete displacement.
Compared with the primal mechanical block of \cite{gir-CAMWA-2025},
this choice trades size for robustness. The mechanical unknowns are the
pseudostress, with its Raviart--Thomas row-wise space, and a discontinuous
displacement, instead of a single continuous displacement field, so the
coupled saddle-point system has more degrees of freedom and a
larger linear algebra cost per Newton step. In exchange, the stability of
the poroelastic block is uniform in the first Lam\'e parameter
(Sections~\ref{section4} and~\ref{section5}), whereas standard
conforming primal discretizations of the displacement are known to
degrade in the nearly incompressible regime, when this parameter
grows, unless they are stabilized or modified. We have not
carried out a head-to-head comparison of computational times with
\cite{gir-CAMWA-2025}.
The coupled and nonlinear nature of the model leads to a Banach
space formulation, as in poroelastic heat transport
\cite{cgir-2024}, chemotaxis--Navier--Stokes
\cite{ccgi-CAMWA-2023}, and Boussinesq flow
\cite{cgm-M2AN-2020}.
We establish stability and solvability of the resulting fully
mixed system by combining block estimates with a fixed-point
argument.
The mechanical and generalized Darcy saddle-point blocks,
generalized Poisson block, and two ionic twofold saddle-point
problems are analyzed using classical and generalized
Babu\v{s}ka--Brezzi theories and Banach-space perturbation
arguments \cite{cg-CAMWA-2022}.
Under block-stability and small-data assumptions, the coupled
problem has a unique solution in a prescribed ball, with
poroelastic stability bounds uniform with respect to the first
Lam\'e parameter.
A conforming discretization based on Raviart--Thomas elements
and discontinuous polynomials yields discrete solvability,
a priori error estimates, and convergence rates under uniform
block stability and discrete small-data conditions.

We also construct and analyze a residual-based a posteriori
error estimator for the fully coupled scheme in two dimensions.
Residual estimators for mixed methods in Banach spaces have
been studied for coupled flow problems and
momentum-conservative Navier--Stokes formulations
\cite{cgos-M2AN-2021,ccov2022}.
For the present system, localization requires decompositions
adapted to test spaces with different field and divergence
exponents.
Using the bounded divergence lifting of
\cite{AcostaDuranMuschietti2006} and the Sobolev potential theorem
of \cite{Mardare2008}, we establish the decompositions needed
for the Darcy and Poisson test spaces, with and without
homogeneous normal trace, respectively.
A mixed-exponent decomposition handles the mechanical and ionic
flux tests \cite{ccov2022}.
Raviart--Thomas interpolation estimates \cite{ccov2022}
and Cl\'ement-type quasi-interpolation
\cite{clement1975,cgos-M2AN-2021}
yield computable element and edge contributions with the
corresponding dual exponents.
Reliability follows from continuous block stability and
residual localization under an additional smallness condition
on the prescribed data.
Bubble-function arguments and polynomial inverse and trace
inequalities \cite{Verfurth,ccov2022} establish local efficiency
for polynomial boundary data, with boundary-data oscillation
for general traces.
We also describe the estimator modifications for nonhomogeneous
displacement and Darcy normal-flux data under corresponding
well-posedness and uniform stability hypotheses.
Numerical experiments assess convergence of the mixed and
algebraically reconstructed fields, robustness in the nearly
incompressible regime, and estimator-driven adaptive refinement.
An additional example examines ionic polarization and the
mechanical and fluid response of a charged porous medium to
an imposed electric potential.

The remainder of this section introduces the notation and
preliminaries used throughout the paper.
Section~\ref{section2} presents the model and its first-order
pseudostress reformulation.
Section~\ref{section3} derives the variational formulations and
their functional setting.
The continuous stability and solvability analysis is given in
Section~\ref{section4}, followed by the Galerkin discretization
and a priori error estimates in Section~\ref{section5}.
Section~\ref{sec:posteriori-2D} develops the residual estimator,
its reliability and efficiency bounds, and the modifications
for nonhomogeneous boundary data.
Finally, Section~\ref{sec:numerical-results} presents the
numerical experiments.

\subsection*{Notation conventions and preliminaries.} 
Throughout, $\O\subset\R^d$, $d\in\{2,3\}$, is a bounded
connected Lipschitz domain with boundary $\Gamma$ and outward unit normal $\bnu$.
We use the standard Lebesgue and Sobolev spaces $\rL^t(\O)$ and
$\rW^{\ell,t}(\O)$, $1\le t<\infty$, $\ell\ge0$, with norms
$\|\cdot\|_{0,t;\O}$, $\|\cdot\|_{\ell,t;\O}$ and seminorm
$|\cdot|_{\ell,t;\O}$, respectively. Thus
$\rW^{0,t}(\O)=\rL^t(\O)$; for $t=2$ we write $\rH^\ell(\O)$
and omit $t$ from the norm notation.
Bold and blackboard symbols denote vector and tensor counterparts.
For conjugate exponents $t,t'\in(1,\infty)$, the trace space
$\rW^{1/t',t}(\Gamma)$ and its dual $\rW^{-1/t',t'}(\Gamma)$
are paired by $\langle\cdot,\cdot\rangle$; for $t=2$ these are
$\rH^{1/2}(\Gamma)$ and $\rH^{-1/2}(\Gamma)$.
We use $\|\cdot\|$ when the underlying space is clear,
$\I$ for the identity tensor, and $|\cdot|$ for the Euclidean norm.
For vector fields $\bv=(v_i)_{i=1}^d$ and $\bw=(w_i)_{i=1}^d$, set
\[
\nabla \boldsymbol{v} := \left(\frac{\partial v_i}{\partial x_j}\right)_{i,j=1}^d\,,\quad
\div(\boldsymbol{v}) := \sum_{j=1}^{d}\frac{\partial v_j}{\partial x_j}
\qan \bv\otimes\bw \,:=\, (v_iw_j)_{i,j=1}^d\,.
\]
For a tensor $\bta=(\tau_{ij})_{i,j=1}^d$, $\bdiv(\bta)$ acts
row-wise; its transpose, trace, inner product with $\bze$, and
deviatoric part are
\begin{equation}\label{transp-trace-inner-devia}
\bta^{\mathtt{t}}\,:=\,\left(\tau_{j i}\right)_{i,j=1}^d\,,\quad \operatorname{tr}(\bta)\,:=\, \sum_{i=1}^{d} \tau_{i i}
\qan \bta: \bze \,:=\, \sum_{i, j=1}^{d} \tau_{i j} \zeta_{i j}\,\qan \bta^{\tt d} = \bta-\dfrac{1}{d}\tr(\bta)\I\,.
\end{equation}
Furthermore, we recall the Hilbert space 
\begin{equation}\label{def-Hdiv}
\mathbb H(\bdiv;\Omega)\,:=\, \Big\{ \bta\in \mathbb L^2(\O): \quad \bdiv(\bta)\in \bL^2(\O)\Big\}\,
\end{equation}
with $\|\bta\|_{\bdiv;\O}:=\|\bta\|_{0;\O}
+\|\bdiv(\bta)\|_{0;\O}$.
We recall the tensorial case of \cite[Lemma 1.4]{g-SPRINGER-2014}, that is
\begin{equation}\label{int-by-parts-Hdiv-H1}
\langle \bta\bnu,\bv \rangle\,=\, \int_\O \Big\{ \bta:\nabla \bv \,\,+\,\, \bv \cdot\bdiv(\bta)\Big\}
\qquad \forall\,(\bta,\bv)\in \mathbb H(\bdiv;\O)\times \bH^1(\O)\,.
\end{equation}
For $1<t<\infty$, we also introduce the Banach spaces
\begin{equation*}%\label{def:graph-spaces}
\begin{array}{c}
\bH(\div_t;\O)
:=\bigl\{\bv\in\bL^2(\O):\quad
             \div(\bv)\in\rL^t(\O)\bigr\},
\qquad
\bH^t(\div;\O)
:=\bigl\{\bv\in\bL^t(\O):\quad
             \div(\bv)\in\rL^2(\O)\bigr\},
\\[1ex]
\bH^t(\div_t;\O)
:=\bigl\{\bv\in\bL^t(\O):\quad
\div(\bv)\in\rL^t(\O)\bigr\},
\end{array}
\end{equation*}
endowed with the respective norms
\begin{equation}\label{def:graph-norms}
\begin{array}{c}
\|\bv\|_{\div_t;\O}
:=\|\bv\|_{0;\O}+\|\div(\bv)\|_{0,t;\O},
\qquad
\|\bv\|_{t,\div;\O}
:=\|\bv\|_{0,t;\O}+\|\div(\bv)\|_{0;\O},
\\[1ex]
\|\bv\|_{t,\div_t;\O}
:=\|\bv\|_{0,t;\O}+\|\div(\bv)\|_{0,t;\O}.
\end{array}
\end{equation}
As in \cite[Section~1.3.4, eq.~(1.43)]{g-SPRINGER-2014}
and \cite[Section~3.1]{cgm-M2AN-2020}, for $t>1$ in two
dimensions and $t\ge6/5$ in three dimensions,
\begin{equation}\label{int-by-parts-Hdivp-H1}
\langle \bta\cdot\bnu,v\rangle\,=\, \int_\O \Big\{ \bta\cdot\nabla v \,\,+\,\, v \,\div(\bta)\Big\}
\qquad \forall\,(\bta,v)\in \mathbf H(\div_t;\O)\times \rH^1(\O)\,,
\end{equation}
where the boundary pairing is between $\rH^{-1/2}(\Gamma)$
and $\rH^{1/2}(\Gamma)$. For conjugate $t,t'\in(1,\infty)$,
\cite[Corollary~B.57]{eg-SPRINGER-2004} also gives
\begin{equation}\label{int-by-parts-Hpdivp-W1q}
\langle \bta\cdot\bnu,v \rangle\,=\, \int_\O \Big\{ \bta\cdot\nabla v \,\,+\,\, v \,\div(\bta)\Big\}
\qquad \forall\,(\bta,v)\in \mathbf H^t(\div_t;\O)\times \mathrm{W}^{1,t'}(\O)\,,
\end{equation}
where $\langle\cdot,\cdot\rangle$ denotes in \eqref{int-by-parts-Hpdivp-W1q} the duality pairing between 
$\mathrm{W}^{-1/t,t}(\Gamma)$ and $\mathrm{W}^{1/t,t'}(\Gamma)$.

\section{The model problem}\label{section2}

Following \cite[Section~2]{gir-CAMWA-2025}, we consider a homogeneous porous medium, consisting of incompressible grains and a charged interstitial fluid, occupying the domain $\O$. The mixture contains both positively and negatively charged ions. Given a body force $\ff$, a mass source $g$, a fixed-charge density $f_\chi$, and ionic source terms $f_i$, $i\in\{1,2\}$, and neglecting convective, gravitational, and inertial effects in the poromechanical equations, the stationary Biot--Poisson--Nernst--Planck system reads:
\begin{subequations}\label{eq:model}
	\begin{alignat}{2}
		\bsi_{\rm p} \,:=\, 2\mu\, \bep(\u) + \lambda\, \div(\u) \,\I - \alpha\, p \,\I  + \varepsilon\,\nabla \chi \otimes \nabla \chi  
		-\frac{\varepsilon}{2} |\nabla\chi|^2  \I - \delta (\xi_1-\xi_2)  \I &\qin\O \,,\label{eq:const} \\[0.5ex]
		- \bdiv (\bsi_{\rm p}) \,=\, \ff\,,\quad \z = -\dfrac{\kappa}{\nu}\nabla p\,,\quad c_0\, p \,+\, \alpha\, \div (\u) + \div (\z) \,=\, g &\qin \Omega\,,\label{eq:primal-strong} \\[0.5ex]
		\bvphi = \varepsilon \,\nabla \chi \,, \quad -\, \div(\bvphi)  = \xi_1 - \xi_2 + f_\chi &\qin \Omega\,,\\[0.5ex]
		\bt_i := \nabla \xi_i \,, \quad
		\bsi_i \,:=\, \kappa_i(\bt_i\,+\,\iota_i\varepsilon^{-1}\xi_i\,\bvphi)-\xi_i\,\z \,,\quad 
		\xi_i - \div(\bsi_i)-\xi_i\,\div(\z)\,=\, f_i  &\qin \Omega \,,\\[0.5ex]
		\u=\mathbf{0}\,,\quad  \z\cdot\bnu=0\,,\quad \chi=\chi_{\tt D} \qan \xi_i=\xi_{i,{\tt D}} &\qon \Gamma\,.
		% Boundary conditions are included in the model system.
	\end{alignat}
\end{subequations}
Here $\bsi_{\rm p}$ is the overall poroelastic (Cauchy) stress, $\u$ the solid displacement,
$\bep(\u):=\tfrac12(\nabla\u+\nabla\u^{\tt t})$ the infinitesimal
strain, $p$ the pore pressure, and $\z$ the Darcy velocity.
The dimensionless potential $\chi$ and electric displacement $\bvphi$
are coupled to positive and negative ion concentrations $\xi_1,\xi_2$,
their gradients $\bt_i$, and total ionic fluxes $\bsi_i$.
The parameters $\kappa,c_0,\alpha,\nu,\varepsilon,\delta$ denote
permeability, constrained specific storage, Biot--Willis coefficient,
fluid viscosity, permittivity, and osmotic coefficient, respectively;
$\kappa_i$ are ionic diffusivities, and $\iota_1=1$, $\iota_2=-1$ are the ionic valences.
The Lam\'e parameters satisfy $\mu>0$ and $\lambda\ge0$.
We assume $\kappa,\nu,\varepsilon,\kappa_1,\kappa_2>0$ and
$c_0,\alpha,\delta\ge0$.
The prescribed density $f_\chi$ accounts for immobile charges of the
solid matrix, so that the charge driving the electric field is
$\xi_1-\xi_2+f_\chi$.

To avoid imposing a symmetry constraint on the tensor unknown in \eqref{eq:model},
we follow the approach developed in \cite{gis-JSC-2022} and
reformulate the mechanical equations using the non-symmetric
pseudostress tensor $\bro$ defined by
\begin{equation}\label{def:brho}
	\bro \,:=\, \mu\,\nabla \u \,+\, (\lambda+\mu)\,\div(\u)\,\I \qin \Omega .
\end{equation}
The definition \eqref{def:brho} rewrites the elastic equilibrium
operator in \eqref{eq:primal-strong} without imposing stress symmetry.
In turn, as a direct consequence of \eqref{def:brho}, taking the trace gives
\begin{equation}\label{prop:divu:trrho}
	{\div(\u) \,=\, c_\lambda\,\tr(\bro), \qquad \text{where } c_\lambda:=\bigl(d\lambda+(d+1)\mu\bigr)^{-1}.}
\end{equation}
The subscript emphasizes the nearly incompressible limit
$\lambda\to+\infty$. Substitution of \eqref{prop:divu:trrho}
into \eqref{def:brho} gives the constitutive equation
\eqref{Biot-NPN-mixed-1} below, using the deviatoric operator
\eqref{transp-trace-inner-devia}.
In addition, we notice from \eqref{prop:divu:trrho}, the divergence theorem and the {Dirichlet boundary condition} satisfied by $\u$ that 
\begin{equation}\label{prop:tr:mean:zero}
c_\lambda\int_\Omega \tr(\bro) = \int_\Omega \div(\u) = \int_\Gamma \u\cdot\bnu = 0\,.
\end{equation}
Finally, integrating the mass equation in \eqref{eq:primal-strong} and using
$\z\cdot\bnu=0$ and \eqref{prop:tr:mean:zero} gives
\begin{equation}\label{eq:pressure-mean}
	c_0\int_\O p \,=\, \int_\O g\,.
\end{equation}
If $c_0>0$, the mean of $p$ is fixed by the data and $g\in\rL^2(\O)$ is arbitrary;
if $c_0=0$, \eqref{eq:pressure-mean} is the compatibility condition $\int_\O g=0$ and
we normalize
\begin{equation}\label{uniq-p-cond}
	\int_\O p \,=\, 0 .
\end{equation}
In both cases we write $p=p_0+\overline p$, with
\begin{equation}\label{def:p-mean}
	\int_\O p_0=0\,,
	\qquad
	\overline p:=\overline g/c_0\ \text{ if } c_0>0\,,
	\quad
	\overline p:=0\ \text{ if } c_0=0\,,
\end{equation}
where $\overline g:=|\O|^{-1}\int_\O g$. Since $\nabla p=\nabla p_0$ and
$c_0\overline p=\overline g$, $p_0$ solves the same system with $g-\overline g$ in place
of $g$. We use $p_0$ as the unknown (cf. Lemma~\ref{lem:pressure-mean}); the constant
$\overline p$ enters only the Cauchy stress reconstruction and the estimator.

Collecting these identities gives the first-order system, for
$i\in\{1,2\}$:
\begin{subequations}\label{eq:strong-complete}
	\begin{alignat}{2}
		\dfrac{1}{\mu}\,\bro^{\tt d} \,+\, \dfrac{c_\lambda}{ d}\,\tr(\bro)\,\I \,-\, \nabla \u &\,=\, 0 && \quad\qin \Omega\,, \label{Biot-NPN-mixed-1}\\[0.5ex]
		-\,\bdiv(\bro) \,-\, \dfrac{\alpha\,\nu}{\kappa}\,\z \,+\, \varepsilon^{-1}(\xi_1-\xi_2+f_\chi)\,\bvphi \,+\, \delta\,(\bt_1-\bt_2) &\,=\, \ff  && \quad\qin \Omega\,, \label{Biot-NPN-mixed-2}\\[0.5ex]
		\dfrac{\nu}{\kappa}\,\z \,+\, \nabla p & \,=\, \mathbf{0} && \quad\qin \Omega\,, \label{Biot-NPN-mixed-3}\\[0.5ex]
		c_0\,p \,+\, \alpha\,c_\lambda\,\tr(\bro) \,+\, \div(\z) &\,=\, g && \quad\qin \Omega\,,  \label{Biot-NPN-mixed-4} \\[0.5ex]
		\varepsilon^{-1}\,\bvphi \,-\, \nabla\chi & \,=\, \mathbf{0} && \quad\qin \Omega\,,\label{Biot-NPN-mixed-5-1} \\[0.5ex]
		-\,\div(\bvphi) & \,=\, \xi_1-\xi_2+f_\chi && \quad\qin \Omega\,, \label{Biot-NPN-mixed-5-2}\\[0.5ex]
		\bt_i \,-\, \nabla\xi_i & \,=\, \mathbf{0} && \quad\qin \Omega\,, \label{Biot-NPN-mixed-6}\\[0.5ex]
		-\,\bsi_i \,+\, \kappa_i\,\bt_i \,+\, \iota_i\,\kappa_i\,\varepsilon^{-1}\,\xi_i\,\bvphi \,-\, \xi_i\,\z & \,=\, \mathbf{0} && \quad\qin \Omega\,,\label{Biot-NPN-mixed-7}\\[0.5ex]
		\xi_i \,-\, \div(\bsi_i) \,-\, \xi_i\,\div(\z) & \,=\, f_i && \quad\qin \Omega\,. \label{Biot-NPN-mixed-8}\\[0.5ex]
		\u=\mathbf{0}\,,\quad  \z\cdot\bnu=0\,,\quad \chi=\chi_{\tt D} \qan \xi_i & =\xi_{i,{\tt D}} &&\quad\qon \Gamma\,.\label{Biot-NPN-mixed-bc}
	\end{alignat}
\end{subequations}
The first-order system \eqref{eq:strong-complete} also gives the
gradients directly in terms of its unknowns:
\begin{equation}\label{eq:gradient-reconstruction}
\begin{aligned}
\nabla\u
&=\mu^{-1}\bro^{\tt d}
  +\frac{c_\lambda}{d}\tr(\bro)\I,
&\qquad \nabla p&=-\frac{\nu}{\kappa}\z,
&\qquad
\nabla\chi&=\varepsilon^{-1}\bvphi,
&\qquad \nabla\xi_i&=\bt_i,\quad i\in\{1,2\}.
\end{aligned}
\end{equation}
Using these identities (cf. \eqref{eq:gradient-reconstruction})
in \eqref{eq:const}, the Cauchy stress is recovered as
\begin{equation}\label{eq:Cauchy-reconstruction}
	\bsi_{\rm p} \,=\, \bro + \bro^{\tt t} \,-\, \dfrac{\lambda + 2\mu}{{d\lambda + (d+1)\mu}}\,\tr(\bro)\,\I \,-\, \alpha\,p\,\I \,+\, \varepsilon^{-1}\!\left(\bvphi\otimes\bvphi 
	\,-\, \dfrac{1}{2}\, |\bvphi|^2\,\I\right) \,-\, \delta\,(\xi_1-\xi_2)\,\I .
\end{equation}
Here $p=p_0+\overline p$ is the full pressure of \eqref{def:p-mean}.

\section{The fully mixed formulation}\label{section3}

We derive a fully mixed formulation of \eqref{eq:strong-complete}
in Banach spaces, treating the pseudostress and fluxes as independent
unknowns and recovering gradients through \eqref{eq:gradient-reconstruction}.
The new pseudostress--displacement block is coupled to the Darcy,
electrostatic, and ionic formulations of
\cite[Sections~3.2--3.4]{gir-CAMWA-2025}.

\subsection{Preliminaries}\label{section31}

As in \cite[Section~3.1]{gir-CAMWA-2025}, the coupled products
determine the Lebesgue exponents: $(\xi_1-\xi_2+f_\chi)\bvphi$ in
\eqref{Biot-NPN-mixed-2} and $\xi_i\bvphi$, $\xi_i\z$ in
\eqref{Biot-NPN-mixed-7} are paired with $\bL^2$ tests, whereas
$\xi_i\div(\z)$ in \eqref{Biot-NPN-mixed-8} is paired with a
concentration test.
For $r>2$, set
\begin{equation}\label{def-r-s-rho-varrho}
s:=\frac{r}{r-1},
\qquad
\rho:=\frac{2r}{r-2},
\qquad
\varrho:=\frac{\rho}{\rho-1}.
\end{equation}
The pairs $(r,s)$ and $(\rho,\varrho)$ are conjugate, with
\begin{equation}\label{relation-r-rho}
\frac1r+\frac1\rho=\frac12.
\end{equation}
Throughout, the plain symbol $\rho$ denotes this Lebesgue exponent, whereas the boldface $\bro$ denotes the pseudostress.
Thus, for $\xi_i\in\rL^\rho(\O)$ and
$\bvphi,\z\in\bL^r(\O)$, H\"older's inequality gives
\[
\begin{aligned}
\|(\xi_1-\xi_2)\bvphi\|_{0;\O}
&\le
\|\xi_1-\xi_2\|_{0,\rho;\O}\|\bvphi\|_{0,r;\O},
\\[1ex]
\|\xi_i\bvphi\|_{0;\O}+\|\xi_i\z\|_{0;\O}
&\le
\|\xi_i\|_{0,\rho;\O}
\bigl(\|\bvphi\|_{0,r;\O}+\|\z\|_{0,r;\O}\bigr).
\end{aligned}
\]
Likewise, for the fixed charge we assume $f_\chi\in\rL^\rho(\O)$ and set
$\Lambda_\chi:=\|f_\chi\|_{0,\rho;\O}$, so that
$\|f_\chi\bvphi\|_{0;\O}\le\Lambda_\chi\|\bvphi\|_{0,r;\O}$.
Testing \eqref{Biot-NPN-mixed-4} in $\rL^2$ requires
$\div(\z)\in\rL^2(\O)$. Thus concentration tests
$\eta_i\in\rL^\rho(\O)$ require $1/2+2/\rho\le1$, i.e.
$\rho\ge4$, for convection in \eqref{Biot-NPN-mixed-8}:
\begin{equation}\label{eq:NP-convection-bound}
\left|\int_\O\xi_i\,\div(\z)\,\eta_i\right|
\le
|\O|^{1/2-2/\rho}
\|\xi_i\|_{0,\rho;\O}
\|\div(\z)\|_{0;\O}
\|\eta_i\|_{0,\rho;\O}.
\end{equation}
By \eqref{def-r-s-rho-varrho}, $\rho\ge4$ is
equivalent to $r\le4$ and also gives $\rho\ge r$.
Consequently, $\rL^\rho(\O)\hookrightarrow\rL^r(\O)$,
as required by the charge equation
\eqref{Biot-NPN-mixed-5-2}.

The ionic boundary-load estimates also use
$\rH^1(\O)\hookrightarrow\rL^\rho(\O)$.
This embedding holds for every finite $\rho$ in two
dimensions and for $\rho\le6$ in three dimensions.
Combining these requirements with the exponent ranges
used in the auxiliary Poisson estimates of
\cite[Section~3.1]{gir-CAMWA-2025}, we retain
\begin{equation}\label{r-rho-range}
\left\{
\begin{array}{l@{\quad}l@{\quad}l@{\quad}l@{\qquad}l}
2<r\le4,
&
s=\dfrac{r}{r-1},
&
\rho=\dfrac{2r}{r-2},
&
\varrho=\dfrac{\rho}{\rho-1},
&
d=2,
\\[2ex]
r=3,
&
s=\dfrac32,
&
\rho=6,
&
\varrho=\dfrac65,
&
d=3.
\end{array}
\right.
\end{equation}
We fix one such choice throughout the paper.
In particular, the coupled fields satisfy
\begin{equation}\label{xi-z-spaces}
\xi_i\in\rL^\rho(\O),
\qquad
\bvphi\in\bL^r(\O),
\qquad
\z\in\bH^r(\div;\O),
\qquad
\bt_i\in\bL^2(\O),
\quad i\in\{1,2\}.
\end{equation}
We incorporate the homogeneous normal-trace condition
for the Darcy flux and its tests through
\begin{align}
\bH_0^s(\div_s;\O)
&:=
\Big\lbrace
\bw\in\bH^s(\div_s;\O):\quad
\bw\cdot\bnu=0\qon\Gamma
\Big\rbrace,
\label{H0rdiv}
\\[1ex]
\bH_0^r(\div;\O)
&:=
\Big\lbrace
\bw\in\bH^r(\div;\O):\quad
\bw\cdot\bnu=0\qon\Gamma
\Big\rbrace.
\label{H0sdivs}
\end{align}
The boundary equalities in
\eqref{H0rdiv}--\eqref{H0sdivs} are understood
in the normal-trace sense.
For $1\le t<\infty$, we also set
\begin{equation}\label{def-L0t}
\rL_0^t(\O)
:=
\left\{
q\in\rL^t(\O):\quad \int_\O q=0
\right\}.
\end{equation}
The pseudostress normalization
\eqref{prop:tr:mean:zero} leads to
\begin{equation}\label{H0-decomp}
\begin{gathered}
\bbH(\bdiv;\O)=\bbH_0(\bdiv;\O)\oplus\R\,\I,
\\[1ex]
\bbH_0(\bdiv;\O)
:=
\left\{
\bta\in\bbH(\bdiv;\O):\quad
\int_\O\tr(\bta)=0
\right\}.
\end{gathered}
\end{equation}
Here, the subscript $0$ denotes a zero-mean trace
condition for tensors, whereas it denotes a vanishing
normal trace in \eqref{H0rdiv}--\eqref{H0sdivs}.

\subsection{The poroelasticity equations}\label{section32}

Given $\ff\in\bL^2(\O)$, $g\in\rL^2(\O)$ (with $\int_\O g=0$ if $c_0=0$), and
prescribed electrochemical fields $\bvphi$, $\xi_i$, and
$\bt_i$, $i\in\{1,2\}$, satisfying \eqref{xi-z-spaces},
we derive the mixed formulation of
\eqref{Biot-NPN-mixed-1}--\eqref{Biot-NPN-mixed-4}.
For sufficiently regular fields, testing the constitutive
equation \eqref{Biot-NPN-mixed-1} against
$\bta\in\bbH(\bdiv;\O)$ and integrating by parts,
using $\u=\mathbf0$ on $\Gamma$
(cf. \eqref{int-by-parts-Hdiv-H1}), gives
\begin{equation}\label{eq:Biot-constitutive-identity}
\frac1\mu\int_\O\bro^{\tt d}:\bta^{\tt d}
+\frac{c_\lambda}{d}\int_\O\tr(\bro)\tr(\bta)
+\int_\O\u\cdot\bdiv(\bta)=0.
\end{equation}
Testing the equilibrium equation
\eqref{Biot-NPN-mixed-2} against $\bv\in\bL^2(\O)$
and rearranging yields
\begin{equation}\label{eq:Biot-equilibrium-identity}
\begin{aligned}
\int_\O\bv\cdot\bdiv(\bro)
+\frac{\alpha\nu}{\kappa}\int_\O\z\cdot\bv
={}&
\varepsilon^{-1}\int_\O(\xi_1-\xi_2+f_\chi)\bvphi\cdot\bv
+\delta\int_\O(\bt_1-\bt_2)\cdot\bv
-\int_\O\ff\cdot\bv.
\end{aligned}
\end{equation}
The right-hand side of \eqref{eq:Biot-equilibrium-identity}
defines a bounded functional on $\bL^2(\O)$ by the
product estimates of Section~\ref{section31}
(cf. \eqref{relation-r-rho} and \eqref{xi-z-spaces}).
Thus, both mechanical identities are well defined in the tensor
space \eqref{def-Hdiv}, with
$\bro,\bta\in\bbH(\bdiv;\O)$ and
$\u,\bv\in\bL^2(\O)$.

Next, proceeding as in \cite[Section~3.2]{gir-CAMWA-2025},
we test the Darcy and mass equations
\eqref{Biot-NPN-mixed-3}--\eqref{Biot-NPN-mixed-4}
against $\bw\in\bH_0^s(\div_s;\O)$ and
$q\in\rL_0^2(\O)$, respectively.
From here on, $p$ denotes $p_0$ in \eqref{def:p-mean}; the constant
$\overline p$ drops out of both tests, and the full pressure is $p+\overline p$.
Integrating the pressure gradient by parts, using
$\bw\cdot\bnu=0$ on $\Gamma$
(cf. \eqref{int-by-parts-Hpdivp-W1q}),
and multiplying the mass identity by $-1$, we obtain
\begin{subequations}\label{eq:Biot-fluid-identities}
\begin{align}
\frac{\nu}{\kappa}\int_\O\z\cdot\bw
-\int_\O p\,\div(\bw)
&=0,
\label{eq:Biot-Darcy-identity}
\\[1ex]
-\int_\O q\,\div(\z)
-c_0\int_\O p\,q
-\alpha c_\lambda\int_\O\tr(\bro)\,q
&=-\int_\O g\,q.
\label{eq:Biot-mass-identity}
\end{align}
\end{subequations}
By conjugacy of $r$ and $s$, the Darcy identity
\eqref{eq:Biot-Darcy-identity} requires only
$p\in\rL^r(\O)$ and $\z\in\bL^r(\O)$.
The mass identity \eqref{eq:Biot-mass-identity} requires
$\div(\z)\in\rL^2(\O)$, while its storage term, namely the second term on the
left-hand side of \eqref{eq:Biot-mass-identity}, is finite
by the Cauchy--Schwarz inequality, since
$p\in\rL^r(\O)\hookrightarrow\rL^2(\O)$ and
$q\in\rL^2(\O)$.
These requirements agree with the Darcy flux space
in \eqref{xi-z-spaces}.

The normalization \eqref{prop:tr:mean:zero} imposes zero mean trace
on $\bro$, and $p=p_0$ has zero mean by \eqref{def:p-mean}. We impose the same
constraints on the tensor and mass tests. No equation is lost:
\eqref{H0-decomp} omits only multiples of $\I$, against which
\eqref{eq:Biot-constitutive-identity} holds; and
\eqref{eq:Biot-mass-identity} tested with the constant $q\equiv1$ reduces to
\eqref{eq:pressure-mean}, which holds by the definition of $\overline p$,
because $\tr(\bro)$ has zero mean and $\z$ has zero normal trace.
Using the normal-trace spaces \eqref{H0rdiv}--\eqref{H0sdivs},
the zero-mean spaces \eqref{def-L0t}, and the decomposition
\eqref{H0-decomp}, we therefore define
\begin{equation}\label{def-bXi-bQi}
\begin{aligned}
\bX&:=\bbH_0(\bdiv;\O),
&\qquad
\bX_2&:=\bH_0^r(\div;\O),
&\qquad
\bX_1&:=\bH_0^s(\div_s;\O),
\\[1ex]
\bQ&:=\bL^2(\O),
&\qquad
\bQ_2&:=\rL_0^2(\O),
&\qquad
\bQ_1&:=\rL_0^r(\O),
\end{aligned}
\end{equation}
endowed with the corresponding norms
(cf. \eqref{def:graph-norms}).
All product spaces carry the sum norm.

To express the mechanical identities
(cf. \eqref{eq:Biot-constitutive-identity} and
\eqref{eq:Biot-equilibrium-identity})
in terms of bilinear forms and a linear functional, define
$\bfa_s:\bX\times\bX\to\R$,
$\bfb_s:\bX\times\bQ\to\R$,
$\bfe_s:\bX_2\times\bQ\to\R$, and
$\bF_{\bvphi,\bxi,\bt}\in\bQ'$ by
\begin{align}
\bfa_s(\bro,\bta)
&:=
\frac1\mu\int_\O\bro^{\tt d}:\bta^{\tt d}
+\frac{c_\lambda}{d}\int_\O\tr(\bro)\tr(\bta),
\qquad
\bfb_s(\bta,\bv):=\int_\O\bv\cdot\bdiv(\bta),
\nonumber\\[1ex]
\bfe_s(\z,\bv)
&:=\frac{\alpha\nu}{\kappa}\int_\O\z\cdot\bv,
\qquad
\bF_{\bvphi,\bxi,\bt}(\bv)
:=
\varepsilon^{-1}\int_\O(\xi_1-\xi_2+f_\chi)\bvphi\cdot\bv
+\delta\int_\O(\bt_1-\bt_2)\cdot\bv
-\int_\O\ff\cdot\bv,
\label{def-as-bs-es-bF}
\end{align}
where $\bxi:=(\xi_1,\xi_2)$ and
$\bt:=(\bt_1,\bt_2)$.
The form $\bfe_s$ contains the coupling with the Darcy
flux, while $\bF_{\bvphi,\bxi,\bt}$ collects the
electrochemical contributions and the body force.

Similarly, the fluid identities \eqref{eq:Biot-fluid-identities}
define
$\bfa_f:\bX_2\times\bX_1\to\R$,
$\bfd_i:\bX_i\times\bQ_i\to\R$, $i\in\{1,2\}$,
$\bfe_f:(\bX\times\bQ_1)\times\bQ_2\to\R$, and
$\bG\in\bQ_2'$ by
\begin{align}
\bfa_f(\z,\bw)
&:=\frac{\nu}{\kappa}\int_\O\z\cdot\bw,
\qquad
\bfd_i(\bw,q):=-\int_\O q\,\div(\bw),
\quad i\in\{1,2\},
\nonumber\\[1ex]
\bfe_f((\bro,p),q)
&:=c_0\int_\O p\,q
+\alpha c_\lambda\int_\O\tr(\bro)\,q,
\qquad
\bG(q):=-\int_\O g\,q.
\label{def-af-d1-d2-ef}
\end{align}
The forms $\bfa_f$ and $\bfd_i$ coincide with those in
\cite[Section~3.2]{gir-CAMWA-2025}, whereas $\bfe_f$
uses the pseudostress identity \eqref{prop:divu:trrho}.

With the forms \eqref{def-as-bs-es-bF} and \eqref{def-af-d1-d2-ef},
the poroelastic problem reads: find $(\bro,\u)\in\bX\times\bQ$
and $(\z,p)\in\bX_2\times\bQ_1$ such that
\begin{equation}\label{Biot-fully-mixed}
\begin{array}{lllcll}
\bfa_s(\bro,\bta)+\bfb_s(\bta,\u)
&&
&=&0
&\quad\forall\,\bta\in\bX,
\\[1ex]
\bfb_s(\bro,\bv)
&&+\bfe_s(\z,\bv)
&=&\bF_{\bvphi,\bxi,\bt}(\bv)
&\quad\forall\,\bv\in\bQ,
\\[1ex]
&\bfa_f(\z,\bw)+\bfd_1(\bw,p)
&
&=&0
&\quad\forall\,\bw\in\bX_1,
\\[1ex]
&\bfd_2(\z,q)
&-\bfe_f((\bro,p),q)
&=&\bG(q)
&\quad\forall\,q\in\bQ_2.
\end{array}
\end{equation}
The arrangement in \eqref{Biot-fully-mixed} separates the mechanical
and generalized Darcy saddle-point blocks from the perturbations
$\bfe_s$ and $\bfe_f$.

Finally, the Cauchy--Schwarz and H\"older inequalities give
\begin{equation}\label{Biot-bounds}
\begin{aligned}
\|\bfa_s\|&\le\frac1\mu,
&\qquad
\|\bfb_s\|&\le1,
&\qquad 
\|\bfe_s\|&\le C_r(\O)\frac{\alpha\nu}{\kappa},
\\[1ex]
\|\bfa_f\|&\le\frac{\nu}{\kappa},
&\qquad
\|\bfd_i\|&\le1,\quad i\in\{1,2\},
&\qquad
\|\bfe_f\|&\le
\max\{C_r(\O)c_0,\sqrt d\,\alpha c_\lambda\},
\\[1ex]
\|\bG\|&\le\|g\|_{0;\O},
\end{aligned}
\end{equation}
where $C_r(\O):=|\O|^{1/2-1/r}$. By \eqref{relation-r-rho}, the mechanical load satisfies
\begin{equation}\label{Biot-F-bounds}
\begin{aligned}
|\bF_{\bvphi,\bxi,\bt}(\bv)|
\le{}&
\Bigl\{
\varepsilon^{-1}\bigl(\|\xi_1-\xi_2\|_{0,\rho;\O}+\Lambda_\chi\bigr)
\|\bvphi\|_{0,r;\O}
+\delta\|\bt_1-\bt_2\|_{0;\O}
+\|\ff\|_{0;\O}
\Bigr\}\|\bv\|_{\bQ}
\qquad\forall\,\bv\in\bQ.
\end{aligned}
\end{equation}

\begin{lemma}[Natural pressure mean]\label{lem:pressure-mean}
Let $c_0>0$ and $g\in\rL^2(\O)$, and consider the problem obtained from
\eqref{Biot-fully-mixed} by replacing $\bQ_1$ and $\bQ_2$ with the full
spaces $\rL^r(\O)$ and $\rL^2(\O)$. If $(\bro,\u,\z,\widetilde p)$ solves it, then
$\widetilde p-\overline p$ has zero mean and $(\bro,\u,\z,\widetilde p-\overline p)$
solves \eqref{Biot-fully-mixed}. Conversely, if $(\bro,\u,\z,p)$ solves
\eqref{Biot-fully-mixed}, then $(\bro,\u,\z,p+\overline p)$ solves the
full-space problem. The same holds on the finite element spaces of
Section~\ref{section51}, with $\rP_\ell(\Th)$ in place of $\rL^r(\O)$ and
$\rL^2(\O)$. Moreover,
$\|\overline p\|_{0,r;\O}\le\|g\|_{0;\O}/(C_r(\O)\,c_0)$.
\end{lemma}
\begin{proof}
Testing the mass equation with $q\equiv1$, and using $\int_\O\div(\z)=0$ and
$\int_\O\tr(\bro)=0$, gives $c_0\int_\O\widetilde p=\int_\O g$, so
$\widetilde p-\overline p$ has zero mean. The Darcy equation is unaffected by
constants since $\int_\O\div(\bw)=0$, and so is the mass equation tested with
zero-mean $q$; the constant part of the latter holds by the choice of
$\overline p$. The bound follows from
$\|\overline p\|_{0,r;\O}=|\O|^{1/r}|\overline g|/c_0
\le|\O|^{1/r-1/2}\|g\|_{0;\O}/c_0$.
\end{proof}
Hence Sections~\ref{section4}--\ref{sec:posteriori-2D} treat the zero-mean unknown $p$ with no
restriction on the mean of $g$ when $c_0>0$. The errors satisfy
$(p+\overline p)-(p_h+\overline p)=p-p_h$, and the full pressure obeys
$\|p+\overline p\|_{0,r;\O}\le\|p\|_{0,r;\O}+\|g\|_{0;\O}/(C_r(\O)c_0)$, where the
dependence on $1/c_0$ is unavoidable since the mean of $p$ is $\overline g/c_0$.

\subsection{Mixed formulation of the Poisson equations}
\label{section33}

The Poisson equations
\eqref{Biot-NPN-mixed-5-1}--\eqref{Biot-NPN-mixed-5-2}
retain the formulation of
\cite[Section~3.3]{gir-CAMWA-2025}, with the fixed charge $f_\chi$ added to the
load $F_{\bxi}$ below.
With the exponents chosen in Section~\ref{section31}
(cf. \eqref{r-rho-range}), its trial and test spaces are
\begin{equation}\label{def:Poisson-spaces}
\begin{aligned}
\rX_2&:=\bH^r(\div_r;\O),
&\qquad \rM_1&:=\rL^r(\O),
&\qquad 
\rX_1&:=\bH^s(\div_s;\O),
&\qquad \rM_2&:=\rL^s(\O),
\end{aligned}
\end{equation}
endowed with their norms
(cf. \eqref{def:graph-norms}).
Given $\chi_{\tt D}\in\rW^{1/s,r}(\Gamma)$
and prescribed
$\bxi\in\rL^\rho(\O)\times\rL^\rho(\O)$,
and the fixed charge $f_\chi\in\rL^\rho(\O)$, define $a:\rX_2\times\rX_1\to\R$,
$b_i:\rX_i\times\rM_i\to\R$, $i\in\{1,2\}$,
$G\in\rX_1'$, and $F_{\bxi}\in\rM_2'$ by
\begin{align}
a(\bvphi,\bpsi)
&:=\varepsilon^{-1}\int_\O\bvphi\cdot\bpsi,
\qquad
b_i(\bpsi,\gamma):=\int_\O\gamma\,\div(\bpsi),
\quad
(\bpsi,\gamma)\in\rX_i\times\rM_i,\quad i\in\{1,2\},
\nonumber\\[1ex]
G(\bpsi)
&:=\langle\bpsi\cdot\bnu,\chi_{\tt D}\rangle_\Gamma,
\qquad
F_{\bxi}(\gamma):=-\int_\O(\xi_1-\xi_2+f_\chi)\gamma.
\label{def:Poisson-forms}
\end{align}
With the forms \eqref{def:Poisson-forms}, the Poisson problem reads:
find $(\bvphi,\chi)\in\rX_2\times\rM_1$ such that
\begin{equation}\label{Poisson-mixed-formulation}
\begin{array}{lllcll}
a(\bvphi,\bpsi)
&+&b_1(\bpsi,\chi)
&=&G(\bpsi)
&\quad\forall\,\bpsi\in\rX_1,
\\[1ex]
b_2(\bvphi,\gamma)
&&
&=&F_{\bxi}(\gamma)
&\quad\forall\,\gamma\in\rM_2.
\end{array}
\end{equation}
The continuity estimates of
\cite[Section~3.3]{gir-CAMWA-2025} give
$\|a\|\le\varepsilon^{-1}$ and $\|b_i\|\le1$,
$i\in\{1,2\}$, together with
\begin{align}
\|F_{\bxi}\|_{\rM_2'}
&\le\|F\|\,\bigl(\|\xi_1-\xi_2\|_{0,\rho;\O}+\Lambda_\chi\bigr),
\qquad
\|F\|:=|\O|^{1/r-1/\rho},
\label{bound-F-xi}\\[1ex]
\|G\|_{\rX_1'}
&\le c_r\|\chi_{\tt D}\|_{1/s,r;\Gamma},
\label{bound-G}
\end{align}
where $c_r>0$ depends only on $\O$ and $r$.
The charge-load estimate \eqref{bound-F-xi}
uses $\rL^\rho(\O)\hookrightarrow\rL^r(\O)$,
and \eqref{bound-G} uses the normal-trace estimate.

\subsection{Mixed formulation of the Nernst--Planck equations}
\label{section34}

With the exponents chosen in Section~\ref{section31}
(cf. \eqref{r-rho-range}), we introduce the following spaces
for the Nernst--Planck equations
\eqref{Biot-NPN-mixed-6}--\eqref{Biot-NPN-mixed-8}
(cf. \cite[Section~3.4]{gir-CAMWA-2025}):
\begin{equation}\label{def:bH-H1-H2-M}
\calH_1:=\bL^2(\O),
\qquad
\calH_2:=\bH(\div_\varrho;\O),
\qquad
\bcalH:=\calH_1\times\calH_2,
\qquad
\calM:=\rL^\rho(\O).
\end{equation}
These spaces carry their respective norms, with the sum norm
on $\bcalH$.
For each $i\in\{1,2\}$, we write
$\vec{\bt}_i:=(\bt_i,\bsi_i)\in\bcalH$ and
$\vec{\bs}_i:=(\bs_i,\bta_i)\in\bcalH$.

For prescribed $(\z,\bvphi)\in\bX_2\times\rX_2$,
define
$\calA:\bcalH\times\bcalH\to\R$,
$\calB:\bcalH\times\calM\to\R$,
$\calC,\calD_{\z}:\calM\times\calM\to\R$,
and $\calE_{\z,\bvphi}:\bcalH\times\calM\to\R$ by
\begin{align}
\calA(\vec{\bt}_i,\vec{\bs}_i)
&:=\kappa_i\int_\O\bt_i\cdot\bs_i
-\int_\O\bsi_i\cdot\bs_i
+\int_\O\bta_i\cdot\bt_i,
\qquad
\calB(\vec{\bs}_i,\eta_i)
:=\int_\O\eta_i\,\div(\bta_i),
\nonumber\\[1ex]
\calC(\xi_i,\eta_i)
&:=\int_\O\xi_i\eta_i,
\qquad
\calD_{\z}(\xi_i,\eta_i)
:=\int_\O\xi_i\,\div(\z)\,\eta_i,
\nonumber\\[1ex]
\calE_{\z,\bvphi}(\vec{\bs}_i,\eta_i)
&:=-\int_\O\eta_i\z\cdot\bs_i
+\iota_i\kappa_i\varepsilon^{-1}
\int_\O\eta_i\bvphi\cdot\bs_i.
\label{def-A-B-C-E-D}
\end{align}
Given $\xi_{i,\tt D}\in\rH^{1/2}(\Gamma)$ and
$f_i\in\rL^\varrho(\O)$, define
$\mathcal G\in\bcalH'$ and $\mathcal F\in\calM'$ by
\begin{equation}\label{def:NP-load-functionals}
\mathcal G(\vec{\bs}_i)
:=\langle\bta_i\cdot\bnu,\xi_{i,\tt D}\rangle_\Gamma,
\qquad
\mathcal F(\eta_i):=-\int_\O f_i\eta_i.
\end{equation}
For each $i\in\{1,2\}$, the ionic problem with the forms
\eqref{def-A-B-C-E-D} and loads \eqref{def:NP-load-functionals}
reads: find $(\vec{\bt}_i,\xi_i)\in\bcalH\times\calM$ such that
\begin{equation}\label{NP-mixed-formulation}
\begin{array}{lllcll}
\calA(\vec{\bt}_i,\vec{\bs}_i)
&+\calB(\vec{\bs}_i,\xi_i)
&+\calE_{\z,\bvphi}(\vec{\bs}_i,\xi_i)
&=&\mathcal G(\vec{\bs}_i)
&\quad\forall\,\vec{\bs}_i\in\bcalH,
\\[1ex]
\calB(\vec{\bt}_i,\eta_i)
&-\calC(\xi_i,\eta_i)
&+\calD_{\z}(\xi_i,\eta_i)
&=&\mathcal F(\eta_i)
&\quad\forall\,\eta_i\in\calM.
\end{array}
\end{equation}
The terms $\calE_{\z,\bvphi}$ and $\calD_{\z}$ in
\eqref{NP-mixed-formulation} perturb the ionic twofold
saddle-point structure.

The continuity estimates in
\cite[Section~3.4]{gir-CAMWA-2025}, with the
domain factors written explicitly, give
\begin{equation}\label{bounds-NP-A-B-C-G-F}
\begin{aligned}
\|\calA\|&\le\max\{\kappa_i,1\},
&\qquad
\|\calB\|&\le1,
&\qquad
\|\calC\|&\le|\O|^{1-2/\rho},
\\[1ex]
\|\mathcal G\|&\le C\|\xi_{i,\tt D}\|_{1/2;\Gamma},
&\qquad
\|\mathcal F\|&\le\|f_i\|_{0,\varrho;\O},
\end{aligned}
\end{equation}
where $C>0$ depends only on $\O$ and $\rho$.
The boundary-load estimate in
\eqref{bounds-NP-A-B-C-G-F} uses the integration-by-parts identity
\eqref{int-by-parts-Hdivp-H1}, an
$\rH^1(\O)$ lifting and
$\rH^1(\O)\hookrightarrow\rL^\rho(\O)$.
The coupling estimates needed in the error analysis are
\begin{equation}\label{bound-E-D}
\begin{aligned}
|\calD_{\bw}(\xi_i,\eta_i)|
&\le\|\calD\|\,\|\bw\|_{\bX_2}
\|\xi_i\|_{\calM}\|\eta_i\|_{\calM},
\\[1ex]
|\calE_{\bw,\bphi}(\vec{\bs}_i,\eta_i)|
&\le\|\calE\|\,\|(\bw,\bphi)\|_{\bX_2\times\rX_2}
\|\vec{\bs}_i\|_{\bcalH}\|\eta_i\|_{\calM},
\end{aligned}
\end{equation}
for $(\bw,\bphi)\in\bX_2\times\rX_2$,
$\xi_i,\eta_i\in\calM$, and
$\vec{\bs}_i\in\bcalH$, with
\begin{equation}\label{bounds-NP-D-E}
\|\calD\|:=|\O|^{1/2-2/\rho},
\qquad
\|\calE\|:=\max\{\kappa_i\varepsilon^{-1},1\}.
\end{equation}
The constants in \eqref{bounds-NP-D-E} follow from
\eqref{eq:NP-convection-bound} and the product estimates
associated with \eqref{relation-r-rho}.

\subsection{The fully mixed formulation of the coupled problem}
\label{section35}

For the two species, set $\bcalM:=\calM\times\calM$ and
$\bcalH_1:=\calH_1\times\calH_1$, with sum norms
\[
\|\bxi\|_{\bcalM}:=\sum_{i=1}^2\|\xi_i\|_{\calM},
\qquad
\|\bt\|_{\bcalH_1}:=\sum_{i=1}^2\|\bt_i\|_{\calH_1}.
\]
We now combine the poroelastic, Poisson, and ionic formulations
developed in Sections~\ref{section32}, \ref{section33},
and~\ref{section34}, respectively.
Assume that $\ff\in\bL^2(\O)$,
$g\in\rL^2(\O)$ (with $\int_\O g=0$ if $c_0=0$),
$f_\chi\in\rL^\rho(\O)$,
$\chi_{\tt D}\in\rW^{1/s,r}(\Gamma)$,
$\xi_{i,\tt D}\in\rH^{1/2}(\Gamma)$, and
$f_i\in\rL^\varrho(\O)$, $i\in\{1,2\}$.
The problem reads: find
$(\bro,\u)\in\bX\times\bQ$,
$(\z,p)\in\bX_2\times\bQ_1$,
$(\bvphi,\chi)\in\rX_2\times\rM_1$, and
$(\vec{\bt}_i,\xi_i)\in\bcalH\times\calM$,
$i\in\{1,2\}$, such that
\begin{equation}\label{Biot-PNP-fully-mixed-formulation}
\begin{alignedat}{4}
\bfa_s(\bro,\bta)
&{}+\bfb_s(\bta,\u)&&&&=0
&\quad&\forall\,\bta\in\bX,
\\[0.75ex]
\bfb_s(\bro,\bv)
&&&&{}+\bfe_s(\z,\bv)&=\bF_{\bvphi,\bxi,\bt}(\bv)
&\quad&\forall\,\bv\in\bQ,
\\[0.75ex]
&&&\bfa_f(\z,\bw)+\bfd_1(\bw,p)
&&=0&\quad&\forall\,\bw\in\bX_1,
\\[0.75ex]
&&&\bfd_2(\z,q)&{}-\bfe_f((\bro,p),q)&=\bG(q)
&\quad&\forall\,q\in\bQ_2,
\\[1.5ex]
a(\bvphi,\bpsi)&{}+b_1(\bpsi,\chi)
&&&&=G(\bpsi)&\quad&\forall\,\bpsi\in\rX_1,
\\[0.75ex]
b_2(\bvphi,\gamma)&&&&&=F_{\bxi}(\gamma)
&\quad&\forall\,\gamma\in\rM_2,
\\[1.5ex]
\calA(\vec{\bt}_i,\vec{\bs}_i)&{}+\calB(\vec{\bs}_i,\xi_i)
&&&{}+\calE_{\z,\bvphi}(\vec{\bs}_i,\xi_i)
&=\mathcal G(\vec{\bs}_i)&\quad&\forall\,\vec{\bs}_i\in\bcalH,
\\[0.75ex]
\calB(\vec{\bt}_i,\eta_i)&{}-\calC(\xi_i,\eta_i)
&&&{}+\calD_{\z}(\xi_i,\eta_i)
&=\mathcal F(\eta_i)&\quad&\forall\,\eta_i\in\calM.
\end{alignedat}
\end{equation}
Here, $\bxi=(\xi_1,\xi_2)$ and
$\bt=(\bt_1,\bt_2)$ are formed from the
concentrations and gradients sought in the same system.

\section{Continuous stability and solvability}\label{section4}

Following \cite[Section~4]{gir-CAMWA-2025}, we establish block
stability for prescribed coupling fields, then impose the coupling
through a fixed point in the Darcy flux and electric displacement.
The mechanical analysis uses the new pseudostress formulation;
see also \cite{cgir-2024,ccgi-CAMWA-2023,gmr-IMANUM-2022}.

\subsection{Stability of the poroelastic block}\label{section41}

The poroelastic spaces defined in \eqref{def-bXi-bQi}
lead to distinct trial and test spaces for the Darcy
components. We collect them as
\begin{equation}\label{def-bbX-bbQ}
\mathbb X
:=
\bX\times\bQ\times\bX_2\times\bQ_1
\qan
\mathbb Q
:=
\bX\times\bQ\times\bX_1\times\bQ_2\,,
\end{equation}
endowed with the corresponding sum norms.
For $\vec{\bro}:=(\bro,\u,\z,p)\in\mathbb X$ and
$\vec{\bta}:=(\bta,\bv,\bw,q)\in\mathbb Q$
(cf. \eqref{def-bbX-bbQ}),
omitting the coupling forms $\bfe_s$ and $\bfe_f$
from \eqref{Biot-fully-mixed} gives the bilinear form
$\bfA:\mathbb X\times\mathbb Q\longrightarrow\mathbb R$
defined by
\begin{equation}\label{def-bfA}
\begin{aligned}
\bfA(\vec{\bro},\vec{\bta})
&:=\bfa_s(\bro,\bta)+\bfb_s(\bta,\u)
+\bfb_s(\bro,\bv)
+\bfa_f(\z,\bw)+\bfd_1(\bw,p)+\bfd_2(\z,q).
\end{aligned}
\end{equation}
The operator associated with this bilinear form
(cf. \eqref{def-bfA}) has the block structure
\begin{equation}\label{bfA-matrix-representation}
\bfA
=
\left(
\begin{array}{cc|cc}
\bfa_s & \bfb_s^\prime & 0 & 0
\\
\bfb_s & 0 & 0 & 0
\\ \hline
0 & 0 & \bfa_f & \bfd_1^\prime
\\
0 & 0 & \bfd_2 & 0
\end{array}
\right),
\end{equation}
where the entries denote the induced operators and primes their
transposes. The diagonal blocks in \eqref{bfA-matrix-representation}
are the uncoupled mechanical and generalized Darcy problems.
The pseudostress stability result in
\cite[Theorem~2.1]{ggs-CAMWA-2016} provides kernel coercivity of
$\bfa_s$ and the divergence inf--sup condition for $\bfb_s$,
uniformly in $\lambda$.
For the generalized Darcy block, the kernel estimate for $\bfa_f$
and the inf--sup condition for $\bfd_1$ follow from
\cite[Lemma~4.1]{gir-CAMWA-2025} and
\cite[Lemma~2.7]{gmr-IMANUM-2022}, respectively.
The inf--sup condition for $\bfd_2$ can be verified directly
under the Lipschitz-domain assumption used here.
Indeed, given $q\in\bQ_2\setminus\{0\}$, the bounded
divergence lifting of
\cite[Chapter~I, Corollary~2.4]{GiraultRaviart1986},
applied to $-q\in\rL_0^2(\O)$, provides
$\bw\in\bH_0^1(\O)$ such that
\[
\div(\bw)=-q,
\qquad
\|\bw\|_{1;\O}\le C\|q\|_{0;\O},
\]
where $C>0$ depends only on $\O$.
The embedding $\bH^1(\O)\hookrightarrow\bL^r(\O)$ for the
admissible exponents (cf. \eqref{r-rho-range}) implies that
$\bw\in\bX_2$ and $\|\bw\|_{\bX_2}\le C\|q\|_{0;\O}$.
Consequently, using the definition of $\bfd_2$
(cf. \eqref{def-af-d1-d2-ef}),
\[
\sup_{\bzeta\in\bX_2\setminus\{0\}}
\frac{\bfd_2(\bzeta,q)}{\|\bzeta\|_{\bX_2}}
\ge
\frac{\bfd_2(\bw,q)}{\|\bw\|_{\bX_2}}
=\frac{\|q\|_{0;\O}^2}{\|\bw\|_{\bX_2}}
\ge C^{-1}\|q\|_{\bQ_2}.
\]
The classical and generalized Babu\v{s}ka--Brezzi theories
\cite[Theorem~2.34]{eg-SPRINGER-2004},
\cite[Theorem~2.1 and Corollary~2.1]{BerCanMad}
yield the following global inf--sup estimates for
$\bfA$ (cf. \eqref{def-bfA}):
\begin{equation}\label{global-inf-sup-bfA}
\begin{aligned}
\sup_{\vec{\bta}\in\mathbb Q\setminus\{\mathbf0\}}
\frac{\bfA(\vec{\bro},\vec{\bta})}
     {\|\vec{\bta}\|_{\mathbb Q}}
&\ge
\balph_{\bfA}\|\vec{\bro}\|_{\mathbb X}
\quad
\forall\,\vec{\bro}\in\mathbb X
\quad\qan\quad
\sup_{\vec{\bro}\in\mathbb X\setminus\{\mathbf0\}}
\frac{\bfA(\vec{\bro},\vec{\bta})}
     {\|\vec{\bro}\|_{\mathbb X}}
&\ge
\balph_{\bfA}\|\vec{\bta}\|_{\mathbb Q}
\quad
\forall\,\vec{\bta}\in\mathbb Q.
\end{aligned}
\end{equation}
The common constant in \eqref{global-inf-sup-bfA}
is independent of $\lambda$ and may be chosen as
$\balph_{\bfA}=\frac12\min\{\balph_s,\balph_f\}$,
where $\balph_s,\balph_f>0$ denote the global inf--sup
lower bounds, in both orientations, for the mechanical
and generalized Darcy blocks, respectively.

We now restore the terms $\bfe_s$ and $-\bfe_f$ in
\eqref{Biot-fully-mixed}, defined in
\eqref{def-as-bs-es-bF} and \eqref{def-af-d1-d2-ef}.
The continuity bounds \eqref{Biot-bounds} ensure that their
combined operator norm is at most half the uncoupled stability
bound whenever
\begin{equation}\label{c0-alpha-lambda-assumption}
\max\left\{
C_r(\O)\frac{\alpha\nu}{\kappa},
C_r(\O)c_0,
\sqrt d\,\alpha c_\lambda
\right\}
\le
\frac{\balph_{\bfA}}{2}.
\end{equation}
The triangle inequality in the dual norm then gives
\begin{equation}\label{global-inf-sup-bfA+es-ef-1}
\sup_{\substack{\vec{\bta}\in\mathbb{Q}\\
                 \vec{\bta}\neq\mathbf{0}}}
\frac{
\bfA(\vec{\bro},\vec{\bta})
+\bfe_s(\z,\bv)
-\bfe_f((\bro,p),q)}
{\|\vec{\bta}\|_{\mathbb{Q}}}
\ge
\frac{\balph_{\bfA}}{2}
\|\vec{\bro}\|_{\mathbb{X}}
\qquad
\forall\,\vec{\bro}\in\mathbb{X}.
\end{equation}
The same bound holds with the trial and test arguments interchanged.
Thus the perturbed Biot operator is an isomorphism by the
Banach--Ne\v{c}as--Babu\v{s}ka theorem
\cite[Theorem~2.6]{eg-SPRINGER-2004}. Its inverse bound
$C_{\rm B}:=2/\balph_{\bfA}$ is independent of $\lambda$ under
\eqref{c0-alpha-lambda-assumption}.

\subsection{Stability of the electrochemical blocks}%\label{section42}

For Poisson, the trial and test pairs belong to
$\rX_2\times\rM_1$ and $\rX_1\times\rM_2$, respectively
(cf. \eqref{def:Poisson-spaces}). The results of
\cite[Section~4.3]{gir-CAMWA-2025}
and \cite[Lemmas~4.3--4.4]{cgrb-M2AN-2023}
yield
\begin{equation}\label{Poisson-global-stability}
\|\bvphi\|_{\rX_2}+\|\chi\|_{\rM_1}
\le
C_{\rm P}
\sup_{\substack{
(\bpsi,\gamma)\in\rX_1\times\rM_2\\
(\bpsi,\gamma)\neq\mathbf0}}
\frac{
a(\bvphi,\bpsi)
+b_1(\bpsi,\chi)
+b_2(\bvphi,\gamma)
}{
\|\bpsi\|_{\rX_1}+\|\gamma\|_{\rM_2}
}
\qquad
\forall\,(\bvphi,\chi)\in\rX_2\times\rM_1.
\end{equation}
Here $C_{\rm P}>0$ depends only on the physical parameters,
domain, exponents, and block stability constants.

For the ionic equations, prescribe $(\bw,\bphi)$ in place of
$(\z,\bvphi)$ and collect both equations of
\eqref{NP-mixed-formulation} into the bilinear form
$\boldsymbol\calA_{\bw,\bphi}:
(\bcalH\times\calM)\times(\bcalH\times\calM)
\longrightarrow\mathbb R$ defined by
\begin{equation}\label{def:ionic-block-form}
\begin{aligned}
\boldsymbol\calA_{\bw,\bphi}
((\vec{\br}_i,\xi_i),(\vec{\bs}_i,\eta_i))
:={}&\calA(\vec{\br}_i,\vec{\bs}_i)
+\calB(\vec{\bs}_i,\xi_i)+\calB(\vec{\br}_i,\eta_i)
-\calC(\xi_i,\eta_i)
+\calE_{\bw,\bphi}(\vec{\bs}_i,\xi_i)
+\calD_{\bw}(\xi_i,\eta_i).
\end{aligned}
\end{equation}
We separate the operator whose stability is known from the terms
depending on the prescribed fields:
\[
\begin{pmatrix}\calA&\calB'\\\calB&-\calC\end{pmatrix}
+\begin{pmatrix}0&\calE_{\bw,\bphi}'\\0&\calD_{\bw}\end{pmatrix}.
\]
The unperturbed stability follows from
\cite[Section~4.4]{gir-CAMWA-2025} and
\cite[Theorem~3.2]{gnr-JNM-2023}.

Let $\balph_{\boldsymbol\calA}>0$ be a common lower bound for both
species and both orientations.
The bounds \eqref{bound-E-D} show that the perturbation is small
enough to preserve stability whenever
\begin{equation}\label{assumption-w-bphi}
\|(\bw,\bphi)\|_{\bX_2\times\rX_2}
\le
\mathtt R
:=
\frac{\balph_{\boldsymbol\calA}}
{2\max\left\{
\|\calD\|,1,
\varepsilon^{-1}\kappa_1,
\varepsilon^{-1}\kappa_2
\right\}}.
\end{equation}
Under this condition,
\begin{equation}\label{global-inf-sup-Awbze-1}
\begin{aligned}
&
\sup_{\substack{(\vec{\bs}_i,\eta_i)\in\bcalH\times\calM\\
                 (\vec{\bs}_i,\eta_i)\neq\mathbf0}}
\frac{
\boldsymbol\calA_{\bw,\bphi}
\bigl((\vec{\br}_i,\xi_i),(\vec{\bs}_i,\eta_i)\bigr)}
{\|(\vec{\bs}_i,\eta_i)\|_{\bcalH\times\calM}}
\ge
\frac{\balph_{\boldsymbol\calA}}{2}
\|(\vec{\br}_i,\xi_i)\|_{\bcalH\times\calM}
\qquad
\forall\,(\vec{\br}_i,\xi_i)\in\bcalH\times\calM,
\end{aligned}
\end{equation}
and the transpose estimate holds as well.
In particular, its inverse bound is
$C_{\rm N}:=2/\balph_{\boldsymbol\calA}$, uniformly for both
species and all coefficient pairs satisfying
\eqref{assumption-w-bphi}.
For the contractivity and coupled error estimates, we introduce
the common coupling bound
\begin{equation}\label{def:ionic-coupling-bound}
D_*:=\|\calD\|+\max_{i\in\{1,2\}}\|\calE_i\|
=|\O|^{\frac12-\frac2\rho}
+\max\{1,\varepsilon^{-1}\kappa_1,\varepsilon^{-1}\kappa_2\},
\end{equation}
where $\calE_i$ denotes $\calE$ for species $i$
(cf. \eqref{bounds-NP-D-E}).
Thus, $D_*$ depends only on the domain, exponents,
and physical coefficients.

\subsection{Subproblem estimates and the coupled solution}
\label{section44}

We enforce the coupling through a fixed-point operator acting
on the Darcy flux and electric displacement.
With the admissible radius $\mathtt R$ defined in
\eqref{assumption-w-bphi}, introduce the closed ball
\begin{equation}\label{def:fixed-point-ball}
\rW(\mathtt R)
:=
\bigl\{
(\bw,\bphi)\in\bX_2\times\rX_2:\quad
\|(\bw,\bphi)\|_{\bX_2\times\rX_2}\le\mathtt R
\bigr\}.
\end{equation}

For given $(\bphi,\bet,\br)\in
\rX_2\times\bcalM\times\bcalH_1$, define
$\bS(\bphi,\bet,\br):=\z$, where $\z\in\bX_2$
is the Darcy flux obtained by solving
\eqref{Biot-fully-mixed} with load $\bF_{\bphi,\bet,\br}$.
For given $\bet\in\bcalM$, define
$\tbS(\bet):=\bvphi$, where $\bvphi\in\rX_2$
is the electric displacement obtained by solving
\eqref{Poisson-mixed-formulation} with load $F_\bet$.
Finally, for $(\bw,\bphi)\in\rW(\mathtt R)$
(cf.~\eqref{def:fixed-point-ball}), define
$\bT(\bw,\bphi):=(\bt_1,\bt_2)\in\bcalH_1$
and $\bXi(\bw,\bphi):=(\xi_1,\xi_2)\in\bcalM$,
where $\bt_i$ and $\xi_i$ are the concentration gradient
and concentration obtained from
\eqref{NP-mixed-formulation} for species $i$
with coefficient fields $(\bw,\bphi)$. These operators define
$\bPi:\rW(\mathtt R)\longrightarrow\bX_2\times\rX_2$ by
\begin{equation}\label{Fixed-point-operator}
\bPi(\bw,\bphi)
:=
\Bigl(
\bS\bigl(\bphi,\bXi(\bw,\bphi),\bT(\bw,\bphi)\bigr),
\tbS\bigl(\bXi(\bw,\bphi)\bigr)
\Bigr).
\end{equation}
Each evaluation solves the ionic problems first, then updates the
Darcy flux and electric displacement through Biot and Poisson.
Thus the coupled problem in $\rW(\mathtt R)$ is equivalent to
$\bPi(\z,\bvphi)=(\z,\bvphi)$.
The following result establishes that these operators are well defined.

\begin{lemma}\label{bS-well-defined}
Under \eqref{c0-alpha-lambda-assumption}, the Biot problem \eqref{Biot-fully-mixed} and
the Poisson problem \eqref{Poisson-mixed-formulation},
with given inputs, are uniquely solvable in their
prescribed spaces.
For each $i\in\{1,2\}$, the ionic problem
\eqref{NP-mixed-formulation} is uniquely solvable
for coefficient pairs satisfying \eqref{assumption-w-bphi}.
Moreover, there exist positive constants
$C_S$, $C_E$, and $C_T$,
independent of $\lambda$, such that
\begin{align}
\|(\bro,\u,\z,p)\|_{\mathbb X}
&\le
C_S
\bigl\{
\|\ff\|_{0;\O}+\|g\|_{0;\O}
+\bigl(\|\bet\|_{\bcalM}+\Lambda_\chi\bigr)\|\bphi\|_{\rX_2}
+\|\br\|_{\bcalH_1}
\bigr\},
\label{bS-apriori-estimation}
\\[1ex]
\|(\bvphi,\chi)\|_{\rX_2\times\rM_1}
&\le
C_E
\bigl\{
\|\chi_{\tt D}\|_{1/s,r;\Gamma}
+\|\bet\|_{\bcalM}+\Lambda_\chi
\bigr\},
\label{bound-apriori-tbS-bsi}
\\[1ex]
\|(\vec{\bt}_i,\xi_i)\|_{\bcalH\times\calM}
&\le
C_T
\bigl\{
\|\xi_{i,\tt D}\|_{1/2;\Gamma}
+\|f_i\|_{0,\varrho;\O}
\bigr\},
\qquad i\in\{1,2\}.
\label{bound-apriori-Xi-bT-i}
\end{align}
These constants are uniform over the prescribed input fields
and independent of the data under the stated stability hypotheses.
\end{lemma}
\begin{proof}
Existence and uniqueness follow from the
Banach--Ne\v{c}as--Babu\v{s}ka theorem
\cite[Theorem~2.6]{eg-SPRINGER-2004},
together with the global inf--sup conditions
\eqref{global-inf-sup-bfA+es-ef-1},
\eqref{Poisson-global-stability}, and
\eqref{global-inf-sup-Awbze-1},
and the corresponding transpose estimates.
Applying these inverse bounds to the loads, then using
\eqref{Biot-bounds}, \eqref{Biot-F-bounds}, \eqref{bound-F-xi},
\eqref{bound-G}, and \eqref{bounds-NP-A-B-C-G-F}, proves
\eqref{bS-apriori-estimation}--\eqref{bound-apriori-Xi-bT-i}.
For Biot one may take $C_S=C_{\rm B}\max\{1,\varepsilon^{-1},\delta\}$.
The two-species bounds use the sum norms.
\end{proof}

For every $(\bw,\bphi)\in\rW(\mathtt R)$, the subproblem
estimates of Lemma~\ref{bS-well-defined}
(cf. \eqref{bS-apriori-estimation},
\eqref{bound-apriori-tbS-bsi}, and
\eqref{bound-apriori-Xi-bT-i}) give
\begin{equation}\label{eq:fixed-point-invariance-bound}
\begin{aligned}
\|\bPi(\bw,\bphi)\|_{\bX_2\times\rX_2}
\le C(\mathtt R)\Biggl\{
&\|\ff\|_{0;\O}+\|g\|_{0;\O}
+\|\chi_{\tt D}\|_{1/s,r;\Gamma}+\Lambda_\chi
+\sum_{i=1}^{2}
\bigl(
\|\xi_{i,\tt D}\|_{1/2;\Gamma}
+\|f_i\|_{0,\varrho;\O}
\bigr)
\Biggr\},
\end{aligned}
\end{equation}
where $C(\mathtt R):=(C_S+C_E)(1+C_T)(1+\mathtt R)$.
For contractivity, splitting the concentration--field product
and applying H\"older's inequality gives
\begin{equation}\label{Biot-F-difference-bound}
\begin{aligned}
\|\bF_{\bvphi,\bxi,\bt}
-\bF_{\bphi,\bet,\br}\|_{\bQ'}
\le{}
\varepsilon^{-1}\bigl(\|\bxi\|_{\bcalM}+\Lambda_\chi\bigr)
\|\bvphi-\bphi\|_{\rX_2}
+\varepsilon^{-1}\|\bphi\|_{\rX_2}
\|\bxi-\bet\|_{\bcalM}
+\delta\|\bt-\br\|_{\bcalH_1}.
\end{aligned}
\end{equation}
The Poisson and ionic difference estimates are those of
\cite[Lemmas~4.8--4.9]{gir-CAMWA-2025}.
Using the ionic bound
\eqref{bound-apriori-Xi-bT-i}
to estimate the concentration norms appearing in these
inequalities yields
\begin{equation}\label{LC-bPi}
\begin{aligned}
&\|\bPi(\z,\bvphi)-\bPi(\bw,\bphi)\|_{\bX_2\times\rX_2}
\le
\Biggl\{L_{\bPi}
\sum_{i=1}^{2}
\bigl(
\|\xi_{i,\tt D}\|_{1/2;\Gamma}
+\|f_i\|_{0,\varrho;\O}
\bigr)
+C_{\rm B}\varepsilon^{-1}\Lambda_\chi\Biggr\}
\|(\z,\bvphi)-(\bw,\bphi)\|_{\bX_2\times\rX_2},
\end{aligned}
\end{equation}
for all
$(\z,\bvphi),(\bw,\bphi)\in\rW(\mathtt R)$,
where one may take
\begin{equation}\label{def:fixed-point-Lipschitz-constant}
L_{\bPi}
:=
C_T\Bigl\{
C_{\rm B}\varepsilon^{-1}
+C_{\rm N}D_*
\bigl[
C_{\rm B}(\varepsilon^{-1}\mathtt R+\delta)
+C_{\rm P}\|F\|
\bigr]
\Bigr\}.
\end{equation}
Here $D_*$ is defined in
\eqref{def:ionic-coupling-bound}.
The fixed charge adds only the term $C_{\rm B}\varepsilon^{-1}\Lambda_\chi$ in
\eqref{LC-bPi}: it enters the Poisson load additively, so it cancels in
differences, and it enters the Maxwell force linearly in $\bvphi$.
The expression \eqref{def:fixed-point-Lipschitz-constant}
follows from the ionic stability estimate
\eqref{global-inf-sup-Awbze-1}, the coefficient-difference
bounds \eqref{bound-E-D}, and the load estimates
\eqref{Biot-F-difference-bound} and \eqref{bound-F-xi}.
Thus, $C(\mathtt R)$ and $L_{\bPi}$ depend only on the
block continuity and stability constants and the chosen
radius, independently of the prescribed data.

\begin{theorem}\label{Thm:unique-solution}
Assume \eqref{c0-alpha-lambda-assumption} and
\begin{equation}\label{small-data-assumption-delta}
\begin{array}{c}
\ds C(\mathtt R)\Biggl\{
\|\ff\|_{0;\O}+\|g\|_{0;\O}
+\|\chi_{\tt D}\|_{1/s,r;\Gamma}+\Lambda_\chi
+\sum_{i=1}^{2}
\bigl(
\|\xi_{i,\tt D}\|_{1/2;\Gamma}
+\|f_i\|_{0,\varrho;\O}
\bigr)
\Biggr\}
\le\mathtt R,
\\
\ds L_{\bPi}
\sum_{i=1}^{2}
\bigl(
\|\xi_{i,\tt D}\|_{1/2;\Gamma}
+\|f_i\|_{0,\varrho;\O}
\bigr)
+C_{\rm B}\varepsilon^{-1}\Lambda_\chi
<1.
\end{array}
\end{equation}
Then \eqref{Biot-PNP-fully-mixed-formulation}
has a unique solution within the class
$(\z,\bvphi)\in\rW(\mathtt R)$.
Moreover, there exist positive constants $C_S$, $C_E$, and $C_T$,
independent of the solution, the prescribed data, and $\lambda$
under the stated stability hypotheses, such that
\begin{equation}\label{bound:a-priori}
\begin{aligned}
\|(\bro,\u,\z,p)\|_{\mathbb X}
&\le
C_S\Biggl\{
\|\ff\|_{0;\O}+\|g\|_{0;\O}
+\mathtt R\Lambda_\chi
+C_T(\mathtt R+1)
\sum_{i=1}^{2}
\bigl(
\|\xi_{i,\tt D}\|_{1/2;\Gamma}
+\|f_i\|_{0,\varrho;\O}
\bigr)
\Biggr\},
\\[1ex]
\|(\bvphi,\chi)\|_{\rX_2\times\rM_1}
&\le
C_E\Biggl\{
\|\chi_{\tt D}\|_{1/s,r;\Gamma}
+\Lambda_\chi
+C_T
\sum_{i=1}^{2}
\bigl(
\|\xi_{i,\tt D}\|_{1/2;\Gamma}
+\|f_i\|_{0,\varrho;\O}
\bigr)
\Biggr\},
\\[1ex]
\|(\vec{\bt}_i,\xi_i)\|_{\bcalH\times\calM}
&\le
C_T
\bigl\{
\|\xi_{i,\tt D}\|_{1/2;\Gamma}
+\|f_i\|_{0,\varrho;\O}
\bigr\},
\qquad i\in\{1,2\}.
\end{aligned}
\end{equation}
\end{theorem}
\begin{proof}
The first condition in
\eqref{small-data-assumption-delta}
ensures that $\bPi$ maps $\rW(\mathtt R)$ into itself
(cf. \eqref{eq:fixed-point-invariance-bound}).
The second condition and \eqref{LC-bPi} make it
contractive.
Banach's fixed-point theorem therefore gives a unique
fixed point, and the three subproblems uniquely
determine the remaining variables.
Substituting the ionic bounds
\eqref{bound-apriori-Xi-bT-i}
and $\|\bvphi\|_{\rX_2}\le\mathtt R$ into
\eqref{bS-apriori-estimation} and
\eqref{bound-apriori-tbS-bsi}
proves \eqref{bound:a-priori}.
\end{proof}

\begin{remark}
Condition \eqref{small-data-assumption-delta} is a sufficient
weak-coupling (small-data) regime: it makes $\bPi$ a contraction on
$\rW(\mathtt R)$, and it is not claimed to be sharp.
Beyond it, the fixed-point argument gives neither existence nor uniqueness.
Accordingly, Example~3 below uses load continuation to move
away from this regime.
\end{remark}

\section{Galerkin approximation and a priori analysis}\label{section5}

We approximate \eqref{Biot-PNP-fully-mixed-formulation} by
the conforming finite element subspaces of Section~\ref{section51}
and use the block structure to establish discrete solvability
and error bounds.

\subsection{Finite element spaces and the discrete problem}\label{section51}
Let $\O$ be polygonal or polyhedral and $\{\Th\}$ a shape-regular
family of conforming simplicial meshes. Write $h_K=\operatorname{diam}K$
and $h=\max_K h_K$. For $\ell\ge0$, let $\rP_\ell(\Th)$ denote
the discontinuous piecewise polynomials of degree at most $\ell$, and
$\mathbf{RT}_\ell(\Th)$ the Raviart--Thomas space, with local space
$\mathbf{RT}_\ell(K)=\bP_\ell(K)
\oplus\widetilde{\rP}_\ell(K)\bx$. Here $\widetilde{\rP}_\ell(K)$ denotes the homogeneous
polynomials of degree $\ell$.
The tensor space $\mathbb{RT}_\ell(\Th)$ is defined row-wise.
The conforming subspaces are
\begin{equation}\label{def-FE}
\begin{array}{l}
\bX_h \,:=\, \mathbb{RT}_{\ell}(\Th)\cap \bbH_0
(\bdiv;\O)\,, \quad \bQ_h \,:=\, \bP_{\ell}(\Th) \,,\\[0.5ex] 
\bX_{2,h} \,:=\, \RTl(\Th)\cap \bH^r_0(\div;\O)\,, \quad \bX_{1,h} \,:=\, \RTl(\Th)\cap \bH^s_0(\div_s;\O)\,,\\[0.5ex]
\bQ_{1,h} \,:=\, \rP_\ell(\Th)\cap \rL_0^r(\O)\,, \quad \bQ_{2,h} \,:=\, \rP_\ell(\Th)\cap \rL_0^2(\O) \,,\\[0.5ex]
\rX_{2,h}\,=\,\rX_{1,h} \,:=\, \RTl(\Th)\,, \quad \rM_{1,h} \,=\, \rM_{2,h} \,:=\, \Pl(\Th)\,, \\[0.5ex]
\calH_{1,h} \,:=\, \bPl(\Th)\,, \quad \calH_{2,h} \,:=\, \RTl(\Th)\,, \qan \calM_h \,:=\, \Pl(\Th)\,.
\end{array}
\end{equation}
We also set
$\bcalH_h:=\calH_{1,h}\times\calH_{2,h}$,
$\bcalH_{1,h}:=\calH_{1,h}\times\calH_{1,h}$, and
$\bcalM_h:=\calM_h\times\calM_h$.
The spaces in \eqref{def-FE} are endowed with the norms
of their continuous counterparts
(cf. \eqref{def-bXi-bQi}, \eqref{def:Poisson-spaces},
and \eqref{def:bH-H1-H2-M}).

The Galerkin approximation of
\eqref{Biot-PNP-fully-mixed-formulation}
reads: find
$(\bro_h,\u_h)\in\bX_h\times\bQ_h$,
$(\z_h,p_h)\in\bX_{2,h}\times\bQ_{1,h}$,
$(\bvphi_h,\chi_h)\in\rX_{2,h}\times\rM_{1,h}$,
and, for each $i\in\{1,2\}$,
$(\vec{\bt}_{i,h},\xi_{i,h})\in\bcalH_h\times\calM_h$,
such that
\begin{equation}\label{Biot-PNP-fully-mixed-formulation-disc}
\begin{array}{rclllcl}
\bfa_s(\bro_h,\bta_h)
&+&\bfb_s(\bta_h,\u_h)
&&
&=&0,
\\[0.75ex]
\bfb_s(\bro_h,\bv_h)
&&
&&+\bfe_s(\z_h,\bv_h)
&=&\bF_{\bvphi_h,\bxi_h,\bt_h}(\bv_h),
\\[0.75ex]
&&
&\bfa_f(\z_h,\bw_h)+\bfd_1(\bw_h,p_h)
&
&=&0,
\\[0.75ex]
&&
&\bfd_2(\z_h,q_h)
&-\bfe_f((\bro_h,p_h),q_h)
&=&\bG(q_h),
\\[0.5ex]
a(\bvphi_h,\bpsi_h)
&+&b_1(\bpsi_h,\chi_h)
&&
&=&G(\bpsi_h),
\\[0.75ex]
b_2(\bvphi_h,\gamma_h)
&&
&&
&=&F_{\bxi_h}(\gamma_h),
\\[0.5ex]
\calA(\vec{\bt}_{i,h},\vec{\bs}_{i,h})
&+&\calB(\vec{\bs}_{i,h},\xi_{i,h})
&&+\calE_{\z_h,\bvphi_h}(\vec{\bs}_{i,h},\xi_{i,h})
&=&\mathcal G(\vec{\bs}_{i,h}),
\\[0.75ex]
\calB(\vec{\bt}_{i,h},\eta_{i,h})
&-&\calC(\xi_{i,h},\eta_{i,h})
&&+\calD_{\z_h}(\xi_{i,h},\eta_{i,h})
&=&\mathcal F(\eta_{i,h}),
\end{array}
\end{equation}
for every
$(\bta_h,\bv_h)\in\bX_h\times\bQ_h$,
$(\bw_h,q_h)\in\bX_{1,h}\times\bQ_{2,h}$,
$(\bpsi_h,\gamma_h)\in\rX_{1,h}\times\rM_{2,h}$,
and
$(\vec{\bs}_{i,h},\eta_{i,h})\in\bcalH_h\times\calM_h$.
Here,
$\vec{\bt}_{i,h}:=(\bt_{i,h},\bsi_{i,h})$,
$\vec{\bs}_{i,h}:=(\bs_{i,h},\bta_{i,h})$,
$\bxi_h:=(\xi_{1,h},\xi_{2,h})$, and
$\bt_h:=(\bt_{1,h},\bt_{2,h})$.

\subsection{Uniform stability and discrete solvability}
\label{section53}

We begin by establishing the discrete stability of the two diagonal
blocks in \eqref{bfA-matrix-representation}, using the
conforming spaces defined in \eqref{def-FE}.
For the mechanical block, the Raviart--Thomas results in
\cite[Lemma~3.2]{gms-CMAME-2010}
and \cite[Section~3]{ggs-CAMWA-2016},
together with the continuous estimates recalled in
Section~\ref{section41}, provide discrete stability
with constants independent of $h$ and $\lambda$.
The generalized Darcy block coincides with that analyzed
in \cite[Section~6.2]{gir-CAMWA-2025}.
The discrete divergence inf--sup condition for $\bfd_1$
is provided by \cite[Lemma~4.4]{gmr-IMANUM-2022}, under the
domain and exponent assumptions of that result.
For $\bfd_2$, we follow the argument of
\cite[Lemma~6.1]{gir-CAMWA-2025}, replacing its initial
divergence lifting with the bounded right inverse of the divergence
in \cite[Chapter~I, Corollary~2.4]{GiraultRaviart1986}.
The required Raviart--Thomas interpolation estimate
from $\bH^1(\O)$ to $\bL^r(\O)$ follows from
\cite[Lemma~A.3]{ccgi-CAMWA-2023}.
For the admissible exponents
(cf. \eqref{r-rho-range}), this establishes the
discrete inf--sup condition for $\bfd_2$ uniformly
in $h$.
In two dimensions, the kernel inf--sup condition for
$\bfa_f$ is established in
\cite[Lemma~4.3]{gmr-IMANUM-2022}.
Under the assumptions ensuring that these block estimates
are uniform in $h$, the generalized
Babu\v{s}ka--Brezzi theory
\cite[Theorem~2.1 and Corollary~2.1]{BerCanMad}
yields
\begin{equation}\label{discrete-inf-sup-Af-hypothesis}
\sup_{\substack{
(\bw_h,q_h)\in\bX_{1,h}\times\bQ_{2,h}\\
(\bw_h,q_h)\neq\mathbf0}}
\frac{
\bfA_f((\z_h,p_h),(\bw_h,q_h))
}{
\|(\bw_h,q_h)\|_{\bX_1\times\bQ_2}
}
\ge
\widehat{\balph}_f
\|(\z_h,p_h)\|_{\bX_2\times\bQ_1}
\qquad
\forall\,(\z_h,p_h)\in\bX_{2,h}\times\bQ_{1,h},
\end{equation}
where $\widehat{\balph}_f>0$ is independent of $h$ and
\[
\bfA_f((\z,p),(\bw,q))
:=
\bfa_f(\z,\bw)
+\bfd_1(\bw,p)
+\bfd_2(\z,q).
\]
\begin{remark}
In two dimensions, the kernel inf--sup condition for $\bfa_f$
follows from the discrete projector estimates in
\cite[Lemmas~4.2--4.3]{gmr-IMANUM-2022},
under the domain and mesh assumptions ensuring uniformity in $h$.
For non-convex domains with $\ell=0$, the cited projector bound
contains a logarithmic meshsize factor, so uniformity in $h$ is
not guaranteed by these references for such domains.
The corresponding uniform three-dimensional kernel estimate
is not supplied by these references either.
In these cases, we take
\eqref{discrete-inf-sup-Af-hypothesis} as an explicit
hypothesis of the analysis.
\end{remark}
With the mechanical and Darcy stability bounds in place, let
$\widehat{\balph}_{\bfA}>0$ be a common discrete inf--sup lower bound
for the uncoupled Biot operator, in both orientations.
The equal dimensions of its trial and test spaces justify the
transpose bound. Under
\begin{equation}\label{c0-alpha-lambda-assumption-disc}
\max\left\{
C_r(\O)\frac{\alpha\nu}{\kappa},
C_r(\O)c_0,
\sqrt d\,\alpha c_\lambda
\right\}
\le
\frac{\widehat{\balph}_{\bfA}}{2},
\end{equation}
the perturbation argument leading to
\eqref{global-inf-sup-bfA+es-ef-1} gives the
lower bound $\widehat{\balph}_{\bfA}/2$ for the coupled Biot block,
uniformly in $h$ and $\lambda$.

For Poisson, assume the uniform discrete stability conditions
of \cite[Sections~5.2 and~6.2]{cgrb-M2AN-2023}, invoked under
their domain and mesh hypotheses; see also
\cite[Lemma~5.2]{ccgi-2022-16} and
\cite[Lemma~4.5]{gmr-IMANUM-2022}.
In three dimensions, the uniform kernel stability of $a$ is
an additional hypothesis.
The ionic spaces coincide with those of
\cite[Section~6.1]{gir-CAMWA-2025}, and satisfy the required
unperturbed stability conditions verified in its Section~6.2.
Let $\widehat{\balph}_{\boldsymbol\calA}>0$ be a common lower bound
for the two unperturbed ionic blocks. Their perturbed versions
are uniformly stable whenever
\begin{equation}\label{assumption-w-bphi-disc}
\|(\bw_h,\bphi_h)\|_{\bX_2\times\rX_2}
\le
\mathtt R_{\dd}
:=
\frac{\widehat{\balph}_{\boldsymbol\calA}}
     {2\max\bigl\{
       \|\calD\|,\,1,\,
       \varepsilon^{-1}\kappa_1,\,
       \varepsilon^{-1}\kappa_2
       \bigr\}}.
\end{equation}

Define $\bS_h$, $\tbS_h$, $\bT_h$, and $\bXi_h$ by the block
problems of Section~\ref{section44} on the conforming subspaces.
Their estimates follow from the uniform bounds above and the
continuous load estimates; see also
\cite[Sections~5.4--5.5]{gir-CAMWA-2025}.

\begin{lemma}\label{bSh-well-defined}
Under the uniform stability assumptions of
Section~\ref{section53} and
\eqref{c0-alpha-lambda-assumption-disc},
the discrete Biot and Poisson problems with given inputs
are uniquely solvable in their prescribed spaces.
The discrete ionic problems are uniquely solvable for
coefficient pairs satisfying
\eqref{assumption-w-bphi-disc}.
Moreover, there exist positive constants
$\widehat{C}_S$, $\widehat{C}_E$, and $\widehat{C}_T$,
independent of $h$ and $\lambda$, such that
\begin{align}
\|(\bro_h,\u_h,\z_h,p_h)\|_{\mathbb X}
&\le
\widehat{C}_S
\bigl\{
\|\ff\|_{0;\O}+\|g\|_{0;\O}
+\bigl(\|\bet_h\|_{\bcalM}+\Lambda_\chi\bigr)\|\bphi_h\|_{\rX_2}
+\|\br_h\|_{\bcalH_1}
\bigr\},
\label{bound-apriori-bSh}
\\[1ex]
\|(\bvphi_h,\chi_h)\|_{\rX_2\times\rM_1}
&\le
\widehat{C}_E
\bigl\{
\|\chi_{\tt D}\|_{1/s,r;\Gamma}
+\|\bet_h\|_{\bcalM}+\Lambda_\chi
\bigr\},
\label{bound-apriori-tbSh-bsi}
\\[1ex]
\|(\vec{\bt}_{i,h},\xi_{i,h})\|_{\bcalH\times\calM}
&\le
\widehat{C}_T
\bigl\{
\|\xi_{i,\tt D}\|_{1/2;\Gamma}
+\|f_i\|_{0,\varrho;\O}
\bigr\},
\qquad i\in\{1,2\}.
\label{bound-apriori-Xiih-bTih-discrete}
\end{align}
These constants are uniform over the prescribed input fields
and independent of the data under the stated stability hypotheses.
\end{lemma}

\begin{proof}
The finite-dimensional Banach--Ne\v{c}as--Babu\v{s}ka theorem
\cite[Theorem~2.6]{eg-SPRINGER-2004}, together with the uniform
discrete inf--sup bounds of Section~\ref{section53} and their
transpose estimates, gives existence, uniqueness, and the
corresponding inverse bounds.
Applying the continuous load estimates on the conforming
subspaces proves the assertions.
For the Biot block, one may take
$\widehat{C}_S=2\widehat{\balph}_{\bfA}^{-1}
\max\{1,\varepsilon^{-1},\delta\}$.
\end{proof}

The discrete coupling is imposed through the same sequence of
block solves as in Section~\ref{section44}.
Let $\rW(\mathtt R_{\dd})$ denote the closed ball of radius
$\mathtt R_{\dd}$ in $\bX_{2,h}\times\rX_{2,h}$, and define
$\bPi_h$ by the composition
\eqref{Fixed-point-operator} with the four discrete solution
operators.
The subproblem estimates of Lemma~\ref{bSh-well-defined}
(cf. \eqref{bound-apriori-bSh},
\eqref{bound-apriori-tbSh-bsi}, and
\eqref{bound-apriori-Xiih-bTih-discrete})
provide the invariance bound with
\[
\widehat C(\mathtt R_{\dd})
:=(\widehat C_S+\widehat C_E)(1+\widehat C_T)(1+\mathtt R_{\dd}).
\]
For contractivity, let
$(\z,\bvphi),(\bw,\bphi)\in\rW(\mathtt R_{\dd})$ and set
$\bxi:=\bXi_h(\z,\bvphi)$, $\bet:=\bXi_h(\bw,\bphi)$,
$\bt:=\bT_h(\z,\bvphi)$, $\br:=\bT_h(\bw,\bphi)$, and
$\mathtt d:=\|(\z,\bvphi)-(\bw,\bphi)\|_{\bX_2\times\rX_2}$,
$\mathtt D_{\rm data}:=\sum_{i=1}^{2}
\bigl(\|\xi_{i,\tt D}\|_{1/2;\Gamma}+\|f_i\|_{0,\varrho;\O}\bigr)$.
Let $\widehat C_{\rm B}:=2/\widehat{\balph}_{\bfA}$ and
$\widehat C_{\rm N}:=2/\widehat{\balph}_{\boldsymbol\calA}$ be the
inverse bounds of the discrete coupled Biot and ionic operators,
and let $\widehat C_{\rm P}$ be the discrete Poisson stability
constant. All three are independent of $h$.

\emph{Ionic step.}
Since $\calE$ and $\calD$ are linear in their coefficient fields,
the difference of the ionic solutions for the two coefficient pairs
satisfies, for each $i\in\{1,2\}$ and every
$(\vec{\bs}_i,\theta_i)\in\bcalH_h\times\calM_h$,
\[
\boldsymbol\calA_{\z,\bvphi}
\bigl((\vec{\bt}_{i,h}-\vec{\br}_{i,h},\xi_{i,h}-\eta_{i,h}),
(\vec{\bs}_i,\theta_i)\bigr)
=
-\calE_{\z-\bw,\bvphi-\bphi}(\vec{\bs}_i,\eta_{i,h})
-\calD_{\z-\bw}(\eta_{i,h},\theta_i).
\]
Using \eqref{bound-E-D}, the definition of $D_*$ in
\eqref{def:ionic-coupling-bound}, the discrete stability
\eqref{assumption-w-bphi-disc} and the discrete ionic bound
\eqref{bound-apriori-Xiih-bTih-discrete} for $\eta_{i,h}$,
\begin{equation}\label{eq:discrete-ionic-difference}
\|(\vec{\bt}_{i,h}-\vec{\br}_{i,h},\xi_{i,h}-\eta_{i,h})\|_{\bcalH\times\calM}
\le
\widehat C_{\rm N}D_*\widehat C_{\rm T}
\bigl(\|\xi_{i,\tt D}\|_{1/2;\Gamma}+\|f_i\|_{0,\varrho;\O}\bigr)\,\mathtt d .
\end{equation}
Summing over the two species,
$\|\bxi-\bet\|_{\bcalM}+\|\bt-\br\|_{\bcalH_1}
\le \widehat C_{\rm N}D_*\widehat C_{\rm T}\,
\mathtt D_{\rm data}\,\mathtt d$.

\emph{Poisson and Biot steps.}
By the discrete Poisson stability and \eqref{bound-F-xi},
$\|\tbS_h(\bxi)-\tbS_h(\bet)\|_{\rX_2}
\le\widehat C_{\rm P}\|F\|\,\|\bxi-\bet\|_{\bcalM}$.
For the Biot step, the discrete analogue of
\eqref{global-inf-sup-bfA+es-ef-1}, with inverse bound
$\widehat C_{\rm B}$, together with
\eqref{Biot-F-difference-bound}, the bound
$\|\bxi\|_{\bcalM}\le\widehat C_{\rm T}\mathtt D_{\rm data}$ from
\eqref{bound-apriori-Xiih-bTih-discrete}, and
$\|\bphi\|_{\rX_2}\le\mathtt R_{\dd}$, gives
\[
\|\bS_h(\bvphi,\bxi,\bt)-\bS_h(\bphi,\bet,\br)\|_{\bX_2}
\le
\widehat C_{\rm B}\Bigl\{
\varepsilon^{-1}\bigl(\widehat C_{\rm T}\mathtt D_{\rm data}+\Lambda_\chi\bigr)\mathtt d
+\varepsilon^{-1}\mathtt R_{\dd}\|\bxi-\bet\|_{\bcalM}
+\delta\|\bt-\br\|_{\bcalH_1}
\Bigr\}.
\]
Adding the two components of $\bPi_h$ and inserting
\eqref{eq:discrete-ionic-difference} yields
\begin{equation}\label{LC-bPi-discrete}
\|\bPi_h(\z,\bvphi)-\bPi_h(\bw,\bphi)\|_{\bX_2\times\rX_2}
\le
\bigl(\widehat L_{\bPi}\,\mathtt D_{\rm data}
+\widehat C_{\rm B}\varepsilon^{-1}\Lambda_\chi\bigr)\mathtt d,
\end{equation}
with
\begin{equation}\label{def:discrete-Lipschitz-constant}
\widehat L_{\bPi}
:=
\widehat C_{\rm T}\Bigl\{
\widehat C_{\rm B}\varepsilon^{-1}
+\widehat C_{\rm N}D_*
\bigl[
\widehat C_{\rm B}(\varepsilon^{-1}\mathtt R_{\dd}+\delta)
+\widehat C_{\rm P}\|F\|
\bigr]
\Bigr\}.
\end{equation}
This is the discrete counterpart of \eqref{def:fixed-point-Lipschitz-constant}.
Both $\widehat C(\mathtt R_{\dd})$ and $\widehat L_{\bPi}$ are
independent of $h$ and of the data under the stated uniform
stability hypotheses.

\begin{theorem}\label{Thm:discrete-unique-solution}
Assume the uniform stability conditions of this subsection,
\eqref{c0-alpha-lambda-assumption-disc}, and

\begin{equation}\label{small-data-assumption-delta0-discrete}
\begin{aligned}
\widehat C(\mathtt R_{\dd})\Biggl\{
&\|\ff\|_{0;\O}+\|g\|_{0;\O}
+\|\chi_{\tt D}\|_{1/s,r;\Gamma}+\Lambda_\chi
+\sum_{i=1}^{2}
\bigl(
\|\xi_{i,\tt D}\|_{1/2;\Gamma}
+\|f_i\|_{0,\varrho;\O}
\bigr)
\Biggr\}\le\mathtt R_{\dd},
\\[1ex]
\widehat{L}_{\bPi}
&\sum_{i=1}^{2}
\bigl(
\|\xi_{i,\tt D}\|_{1/2;\Gamma}
+\|f_i\|_{0,\varrho;\O}
\bigr)+\widehat C_{\rm B}\varepsilon^{-1}\Lambda_\chi<1.
\end{aligned}
\end{equation}
Then \eqref{Biot-PNP-fully-mixed-formulation-disc} has a unique
solution within the class
$(\z_h,\bvphi_h)\in\rW(\mathtt R_{\dd})$.
Moreover, there exist positive constants
$\widehat{C}_S$, $\widehat{C}_E$, and $\widehat{C}_T$,
independent of $h$, the solution, the prescribed data,
and $\lambda$ under the stated uniform stability hypotheses,
such that
\begin{equation}\label{bound:a-priori-discr}
\begin{aligned}
\|(\bro_h,\u_h,\z_h,p_h)\|_{\mathbb X}
&\le \widehat{C}_S\Biggl\{
\|\ff\|_{0;\O}+\|g\|_{0;\O}
+\mathtt R_{\dd}\Lambda_\chi
+\widehat{C}_T(\mathtt R_{\dd}+1)
\sum_{i=1}^{2}
\bigl(
\|\xi_{i,\tt D}\|_{1/2;\Gamma}
+\|f_i\|_{0,\varrho;\O}
\bigr)\Biggr\},
\\[1ex]
\|(\bvphi_h,\chi_h)\|_{\rX_2\times\rM_1}
&\le \widehat{C}_E\Biggl\{
\|\chi_{\tt D}\|_{1/s,r;\Gamma}
+\Lambda_\chi
+\widehat{C}_T
\sum_{i=1}^{2}
\bigl(
\|\xi_{i,\tt D}\|_{1/2;\Gamma}
+\|f_i\|_{0,\varrho;\O}
\bigr)\Biggr\},
\\[1ex]
\|(\vec{\bt}_{i,h},\xi_{i,h})\|_{\bcalH\times\calM}
&\le \widehat{C}_T
\bigl\{
\|\xi_{i,\tt D}\|_{1/2;\Gamma}
+\|f_i\|_{0,\varrho;\O}
\bigr\},\qquad i\in\{1,2\}.
\end{aligned}
\end{equation}
\end{theorem}
\begin{proof}
The Banach fixed-point argument of
Theorem~\ref{Thm:unique-solution} applies under the discrete
smallness assumptions
(cf. \eqref{small-data-assumption-delta0-discrete}), using
Lemma~\ref{bSh-well-defined}
(cf. \eqref{bound-apriori-bSh},
\eqref{bound-apriori-tbSh-bsi}, and
\eqref{bound-apriori-Xiih-bTih-discrete}).
These uniform estimates also give
\eqref{bound:a-priori-discr}.
\end{proof}
Without the contractivity inequality, invariance still gives
existence in this ball by Brouwer's theorem
\cite[Theorem~9.9-2]{Ciarlet}, without a uniqueness assertion;
the discrete operator $\bPi_h$ is continuous, since it is
Lipschitz by \eqref{LC-bPi-discrete}, and the ball is a compact convex
subset of the finite-dimensional discrete space.

\subsection{A priori error estimates}

Let the continuous and discrete solutions be those of
Theorems~\ref{Thm:unique-solution}
and~\ref{Thm:discrete-unique-solution}.
Writing
$\vec{\bro}:=(\bro,\u,\z,p)$ and
$\vec{\bro}_h:=(\bro_h,\u_h,\z_h,p_h)$,
define the block errors by
\begin{equation}\label{def:block-errors}
\mathtt E_1
:=
\|\vec{\bro}-\vec{\bro}_h\|_{\mathbb X},
\qquad
\mathtt E_2
:=
\|(\bvphi-\bvphi_h,\chi-\chi_h)\|_{\rX_2\times\rM_1},
\qquad
\mathtt E_3
:=
\sum_{i=1}^{2}
\|(\vec{\bt}_i-\vec{\bt}_{i,h},\xi_i-\xi_{i,h})
  \|_{\bcalH\times\calM},
\end{equation}
where all product spaces carry their sum norms.
The total error is
$\mathtt E:=\mathtt E_1+\mathtt E_2+\mathtt E_3$.

For $Z_h\subset Z$, write
$\dist(v,Z_h):=\inf_{v_h\in Z_h}\|v-v_h\|_Z$ and define
the best-approximation errors
\begin{equation}\label{def:block-approximation-errors}
\begin{gathered}
\mathtt D_1
:=
\dist((\bro,\u),\bX_h\times\bQ_h)
+\dist((\z,p),\bX_{2,h}\times\bQ_{1,h}),
\qquad
\mathtt D_2
:=
\dist((\bvphi,\chi),\rX_{2,h}\times\rM_{1,h}),
\\[1ex]
\mathtt D_3
:=
\sum_{i=1}^{2}
\dist((\vec{\bt}_i,\xi_i),\bcalH_h\times\calM_h).
\end{gathered}
\end{equation}

A blockwise Strang argument
(cf. \cite[Lemma~2.27]{eg-SPRINGER-2004})
bounds each error by its best-approximation contribution
and the errors in the fields supplied by the other blocks.
For the electrochemical blocks, we follow
\cite[Section~5.7]{gir-CAMWA-2025}.
For the present mechanical block, we use
\eqref{Biot-F-difference-bound} and the uniform discrete
Biot stability established under
\eqref{c0-alpha-lambda-assumption-disc}.
Together with \eqref{bound:a-priori} and \eqref{bound:a-priori-discr},
these arguments yield
\begin{align}
\mathtt E_1
&\le C_1\Biggl\{
\mathtt D_1
+\Biggl[C_T
\sum_{i=1}^{2}
\bigl(
\|\xi_{i,\tt D}\|_{1/2;\Gamma}
+\|f_i\|_{0,\varrho;\O}
\bigr)+\Lambda_\chi\Biggr]\mathtt E_2
\notag\\
&\qquad+
\Biggl[
1+\widehat{C}_E\Biggl(
\|\chi_{\tt D}\|_{1/s,r;\Gamma}
+\Lambda_\chi
+\widehat{C}_T
\sum_{i=1}^{2}
\bigl(
\|\xi_{i,\tt D}\|_{1/2;\Gamma}
+\|f_i\|_{0,\varrho;\O}
\bigr)\Biggr)
\Biggr]\mathtt E_3
\Biggr\},
\label{ineq-E1-2}
\\[1ex]
\mathtt E_2
&\le C_2(\mathtt D_2+\mathtt E_3),
\label{ineq-E2-1}
\\[1ex]
\mathtt E_3
&\le C_3\Biggl\{
\mathtt D_3
+C_T
\sum_{i=1}^{2}
\bigl(
\|\xi_{i,\tt D}\|_{1/2;\Gamma}
+\|f_i\|_{0,\varrho;\O}
\bigr)(\mathtt E_1+\mathtt E_2)
\Biggr\}.
\label{ineq-E3-2}
\end{align}
The constants $C_1$, $C_2$, and $C_3$ in
\eqref{ineq-E1-2}--\eqref{ineq-E3-2}
are independent of $h$ and the prescribed data under
the block stability hypotheses of
Lemmas~\ref{bS-well-defined} and~\ref{bSh-well-defined}
and the radius conditions
\eqref{assumption-w-bphi} and
\eqref{assumption-w-bphi-disc}.

Adding \eqref{ineq-E1-2} and \eqref{ineq-E2-1},
and substituting \eqref{ineq-E2-1} into the resulting
$\mathtt E_2$ term, gives
\begin{equation}\label{eq:apriori-first-two-blocks}
\mathtt E_1+\mathtt E_2
\le
C_1\mathtt D_1
+\Big(1+C_1\Bigl[C_T\sum_{i=1}^{2}
\bigl(
\|\xi_{i,\tt D}\|_{1/2;\Gamma}
+\|f_i\|_{0,\varrho;\O}
\bigr)+\Lambda_\chi\Bigr]\Big)C_2\mathtt D_2
+L\mathtt E_3,
\end{equation}
where
\begin{equation}\label{def:apriori-coupling-constant}
\begin{array}{l}
\ds L
:=
C_1
\Bigl[
1+\widehat{C}_E
\bigl(\|\chi_{\tt D}\|_{1/s,r;\Gamma}+\Lambda_\chi+\widehat{C}_T\sum_{i=1}^{2}
\bigl(
\|\xi_{i,\tt D}\|_{1/2;\Gamma}
+\|f_i\|_{0,\varrho;\O}
\bigr)\bigr)
\Bigr]\\[2ex]
\ds \qquad 
+C_2\Big(1+C_1\Bigl[C_T\sum_{i=1}^{2}
\bigl(
\|\xi_{i,\tt D}\|_{1/2;\Gamma}
+\|f_i\|_{0,\varrho;\O}
\bigr)+\Lambda_\chi\Bigr]\Big).
\end{array}
\end{equation}

The following data condition permits absorption of the coupling
term after combining \eqref{ineq-E3-2} and
\eqref{eq:apriori-first-two-blocks}.

\begin{theorem}\label{thm:Cea}
Under the hypotheses of
Theorems~\ref{Thm:unique-solution}
and~\ref{Thm:discrete-unique-solution},
assume that
\begin{equation}\label{Cea-assumption}
C_3C_T\sum_{i=1}^{2}
\bigl(
\|\xi_{i,\tt D}\|_{1/2;\Gamma}
+\|f_i\|_{0,\varrho;\O}
\bigr) L\le\dfrac12,
\end{equation}
where $L$ is defined in
\eqref{def:apriori-coupling-constant}.
Then there exists $C_{\rm Cea}>0$, independent of $h$,
such that the errors defined in \eqref{def:block-errors}
and \eqref{def:block-approximation-errors} satisfy
\begin{equation}\label{Cea-estimate}
\mathtt E
\le
C_{\rm Cea}
(\mathtt D_1+\mathtt D_2+\mathtt D_3).
\end{equation}
\end{theorem}

\begin{proof}
Substituting \eqref{eq:apriori-first-two-blocks}
into \eqref{ineq-E3-2} and using \eqref{Cea-assumption}
allows the remaining $\mathtt E_3$ term to be absorbed,
with a factor bounded below by $1/2$.
This bounds $\mathtt E_3$ by
$C(\mathtt D_1+\mathtt D_2+\mathtt D_3)$.
Substitution into \eqref{eq:apriori-first-two-blocks}
proves \eqref{Cea-estimate}.
\end{proof}

\begin{theorem}%\label{th:cv}
Under the hypotheses of Theorem~\ref{thm:Cea},
assume that, for some $l\in[1,\ell+1]$, the continuous
solution satisfies
$\bro\in\mathbb H^l(\O)$,
$\bdiv(\bro),\u\in\bH^l(\O)$,
$\z\in\bW^{l,r}(\O)$,
$\div(\z)\in\rH^l(\O)$,
$p\in\rW^{l,r}(\O)$,
$\bvphi\in\bW^{l,r}(\O)$,
$\div(\bvphi),\chi\in\rW^{l,r}(\O)$,
and, for each $i\in\{1,2\}$,
$\bt_i,\bsi_i\in\bH^l(\O)$,
$\div(\bsi_i)\in\rW^{l,\varrho}(\O)$,
and $\xi_i\in\rW^{l,\rho}(\O)$.
Then there exists $C>0$, independent of $h$, such that
\begin{equation}\label{eq:apriori-convergence-rate}
\begin{aligned}
&\mathtt E
\le C h^l
\Bigl\{
\|\bro\|_{l;\O}
+\|\bdiv(\bro)\|_{l;\O}
+\|\u\|_{l;\O}
+\|\z\|_{l,r;\O}
+\|\div(\z)\|_{l;\O}
+\|p\|_{l,r;\O}
\\
&\quad +
\|\bvphi\|_{l,r;\O}
+\|\div(\bvphi)\|_{l,r;\O}
+\|\chi\|_{l,r;\O}
+
\sum_{i=1}^{2}
\Bigl(
\|\bt_i\|_{l;\O}
+\|\bsi_i\|_{l;\O}
+\|\div(\bsi_i)\|_{l,\varrho;\O}
+\|\xi_i\|_{l,\rho;\O}
\Bigr)
\Bigr\}.
\end{aligned}
\end{equation}
\end{theorem}
\begin{proof}
Apply the approximation estimates of
\cite[Section~6.3]{gir-CAMWA-2025},
\cite[Section~4.1]{gmr-IMANUM-2022},
\cite[Appendix~A]{ccgi-CAMWA-2023}, and
\cite[Proposition~1.135]{eg-SPRINGER-2004}
to \eqref{def:block-approximation-errors}. For the pseudostress,
subtract a constant multiple of $\I$ from its Raviart--Thomas
interpolant to impose zero mean trace without changing its
divergence or approximation order. Substitution into
\eqref{Cea-estimate} proves \eqref{eq:apriori-convergence-rate}.
\end{proof}

\section{A residual-based \textit{a posteriori} error estimator}
\label{sec:posteriori-2D}

We now bound the error through computable residuals of
\eqref{Biot-PNP-fully-mixed-formulation-disc}.
The analysis is restricted to two dimensions and follows the
localization techniques of \cite{cgos-M2AN-2021,ccov2022}.

\subsection{Preliminaries}
\label{subsec:posteriori-preliminaries}

Throughout this section, $\O\subset\mathbb R^2$
is a bounded, simply connected polygonal Lipschitz
domain, and the continuous and discrete solutions
satisfy the hypotheses of
Theorems~\ref{Thm:unique-solution}
and~\ref{Thm:discrete-unique-solution}.

Let $\{\Th\}$ be a shape-regular family of conforming
triangular meshes, with fixed polynomial degree $\ell$.
For each $T\in\Th$ and each edge $e$, we set
$h_T:=\operatorname{diam}T$ and $h_e:=|e|$.
We denote by $\cE_h^{\rm int}$ and $\cE_h^\Gamma$
the sets of interior and boundary edges, respectively.
For each $T$, let $\cE_h(T)$ be its edges, with
$\cE_h^{\rm int}(T):=\cE_h(T)\cap\cE_h^{\rm int}$ and
$\cE_h^\Gamma(T):=\cE_h(T)\cap\cE_h^\Gamma$.
The patches $\omega_T$ and $\omega_e$ consist of
the elements sharing a vertex with $T$ and $e$,
respectively, and have uniformly bounded overlap.

On each edge, fix a unit normal $\bnu_e$ and set
$\boldsymbol\tau_e:=(-(\bnu_e)_2,(\bnu_e)_1)^{\tt t}$.
On boundary edges, $\bnu_e$ is outward.
For an interior edge
$e=\partial T^+\cap\partial T^-$, oriented by the
outward normal of $T^+$, set
$\jump{v}_e:=v|_{T^+}-v|_{T^-}$.
For scalar and vector fields, respectively, we use
\[
\operatorname{curl}(\beta)
:=
(\partial_{x_2}\beta,-\partial_{x_1}\beta)^{\tt t},
\qquad
\operatorname{rot}(\bv)
:=
\partial_{x_1}v_2-\partial_{x_2}v_1,
\qquad
\partial_{\tau_e}v:=\nabla v\cdot\boldsymbol\tau_e.
\]
Tensor rotation is understood row-wise.
Derivatives of discrete fields are taken elementwise.

The exponents satisfy
\begin{equation}\label{eq:posteriori-exponents}
2<r\le4,
\qquad
s=\frac r{r-1},
\qquad
\rho=\frac{2r}{r-2},
\qquad
\varrho=\frac{\rho}{\rho-1},
\qquad
\frac1r+\frac1\rho=\frac12.
\end{equation}

Some indicators compare tangential components with tangential
derivatives of the prescribed Dirichlet data.
To evaluate those boundary terms, we require, in addition to
the global trace assumptions of Section~\ref{section35},
\begin{equation}\label{eq:posteriori-boundary-regularity}
\chi_{\tt D}|_e\in \rW^{1,r}(e),
\qquad
\xi_{i,\tt D}|_e\in \rH^1(e),
\quad i\in\{1,2\},
\qquad e\in\cE_h^\Gamma.
\end{equation}
These edgewise restrictions retain the compatibility
of the original global Dirichlet traces.
Constants in this section may depend on the domain, the physical
coefficients, the exponents, the fixed polynomial
degree, and shape regularity, but are independent of $h$.

\subsubsection{Regular decompositions}

The following decompositions permit local interpolation of the
flux tests. The first two treat the Darcy and Poisson spaces
(cf. \eqref{def-bXi-bQi} and \eqref{def:Poisson-spaces}), with and
without homogeneous normal trace. They combine the divergence
lifting of \cite[Theorem~4.1]{AcostaDuranMuschietti2006} and the
potential theorem of \cite[Theorem~3.1]{Mardare2008}.
The third, from \cite[Lemma~4.4]{ccov2022}, treats the mixed
exponents of the pseudostress and ionic flux tests.

\begin{lemma}
\label{lem:regular-Hpdivp-zero-normal-2D}
Let $1<t<\infty$.
For every $\bv\in\mathbf H_0^t(\div_t;\O)$,
there exist
$\bw\in\mathbf W_0^{1,t}(\O)$ and
$\beta\in \rW_0^{1,t}(\O)$ such that
\begin{equation}\label{eq:regular-Hpdivp-zero-normal-2D}
\bv=\bw+\operatorname{curl}(\beta),
\qquad
\div(\bw)=\div(\bv).
\end{equation}
Moreover,
\begin{equation}\label{eq:regular-Hpdivp-zero-normal-stability}
\|\bw\|_{1,t;\O}
\le C\|\div(\bv)\|_{0,t;\O},
\qquad
\|\bw\|_{1,t;\O}+\|\beta\|_{1,t;\O}
\le C\|\bv\|_{t,\div_t;\O},
\end{equation}
where $C>0$ depends only on $\O$ and $t$.
\end{lemma}

\begin{proof}
Green's formula gives
$\int_\O\div(\bv)=\langle\bv\cdot\bnu,1\rangle_\Gamma=0$,
so $\div(\bv)\in\rL_0^t(\O)$.
Applying the bounded divergence lifting
\cite[Theorem~4.1]{AcostaDuranMuschietti2006} with
$p=t$ and $f=\div(\bv)$ yields
$\bw\in\bW_0^{1,t}(\O)$ such that
\[
\div(\bw)=\div(\bv),
\qquad
\|\bw\|_{1,t;\O}
\le C\|\div(\bv)\|_{0,t;\O}.
\]
The remainder $\ba:=\bv-\bw$ therefore satisfies
$\div(\ba)=0$ in $\O$ and $\ba\cdot\bnu=0$ on $\Gamma$.

Set $\ba^\perp:=(-a_2,a_1)$.
Since $\operatorname{rot}(\ba^\perp)=\div(\ba)=0$
in $\mathcal D'(\O)$ and $\O$ is simply connected,
the Sobolev potential theorem
\cite[Theorem~3.1]{Mardare2008}, applied with $m=0$ and $p=t$,
provides $\beta_0\in\rW^{1,t}(\O)$ with
$\nabla\beta_0=\ba^\perp$, hence $\operatorname{curl}(\beta_0)=\ba$.
Writing $\gamma(\beta_0)$ for its boundary trace,
the normal--tangential trace identity yields
\[
\partial_{\bta_\Gamma}\gamma(\beta_0)
=
\operatorname{curl}(\beta_0)\cdot\bnu
=
\ba\cdot\bnu
=
0
\quad\text{in }\mathcal D'(\Gamma),
\]
where $\bta_\Gamma:=(-n_2,n_1)^{\tt t}$ and
$\bnu=(n_1,n_2)^{\tt t}$.
Since $\Gamma$ is connected, $\gamma(\beta_0)$
is constant, say $c$.
Thus $\beta:=\beta_0-c\in\rW_0^{1,t}(\O)$ and
$\bv=\bw+\operatorname{curl}(\beta)$, proving
\eqref{eq:regular-Hpdivp-zero-normal-2D}.

Finally, Poincar\'e's inequality and the lifting bound give
\[
\begin{aligned}
\|\beta\|_{1,t;\O}
&\le C\|\nabla\beta\|_{0,t;\O}
 = C\|\bv-\bw\|_{0,t;\O}
\le C\bigl(
\|\bv\|_{0,t;\O}
+\|\div(\bv)\|_{0,t;\O}
\bigr).
\end{aligned}
\]
Together with the estimate for $\bw$, this proves
\eqref{eq:regular-Hpdivp-zero-normal-stability}.
The constants depend only on $\O$ and $t$.
\end{proof}

\begin{lemma}
\label{lem:regular-Hpdivp-free-2D}
Let $1<t<\infty$.
For every $\bv\in\mathbf H^t(\div_t;\O)$,
there exist
$\bw\in\mathbf W^{1,t}(\O)$ and
$\beta\in \rW^{1,t}(\O)$, with $\int_\O\beta=0$,
such that
\begin{equation}\label{eq:regular-Hpdivp-free-2D}
\bv=\bw+\operatorname{curl}(\beta),
\qquad
\div(\bw)=\div(\bv).
\end{equation}
Moreover,
\begin{equation}\label{eq:regular-Hpdivp-free-stability}
\|\bw\|_{1,t;\O}
\le C\|\div(\bv)\|_{0,t;\O},
\qquad
\|\bw\|_{1,t;\O}+\|\beta\|_{1,t;\O}
\le C\|\bv\|_{t,\div_t;\O},
\end{equation}
where $C>0$ depends only on $\O$ and $t$.
\end{lemma}

\begin{proof}
Let $m:=|\O|^{-1}\int_\O\div(\bv)$. The divergence lifting
\cite[Theorem~4.1]{AcostaDuranMuschietti2006} of $\div(\bv)-m$
gives $\bw_0\in\mathbf W_0^{1,t}(\O)$ with
$\div(\bw_0)=\div(\bv)-m$ and
$\|\bw_0\|_{1,t;\O}\le C\|\div(\bv)-m\|_{0,t;\O}$.
Then $\bw:=\bw_0+\frac m2\bx$ satisfies $\div(\bw)=\div(\bv)$ and, since
$|m|\le|\O|^{-1/t}\|\div(\bv)\|_{0,t;\O}$,
$\|\bw\|_{1,t;\O}\le C\|\div(\bv)\|_{0,t;\O}$.
The remainder $\bv-\bw$ is divergence free, so the potential argument of
Lemma~\ref{lem:regular-Hpdivp-zero-normal-2D}
(cf. \cite[Theorem~3.1]{Mardare2008}) gives a zero-mean
$\beta\in\rW^{1,t}(\O)$ with $\operatorname{curl}(\beta)=\bv-\bw$,
and Poincar\'e's inequality yields
\eqref{eq:regular-Hpdivp-free-stability}.
\end{proof}

For pseudostress and ionic flux tests, the field lies in $\bL^2$
and its divergence in $\rL^t$; we use the following result.

\begin{lemma}
\label{lem:mixed-regular-decomposition-2D}
Let $1<t\le2$ and
$\bv\in\mathbf L^2(\O)$ with $\div(\bv)\in\rL^t(\O)$.
There exist
$\bw\in\mathbf W^{1,t}(\O)$ and $\beta\in\rH^1(\O)$
such that
\begin{equation}\label{eq:mixed-regular-decomposition-2D}
\bv=\bw+\operatorname{curl}(\beta),
\qquad
\|\bw\|_{1,t;\O}+\|\beta\|_{1;\O}
\le
C\bigl(
\|\bv\|_{0;\O}+\|\div(\bv)\|_{0,t;\O}
\bigr),
\end{equation}
where $C>0$ depends only on $\O$ and $t$.
The result applies row-wise to tensor fields.
\end{lemma}

\begin{proof}
See \cite[Lemma~4.4]{ccov2022}.
\end{proof}

\subsubsection{Interpolation and localization tools}

Let $\Pi_h^\ell$ be the canonical Raviart--Thomas
interpolant and let $I_h$ be a locally stable
Cl\'ement-type quasi-interpolant into continuous
piecewise affine functions.
For $1<t<\infty$, the local estimates are
\begin{align}
\|\bw-\Pi_h^\ell\bw\|_{0,t;T}
&\le
C h_T|\bw|_{1,t;T},
\qquad
\|(\bw-\Pi_h^\ell\bw)\cdot\bnu_e\|_{0,t;e}
\le
C h_e^{1-1/t}|\bw|_{1,t;T},
\notag
\\
\|\beta-I_h\beta\|_{0,t;T}
&\le
C h_T\|\beta\|_{1,t;\omega_T},
\qquad
\|\beta-I_h\beta\|_{0,t;e}
\le
C h_e^{1-1/t}\|\beta\|_{1,t;\omega_e}.
\label{eq:posteriori-interpolation}
\end{align}
For the Raviart--Thomas estimates, see
\cite[Lemmas~4.1 and~4.2]{ccov2022};
for quasi-interpolation, see
\cite[Lemma~4.2]{cgos-M2AN-2021}
and \cite{clement1975}.
All estimates extend componentwise.
Moreover, $\operatorname{curl}(I_h\beta)
\in\mathbf{RT}_0(\Th)
\subseteq\mathbf{RT}_\ell(\Th)$,
and its normal trace vanishes when $I_h\beta$
has zero boundary trace. Thus $\Pi_h^\ell(\bw)+\operatorname{curl}(I_h\beta)$
is a conforming Raviart--Thomas field.
The tensor mean-trace constraint is handled separately in
Section~\ref{subsubsec:residual-localization}.

For the efficiency analysis, we use standard
element- and edge-bubble estimates in local $\rL^t$
spaces, together with polynomial inverse and
trace inequalities.
We refer to
\cite[Section~5.2]{ccov2022},
\cite{Verfurth}, and \cite{eg-SPRINGER-2004}.

\subsection{Local indicators and the global estimator}

We organize the estimator according to the Biot,
Poisson, and Nernst--Planck mixed formulations
\eqref{Biot-fully-mixed},
\eqref{Poisson-mixed-formulation}, and
\eqref{NP-mixed-formulation}, respectively.
Following the residual localization strategy in
\cite{cgos-M2AN-2021,ccov2022},
we use two indicator families for the gradient identities
\eqref{eq:gradient-reconstruction}.
They apply to the displacement, pressure, electrostatic potential,
and concentrations, with the exponents required by each block.

For a piecewise polynomial field $G_h$, a discrete
potential $v_h$, boundary datum $d_{\tt D}$,
and exponents $\varpi,\vartheta>1$, define
\begin{align}
&A_T(G_h,v_h,d_{\tt D};\varpi)^\varpi
:={}
h_T^\varpi\|G_h-\nabla v_h\|_{0,\varpi;T}^\varpi
+
\sum_{e\in\cE_h^{\rm int}(T)}
h_e\|\jump{v_h}_e\|_{0,\varpi;e}^\varpi
+
\sum_{e\in\cE_h^\Gamma(T)}
h_e\|v_h-d_{\tt D}\|_{0,\varpi;e}^\varpi,
\label{def:gradient-indicator-families}
\\[1ex]
&B_T(G_h,d_{\tt D};\vartheta)^\vartheta
:={}
h_T^\vartheta\|\operatorname{rot}(G_h)\|_{0,\vartheta;T}^\vartheta
+
\sum_{e\in\cE_h^{\rm int}(T)}
h_e\|\jump{G_h\boldsymbol{\tau}_e}_e\|_{0,\vartheta;e}^\vartheta
	+
\sum_{e\in\cE_h^\Gamma(T)}
h_e\|G_h\boldsymbol{\tau}_e
      -\partial_{\tau_e}d_{\tt D}\|_{0,\vartheta;e}^\vartheta.\notag
\end{align}
All derivatives are taken elementwise.
For a vector field, $G_h\boldsymbol{\tau}_e$
denotes $G_h\cdot\boldsymbol{\tau}_e$;
for a tensor field, rotation is applied row-wise.
The corresponding potential and boundary datum
are scalar or vector-valued, respectively.
For the exact fields, $G=\nabla v$.
The family $A_T$ measures the defect in this gradient relation
and jumps of the discrete potential, whereas $B_T$ measures
rotational defects and tangential jumps of the gradient field
(cf. \eqref{def:gradient-indicator-families}).
Their boundary terms measure the corresponding mismatches
with the prescribed Dirichlet data.
We denote by $A_T^0(G_h,v_h;\varpi)$ and $B_T^0(G_h;\vartheta)$
the corresponding expressions with the boundary
sums omitted.
These versions apply to the Darcy block
\eqref{Biot-fully-mixed}, whose flux trial and test spaces
in \eqref{def-bXi-bQi} have zero normal trace. For the mechanical gradient identity
\eqref{eq:gradient-reconstruction}, with $d=2$, set
\begin{equation}\label{def:constitutive-localization-fields}
M_h
:=
\frac{1}{\mu}\bro_h^{\tt d}
+
\frac{c_\lambda}{2}\tr(\bro_h)\I.
\end{equation}
The other residuals measure the defects in momentum, mass,
and charge balance and in the ionic constitutive and balance
equations. Evaluating
\eqref{Biot-NPN-mixed-2}, \eqref{Biot-NPN-mixed-4},
\eqref{Biot-NPN-mixed-5-2}, and
\eqref{Biot-NPN-mixed-7}--\eqref{Biot-NPN-mixed-8}
at the discrete solution gives
\begin{subequations}\label{def:posteriori-volume-residuals}
\begin{align}
r_{{\rm mom},h}
:={}&
\ff+\bdiv(\bro_h)
+\frac{\alpha\nu}{\kappa}\z_h
-\varepsilon^{-1}(\xi_{1,h}-\xi_{2,h}+f_\chi)\bvphi_h
-\delta(\bt_{1,h}-\bt_{2,h}),
\label{def:momentum-volume-residual}
\\[1ex]
r_{{\rm mass},h}
:={}&
g-\div(\z_h)-c_0(p_h+\overline p)
-\alpha c_\lambda\tr(\bro_h),
\label{def:mass-volume-residual}
\\[1ex]
r_{{\rm ch},h}
:={}&
\div(\bvphi_h)+\xi_{1,h}-\xi_{2,h}+f_\chi,
\label{def:charge-volume-residual}
\\[1ex]
r_{i,{\rm con},h}
:={}&
\bsi_{i,h}-\kappa_i\bt_{i,h}
-\iota_i\kappa_i\varepsilon^{-1}\xi_{i,h}\bvphi_h
+\xi_{i,h}\z_h,
\label{def:ionic-constitutive-volume-residual}
\\[1ex]
r_{i,{\rm bal},h}
:={}&
f_i+\div(\bsi_{i,h})-\xi_{i,h}
+\xi_{i,h}\div(\z_h),
\qquad i\in\{1,2\}.
\label{def:ionic-balance-volume-residual}
\end{align}
\end{subequations}
Here $\overline p$ is the natural pressure mean of \eqref{def:p-mean}
(or of \eqref{eq:nonhomogeneous-mass-compatibility} in
Section~\ref{subsec:posteriori-nonhomogeneous-boundary}), so that
$p_h+\overline p$ is the full discrete pressure of
Lemma~\ref{lem:pressure-mean} and $r_{{\rm mass},h}$ has zero mean.
The data $g$ and $f_\chi$ enter the residuals directly.
The exponents of the volume residual norms and the
gradient indicators are determined by H\"older's
inequality with the corresponding test spaces
(cf. \eqref{def-bXi-bQi}, \eqref{def:Poisson-spaces},
\eqref{def:bH-H1-H2-M}, and
\eqref{eq:posteriori-exponents}).
The momentum, mass, and ionic constitutive residuals
are measured in $\rL^2$, the charge residual in $\rL^r$,
and the ionic balance residuals in $\rL^\varrho$.
These volume norms carry no meshsize factor because
their functionals act directly on Lebesgue test spaces.
For the gradient indicators, the pseudostress decomposition
uses exponent $2$ for both components, giving $(\varpi,\vartheta)=(2,2)$.
The Darcy and Poisson decompositions use $s$, whose
conjugate exponent is $r$, giving $(\varpi,\vartheta)=(r,r)$
(cf. Lemmas~\ref{lem:regular-Hpdivp-zero-normal-2D}
and~\ref{lem:regular-Hpdivp-free-2D}).
For the ionic flux tests, the regular component belongs
to $\bW^{1,\varrho}(\O)$ and the curl potential to
$\rH^1(\O)$, giving $(\varpi,\vartheta)=(\rho,2)$
(cf. Lemma~\ref{lem:mixed-regular-decomposition-2D}).
The meshsize weights follow from the interpolation
estimates \eqref{eq:posteriori-interpolation}.

Combining the volume residuals and gradient indicators
(cf. \eqref{def:posteriori-volume-residuals}
and \eqref{def:gradient-indicator-families})
with these exponents, we define, for every $T\in\Th$,
\begin{align}
\theta_{{\rm B},2,T}^2
:={}&
\|r_{{\rm mom},h}\|_{0;T}^2
+\|r_{{\rm mass},h}\|_{0;T}^2
+
A_T(M_h,\u_h,0;2)^2
+B_T(M_h,0;2)^2,
\notag
\\[1ex]
\theta_{{\rm B},r,T}^r
:={}&
A_T^0\left(-\frac{\nu}{\kappa}\z_h,p_h;r\right)^r
+
B_T^0\left(-\frac{\nu}{\kappa}\z_h;r\right)^r,
\notag
\\[1ex]
\theta_{{\rm P},r,T}^r
:={}&
\|r_{{\rm ch},h}\|_{0,r;T}^r
+
A_T(\varepsilon^{-1}\bvphi_h,
    \chi_h,\chi_{\tt D};r)^r
+
B_T(\varepsilon^{-1}\bvphi_h,\chi_{\tt D};r)^r,
\label{def:posteriori-local-indicators}
\\[1ex]
\theta_{i,2,T}^2
:={}&
\|r_{i,{\rm con},h}\|_{0;T}^2
+
B_T(\bt_{i,h},\xi_{i,\tt D};2)^2,
\notag
\\[1ex]
\theta_{i,\rho,T}^{\rho}
:={}&
A_T(\bt_{i,h},\xi_{i,h},\xi_{i,\tt D};\rho)^\rho,
\qquad
\theta_{i,\varrho,T}^{\varrho}
:=
\|r_{i,{\rm bal},h}\|_{0,\varrho;T}^{\varrho},
\qquad i\in\{1,2\}.
\notag
\end{align}

Each family is summed over the mesh using its own Lebesgue
exponent. The global estimator is
\begin{align}
\Theta_{\rm B}
:={}&
\left(
\sum_{T\in\Th}\theta_{{\rm B},2,T}^2
\right)^{1/2}
+
\left(
\sum_{T\in\Th}\theta_{{\rm B},r,T}^r
\right)^{1/r},
\qquad
\Theta_{\rm P}
:= 
\left(
\sum_{T\in\Th}\theta_{{\rm P},r,T}^r
\right)^{1/r}, \notag
\\[1ex]
\Theta_{\rm N}
:={}&
\sum_{i=1}^2
\left\{
\left(
\sum_{T\in\Th}\theta_{i,2,T}^2
\right)^{1/2}
+
\left(
\sum_{T\in\Th}\theta_{i,\rho,T}^{\rho}
\right)^{1/\rho}
+
\left(
\sum_{T\in\Th}\theta_{i,\varrho,T}^{\varrho}
\right)^{1/\varrho}
\right\}, \label{def:posteriori-global-estimator}
\\[1ex]
\Theta
:={}&
\Theta_{\rm B}+\Theta_{\rm P}+\Theta_{\rm N}. \notag
\end{align}

Under the stability and smallness hypotheses of
Theorem~\ref{thm:posteriori-reliability-2D}, and the polynomial
boundary-data assumption of
Theorem~\ref{thm:posteriori-efficiency-2D}, we prove
\begin{equation}\label{eq:posteriori-estimator-bounds}
C_{\rm eff}\Theta\le\mathtt E\le C_{\rm rel}\Theta,
\end{equation}
where $C_{\rm eff},C_{\rm rel}>0$ are independent of $h$.
The upper and lower bounds establish reliability and
efficiency, respectively.

\subsection{Reliability}

Evaluating the continuous equations
\eqref{Biot-PNP-fully-mixed-formulation} at the discrete solution
leaves residual functionals on the continuous test spaces.
They vanish on discrete tests by
\eqref{Biot-PNP-fully-mixed-formulation-disc}.
We bound the errors \eqref{def:block-errors} by their dual norms
and coupling terms, then localize these norms using
\eqref{def:posteriori-local-indicators} to prove the upper bound in
\eqref{eq:posteriori-estimator-bounds}.

\subsubsection{Residual functionals and block stability}

For the Biot formulation \eqref{Biot-fully-mixed},
we define
\begin{align}
\mathcal R_1^{\rm B}(\bta)
&:=
-\bfa_s(\bro_h,\bta)-\bfb_s(\bta,\u_h),
&&\forall\,\bta\in\bX,\notag
\\[1ex]
\mathcal R_2^{\rm B}(\bv)
&:=
\bF_{\bvphi_h,\bxi_h,\bt_h}(\bv)
-\bfb_s(\bro_h,\bv)-\bfe_s(\z_h,\bv),
&&\forall\,\bv\in\bQ,
\label{def:Biot-residual-functionals}
\\[1ex]
\mathcal R_3^{\rm B}(\bw)
&:=
-\bfa_f(\z_h,\bw)-\bfd_1(\bw,p_h),
&&\forall\,\bw\in\bX_1,
\notag
\\[1ex]
\mathcal R_4^{\rm B}(q)
&:=
\bG(q)-\bfd_2(\z_h,q)+\bfe_f((\bro_h,p_h),q),
&&\forall\,q\in\bQ_2.
\notag
\end{align}
By \eqref{Biot-PNP-fully-mixed-formulation-disc}, each functional
vanishes on its corresponding discrete test space.

\begin{lemma}\label{lemma:abstract-Biot-residual}
Assume the continuous and discrete well-posedness
hypotheses of
Theorems~\ref{Thm:unique-solution}
and~\ref{Thm:discrete-unique-solution}.
Then there exists $\widetilde{C}_{\rm B}>0$,
independent of $h$, such that the block errors
defined in \eqref{def:block-errors} satisfy
\begin{align*}
&\mathtt E_1
\le \widetilde{C}_{\rm B}\Biggl\{
\|\mathcal R_1^{\rm B}\|_{\bX'}
+\|\mathcal R_2^{\rm B}\|_{\bQ'}
+\|\mathcal R_3^{\rm B}\|_{\bX_1'}
+\|\mathcal R_4^{\rm B}\|_{\bQ_2'}
+
\varepsilon^{-1}
\Biggl[
\sum_{i=1}^{2}
\bigl(
\|\xi_{i,\tt D}\|_{1/2;\Gamma}
+\|f_i\|_{0,\varrho;\O}
\bigr)+\Lambda_\chi\Biggr]\mathtt E_2\\
&\qquad
+
\Biggl[
\varepsilon^{-1}\Biggl(
\|\chi_{\tt D}\|_{1/s,r;\Gamma}
+\Lambda_\chi
+
\sum_{i=1}^{2}
\bigl(
\|\xi_{i,\tt D}\|_{1/2;\Gamma}
+\|f_i\|_{0,\varrho;\O}
\bigr)\Biggr)
+\delta
\Biggr]\mathtt E_3
\Biggr\}.
\end{align*}
\end{lemma}
\begin{proof}
Since the approximation is conforming, the Biot error
belongs to $\mathbb X$.
The continuous equations \eqref{Biot-fully-mixed}
and the residual definitions
\eqref{def:Biot-residual-functionals} give
\[
\begin{aligned}
\bfa_s(\bro-\bro_h,\bta)
+\bfb_s(\bta,\u-\u_h)
&=\mathcal R_1^{\rm B}(\bta),
\\[1ex]
\bfb_s(\bro-\bro_h,\bv)
+\bfe_s(\z-\z_h,\bv)
&=\mathcal R_2^{\rm B}(\bv)
+
\bigl(
\bF_{\bvphi,\bxi,\bt}
-\bF_{\bvphi_h,\bxi_h,\bt_h}
\bigr)(\bv),
\\[1ex]
\bfa_f(\z-\z_h,\bw)
+\bfd_1(\bw,p-p_h)
&=\mathcal R_3^{\rm B}(\bw),
\\[1ex]
\bfd_2(\z-\z_h,q)
-\bfe_f((\bro-\bro_h,p-p_h),q)
&=\mathcal R_4^{\rm B}(q),
\end{aligned}
\]
for every $(\bta,\bv,\bw,q)\in\mathbb Q$.
Applying \eqref{global-inf-sup-bfA+es-ef-1}
to these error equations and estimating the load
difference by \eqref{Biot-F-difference-bound} yields
\[
\begin{aligned}
\frac{\balph_{\bfA}}{2}\mathtt E_1
\le{}&
\|\mathcal R_1^{\rm B}\|_{\bX'}
+\|\mathcal R_2^{\rm B}\|_{\bQ'}
+\|\mathcal R_3^{\rm B}\|_{\bX_1'}
+\|\mathcal R_4^{\rm B}\|_{\bQ_2'}
\\
&+
\varepsilon^{-1}\bigl(\|\bxi\|_{\bcalM}+\Lambda_\chi\bigr)
\|\bvphi-\bvphi_h\|_{\rX_2}
+
\varepsilon^{-1}\|\bvphi_h\|_{\rX_2}
\|\bxi-\bxi_h\|_{\bcalM}
+
\delta\|\bt-\bt_h\|_{\bcalH_1}.
\end{aligned}
\]
By \eqref{def:block-errors}, the field error is at most
$\mathtt E_2$ and each ionic error at most $\mathtt E_3$.
Use \eqref{bound:a-priori} and \eqref{bound:a-priori-discr} for
$\|\bxi\|_{\bcalM}$ and $\|\bvphi_h\|_{\rX_2}$, respectively,
keep $\Lambda_\chi$ explicit, and absorb the constants into
\[
\widetilde{C}_{\rm B}
:=
\frac{2}{\balph_{\bfA}}
\max\bigl\{
1,\,
C_T,\,
\widehat{C}_E,\,
\widehat{C}_E\widehat{C}_T
\bigr\}.
\]
The uniform bounds of Lemmas~\ref{bS-well-defined}
and~\ref{bSh-well-defined} ensure independence of $h$.
\end{proof}

The volume definitions \eqref{def:momentum-volume-residual} and
\eqref{def:mass-volume-residual} give
$\mathcal R_2^{\rm B}(\bv)=-\int_\O r_{{\rm mom},h}\cdot\bv$
and $\mathcal R_4^{\rm B}(q)=-\int_\O r_{{\rm mass},h}q$, hence
\begin{equation}\label{eq:Biot-volume-residual-norms}
\|\mathcal R_2^{\rm B}\|_{\bQ'}
=
\|r_{{\rm mom},h}\|_{0;\O},
\qquad
\|\mathcal R_4^{\rm B}\|_{\bQ_2'}
\le
\|r_{{\rm mass},h}\|_{0;\O}.
\end{equation}
The gradient residuals are localized below using the regular
decompositions and \eqref{eq:posteriori-interpolation}.

\subsubsection{Electrostatic and ionic residuals}

For the Poisson formulation \eqref{Poisson-mixed-formulation}, define
\begin{align}
\mathcal R_1^{\rm P}(\bpsi)
&:=
G(\bpsi)-a(\bvphi_h,\bpsi)-b_1(\bpsi,\chi_h),
&&\forall\,\bpsi\in\rX_1,
\label{def:Poisson-residual-functionals}
\\[1ex]
\mathcal R_2^{\rm P}(\gamma)
&:=
F_{\bxi_h}(\gamma)-b_2(\bvphi_h,\gamma),
&&\forall\,\gamma\in\rM_2.
\notag
\end{align}
These functionals vanish on
$\rX_{1,h}$ and $\rM_{2,h}$, respectively.

For later use, we collect the dual residual norms into
one quantity for each block:
\begin{equation}\label{def:abstract-residual-sums}
\begin{aligned}
\mathcal R_{\rm B}
:=
\|\mathcal R_1^{\rm B}\|_{\bX'}
+\|\mathcal R_2^{\rm B}\|_{\bQ'}
+\|\mathcal R_3^{\rm B}\|_{\bX_1'}
+\|\mathcal R_4^{\rm B}\|_{\bQ_2'}
\qan
\mathcal R_{\rm P}
:=
\|\mathcal R_1^{\rm P}\|_{\rX_1'}
+\|\mathcal R_2^{\rm P}\|_{\rM_2'}.
\end{aligned}
\end{equation}
The charge residual (cf. \eqref{def:charge-volume-residual}) gives
\begin{equation}\label{Poisson-balance-residual-norm}
\|\mathcal R_2^{\rm P}\|_{\rM_2'}
=
\|r_{{\rm ch},h}\|_{0,r;\O}.
\end{equation}

Separating the flux, gradient, and balance relations in
\eqref{NP-mixed-formulation} gives, for each species,
\begin{align}
\mathcal R_{i,1}^{\rm N}(\bs)
&:=
\int_\O r_{i,{\rm con},h}\cdot\bs,
&&\forall\,\bs\in\calH_1, \notag
\\[1ex]
\mathcal R_{i,2}^{\rm N}(\bta)
&:=
\langle\bta\cdot\bnu,\xi_{i,\tt D}\rangle_\Gamma
-\int_\O\bt_{i,h}\cdot\bta
-\int_\O\xi_{i,h}\div(\bta),
&&\forall\,\bta\in\calH_2, \label{def:NP-residual-functionals}
\\[1ex]
\mathcal R_{i,3}^{\rm N}(\eta)
&:=
-\int_\O r_{i,{\rm bal},h}\,\eta,
&&\forall\,\eta\in\calM.\notag
\end{align}
Each functional vanishes on its corresponding
discrete test space.
Set
\begin{equation}\label{def:abstract-NP-residual-sum}
{\mathcal R}_{\rm N}
:=
\sum_{i=1}^2
\left\{
\|\mathcal R_{i,1}^{\rm N}\|_{\calH_1'}
+
\|\mathcal R_{i,2}^{\rm N}\|_{\calH_2'}
+
\|\mathcal R_{i,3}^{\rm N}\|_{\calM'}
\right\}.
\end{equation}
The first and third contributions are the volume norms of the
ionic constitutive and balance residuals
(cf. \eqref{def:ionic-constitutive-volume-residual} and
\eqref{def:ionic-balance-volume-residual}):
\begin{equation}\label{NP-volume-residual-norms}
\|\mathcal R_{i,1}^{\rm N}\|_{\calH_1'}
=
\|r_{i,{\rm con},h}\|_{0;\O}
\qan
\|\mathcal R_{i,3}^{\rm N}\|_{\calM'}
=
\|r_{i,{\rm bal},h}\|_{0,\varrho;\O}.
\end{equation}

\begin{lemma}\label{lemma:abstract-PNP-residual}
Assume the continuous and discrete well-posedness
hypotheses of
Theorems~\ref{Thm:unique-solution}
and~\ref{Thm:discrete-unique-solution}.
Then there exist positive constants
$\widetilde{C}_{\rm P}$ and $\widetilde{C}_{\rm N}$,
independent of $h$, such that the block errors
defined in \eqref{def:block-errors} and the residual
sums defined in \eqref{def:abstract-residual-sums}
and \eqref{def:abstract-NP-residual-sum} satisfy
\begin{subequations}\label{abstract-PNP-error-bounds}
\begin{align}
\mathtt E_2
&\le
\widetilde{C}_{\rm P}
\bigl(
\mathcal R_{\rm P}+\mathtt E_3
\bigr),
\label{eq:abstract-Poisson-error-bound}
\\[1ex]
\mathtt E_3
&\le
\widetilde{C}_{\rm N}\Biggl\{
\mathcal R_{\rm N}
+
\sum_{i=1}^{2}
\bigl(
\|\xi_{i,\tt D}\|_{1/2;\Gamma}
+\|f_i\|_{0,\varrho;\O}
\bigr)
(\mathtt E_1+\mathtt E_2)
\Biggr\}.
\label{eq:abstract-ionic-error-bound}
\end{align}
\end{subequations}
\end{lemma}

\begin{proof}
By \eqref{Poisson-mixed-formulation},
\eqref{def:Poisson-residual-functionals}, and
\eqref{Poisson-global-stability}, followed by \eqref{bound-F-xi}
and the sum norms in \eqref{def:block-errors},
\[
\mathtt E_2\le C_{\rm P}\bigl\{
\mathcal R_{\rm P}+\|F_{\bxi}-F_{\bxi_h}\|_{\rM_2'}\bigr\}
\le C_{\rm P}\bigl(\mathcal R_{\rm P}+\|F\|\mathtt E_3\bigr).
\]
Thus \eqref{eq:abstract-Poisson-error-bound}
holds with
$\widetilde{C}_{\rm P}
:=C_{\rm P}\max\{1,\|F\|\}$.

For the ionic equations, we use the operator
$\boldsymbol\calA_{\z,\bvphi}$
(cf. \eqref{def:ionic-block-form})
with the exact coefficient pair
$(\z,\bvphi)\in\rW(\mathtt R)$.
For every
$(\vec{\bs}_i,\eta_i)\in\bcalH\times\calM$,
with $\vec{\bs}_i=(\bs_i,\bta_i)$,
the continuous equations \eqref{NP-mixed-formulation}
and the residual definitions
\eqref{def:NP-residual-functionals} yield
\begin{align*}
&\boldsymbol\calA_{\z,\bvphi}
\bigl(
(\vec{\bt}_i-\vec{\bt}_{i,h},\xi_i-\xi_{i,h}),
(\vec{\bs}_i,\eta_i)
\bigr)
=
\mathcal R_{i,1}^{\rm N}(\bs_i)
+
\mathcal R_{i,2}^{\rm N}(\bta_i)
+
\mathcal R_{i,3}^{\rm N}(\eta_i)
\\
&\qquad+
\bigl(
\calE_{\z_h,\bvphi_h}-\calE_{\z,\bvphi}
\bigr)(\vec{\bs}_i,\xi_{i,h})
+
\bigl(
\calD_{\z_h}-\calD_{\z}
\bigr)(\xi_{i,h},\eta_i).
\end{align*}
The linear dependence of $\calE$ and $\calD$
on their coefficient fields, together with
\eqref{bound-E-D} and the definition of $D_*$
(cf. \eqref{def:ionic-coupling-bound}),
bounds the last two terms by
\[
D_*\|\xi_{i,h}\|_{\calM}
\left(
\|\z-\z_h\|_{\bX_2}
+
\|\bvphi-\bvphi_h\|_{\rX_2}
\right)
\|(\vec{\bs}_i,\eta_i)\|_{\bcalH\times\calM}.
\]
Applying \eqref{global-inf-sup-Awbze-1},
summing over the two species, and using
\eqref{def:block-errors} gives
\[
\mathtt E_3
\le
\frac{2}{\balph_{\boldsymbol\calA}}\Biggl\{
\mathcal R_{\rm N}
+
D_*\sum_{i=1}^{2}\|\xi_{i,h}\|_{\calM}
(\mathtt E_1+\mathtt E_2)
\Biggr\}.
\]
Using the discrete a priori estimate
\eqref{bound:a-priori-discr} for the concentrations
and absorbing the coefficient into a common factor
proves \eqref{eq:abstract-ionic-error-bound} with
$
\widetilde{C}_{\rm N}
:=
\frac{2}{\balph_{\boldsymbol\calA}}
\max\{1,D_*\widehat{C}_T\}.
$
\end{proof}
The remaining dual norms,
$\|\mathcal R_1^{\rm P}\|_{\rX_1'}$ and
$\|\mathcal R_{i,2}^{\rm N}\|_{\calH_2'}$,
will be localized using
Lemmas~\ref{lem:regular-Hpdivp-free-2D}
and~\ref{lem:mixed-regular-decomposition-2D},
respectively, together with the interpolation
estimates \eqref{eq:posteriori-interpolation}.

\subsubsection{Localization of the dual residual norms}
\label{subsubsec:residual-localization}

The gradient residuals share the following localization bound,
expressed through $A_T$ and $B_T$
(cf. \eqref{def:gradient-indicator-families}).

\begin{lemma}\label{lem:gradient-residual-localization}
Let
\[
Y:=
\bigl\{
\bv\in\mathbf L^a(\O):\quad
\div(\bv)\in\rL^b(\O)
\bigr\},
\qquad
\|\bv\|_Y
:=
\|\bv\|_{0,a;\O}
+
\|\div(\bv)\|_{0,b;\O},
\]
where $(a,b)\in\bigl\{(2,2),(s,s),(2,\varrho)\bigr\}$, 
$\varpi:=b'$, and $\vartheta:=a'$. Here, primes denote conjugate exponents.
Consider the functional
\begin{equation}\label{def:generic-gradient-residual}
\mathcal R(\bv)
:=
\langle\bv\cdot\bnu,d_{\tt D}\rangle_\Gamma
-\int_\O G_h\cdot\bv
-\int_\O v_h\div(\bv),
\end{equation}
where $G_h$ and $v_h$ are piecewise polynomial
and $d_{\tt D}$ satisfies the corresponding
boundary assumptions in
Section~\ref{subsec:posteriori-preliminaries}.
Assume that $\mathcal R$ vanishes on the conforming
Raviart--Thomas test space of order $\ell$.
Then
\begin{equation}\label{eq:gradient-residual-localization}
\begin{aligned}
\|\mathcal R\|_{Y'}
\le C\Bigg\{&
\left(
\sum_{T\in\Th}
A_T(G_h,v_h,d_{\tt D};\varpi)^\varpi
\right)^{1/\varpi}
+
\left(
\sum_{T\in\Th}
B_T(G_h,d_{\tt D};\vartheta)^\vartheta
\right)^{1/\vartheta}
\Bigg\},
\end{aligned}
\end{equation}
where $C>0$ is independent of $h$.

The estimate also holds for
$Y=\mathbf H_0^s(\div_s;\O)$, with the boundary
functional omitted and $A_T,B_T$ replaced by
$A_T^0,B_T^0$, respectively.
The tensor version holds row-wise.
\end{lemma}

\begin{proof}
For $(a,b)=(2,2)$ and $(2,\varrho)$, apply
Lemma~\ref{lem:mixed-regular-decomposition-2D}
(cf. \eqref{eq:mixed-regular-decomposition-2D}).
For $(a,b)=(s,s)$, apply
Lemma~\ref{lem:regular-Hpdivp-free-2D},
or Lemma~\ref{lem:regular-Hpdivp-zero-normal-2D}
in the homogeneous-normal case.
These results provide
\[
\bv=\bw+\operatorname{curl}(\beta),
\qquad
\bw\in\mathbf W^{1,b}(\O),
\qquad
\beta\in \rW^{1,a}(\O),
\]
with $\|\bw\|_{1,b;\O}+\|\beta\|_{1,a;\O}\le C\|\bv\|_Y$.
In the homogeneous-normal case,
$\bw\in\mathbf W_0^{1,s}(\O)$ and
$\beta\in \rW_0^{1,s}(\O)$.

The field $\bv_h:=\Pi_h^\ell\bw+\operatorname{curl}(I_h\beta)$
is an admissible discrete test, with zero normal trace in the
homogeneous case. Galerkin orthogonality gives
\[
\mathcal R(\bv)
=
\mathcal R(\bw-\Pi_h^\ell\bw)
+
\mathcal R\bigl(
\operatorname{curl}(\beta-I_h\beta)
\bigr).
\]

Elementwise integration by parts represents the
first term through the volume defect
$G_h-\nabla v_h$, the interior jumps of $v_h$,
and the boundary mismatch $v_h-d_{\tt D}$.
H\"older's inequality with exponents $\varpi=b'$ and $b$,
the Raviart--Thomas estimates in
\eqref{eq:posteriori-interpolation},
and bounded mesh overlap yield
\[
\bigl|
\mathcal R(\bw-\Pi_h^\ell\bw)
\bigr|
\le
C
\left(
\sum_{T\in\Th}
A_T(G_h,v_h,d_{\tt D};\varpi)^\varpi
\right)^{1/\varpi}
\|\bw\|_{1,b;\O}.
\]

For the curl component, integration by parts produces
$\operatorname{rot}G_h$, tangential jumps, and
$G_h\boldsymbol\tau_e-\partial_{\tau_e}d_{\tt D}$ on the boundary.
The boundary terms use the trace assumptions of
Section~\ref{subsec:posteriori-preliminaries}. H\"older's inequality
with $\vartheta=a'$ and the Cl\'ement estimates
\eqref{eq:posteriori-interpolation} give
\[
\left|
\mathcal R\bigl(
\operatorname{curl}(\beta-I_h\beta)
\bigr)
\right|
\le
C
\left(
\sum_{T\in\Th}
B_T(G_h,d_{\tt D};\vartheta)^\vartheta
\right)^{1/\vartheta}
\|\beta\|_{1,a;\O}.
\]
Combining these bounds with the stability of the
decomposition proves
\eqref{eq:gradient-residual-localization}.
In the homogeneous-normal case, both boundary
contributions vanish.
The tensor estimate follows by applying the
argument to each row.
\end{proof}

To apply Lemma~\ref{lem:gradient-residual-localization} row-wise,
extend $\mathcal R_1^{\rm B}$ from $\bX$ to $\bbH(\bdiv;\O)$
using \eqref{def:Biot-residual-functionals}.
The mean-trace constraint gives $\mathcal R_1^{\rm B}(\I)=0$.
Set $c_h:=(2|\O|)^{-1}\int_\O\tr(\bta_h)$ for any tensor
Raviart--Thomas field $\bta_h$. Then $\bta_h-c_h\I\in\bX_h$, so
\begin{equation}\label{eq:mechanical-extended-orthogonality}
\mathcal R_1^{\rm B}(\bta_h)
=
\mathcal R_1^{\rm B}(\bta_h-c_h\I)
+
c_h\mathcal R_1^{\rm B}(\I)
=
0.
\end{equation}
Thus the extended residual satisfies the required orthogonality,
and its estimate also bounds the dual norm on $\bX$.

The choices for the four gradient residuals are summarized in
\eqref{eq:residual-localization-choices}:
\begin{equation}\label{eq:residual-localization-choices}
\begin{array}{c|c|c|c|c}
\text{Residual}&G_h&v_h&d_{\tt D}&(\varpi,\vartheta)
\\ \hline
\mathcal R_1^{\rm B}&M_h&\u_h&0&(2,2)
\\[1ex]
\mathcal R_3^{\rm B}&-\nu\kappa^{-1}\z_h&p_h&
\text{not used}&(r,r)
\\[1ex]
\mathcal R_1^{\rm P}&\varepsilon^{-1}\bvphi_h&\chi_h&
\chi_{\tt D}&(r,r)
\\[1ex]
\mathcal R_{i,2}^{\rm N}&\bt_{i,h}&\xi_{i,h}&
\xi_{i,\tt D}&(\rho,2)
\end{array}
\end{equation}
For $\mathcal R_3^{\rm B}$ we use the homogeneous-normal
version of Lemma~\ref{lem:gradient-residual-localization}.
This functional is the negative of
\eqref{def:generic-gradient-residual},
which leaves its dual norm unchanged.
For $\mathcal R_1^{\rm B}$ the estimate is applied row-wise,
using the orthogonality of the extended mechanical residual
established in \eqref{eq:mechanical-extended-orthogonality}.

Combining \eqref{eq:gradient-residual-localization}
with \eqref{eq:Biot-volume-residual-norms},
\eqref{Poisson-balance-residual-norm},
\eqref{NP-volume-residual-norms}, and
\eqref{def:posteriori-local-indicators}--
\eqref{def:posteriori-global-estimator}, we obtain
\begin{equation}\label{eq:localized-residual-sums}
\mathcal R_{\rm B}\le C_{\rm loc}\Theta_{\rm B},
\qquad
\mathcal R_{\rm P}\le C_{\rm loc}\Theta_{\rm P},
\qquad
\mathcal R_{\rm N}\le C_{\rm loc}\Theta_{\rm N},
\end{equation}
where $C_{\rm loc}>0$ is independent of $h$.

\subsubsection{Closure of the nonlinear coupling}

We combine Lemma~\ref{lemma:abstract-Biot-residual}
with the Poisson and ionic estimates of
Lemma~\ref{lemma:abstract-PNP-residual}
(cf. \eqref{abstract-PNP-error-bounds}).
The ionic bound \eqref{eq:abstract-ionic-error-bound}
provides the data-dependent coupling coefficient
\begin{equation}\label{def:coupled-residual-constants}
d_*
:=
\widetilde{C}_{\rm N}
\sum_{i=1}^{2}
\bigl(
\|\xi_{i,\tt D}\|_{1/2;\Gamma}
+\|f_i\|_{0,\varrho;\O}
\bigr).
\end{equation}
Adding the Biot bound of Lemma~\ref{lemma:abstract-Biot-residual}
to \eqref{eq:abstract-Poisson-error-bound} and using the latter
to eliminate $\mathtt E_2$ gives
\begin{equation}\label{eq:reliability-first-two-blocks}
\begin{aligned}
\mathtt E_1+\mathtt E_2
\le{}&
\widetilde{C}_{\rm BP}\Biggl\{
\mathcal R_{\rm B}
+
\Biggl[
1+\Lambda_\chi+
\sum_{i=1}^{2}
\bigl(
\|\xi_{i,\tt D}\|_{1/2;\Gamma}
+\|f_i\|_{0,\varrho;\O}
\bigr)
\Biggr]\mathcal R_{\rm P}
\Biggr\}
\\[1ex]
&+
\widetilde{C}_{\rm BP}\Biggl\{
1+\varepsilon^{-1}
\|\chi_{\tt D}\|_{1/s,r;\Gamma}+\delta+\Lambda_\chi
+
\sum_{i=1}^{2}
\bigl(
\|\xi_{i,\tt D}\|_{1/2;\Gamma}
+\|f_i\|_{0,\varrho;\O}
\bigr)
\Biggr\}\mathtt E_3.
\end{aligned}
\end{equation}
The constant $\widetilde{C}_{\rm BP}>0$ depends only on
$\widetilde{C}_{\rm B}$, $\widetilde{C}_{\rm P}$,
and $\varepsilon$, independently of $h$ and the
prescribed data.

\begin{theorem}[Reliability in two dimensions]
\label{thm:posteriori-reliability-2D}
Assume the continuous and discrete well-posedness
hypotheses of
Theorems~\ref{Thm:unique-solution}
and~\ref{Thm:discrete-unique-solution},
and the domain, mesh, and boundary assumptions
of Section~\ref{subsec:posteriori-preliminaries}.
There exists a constant $C_{\rm coup}>0$,
independent of $h$ and the prescribed data,
such that, if
\begin{equation}\label{eq:reliability-small-data-condition}
\begin{aligned}
&C_{\rm coup}
\sum_{i=1}^{2}
\bigl(
\|\xi_{i,\tt D}\|_{1/2;\Gamma}
+\|f_i\|_{0,\varrho;\O}
\bigr)
\Biggl\{
1+\varepsilon^{-1}\|\chi_{\tt D}\|_{1/s,r;\Gamma}
+\delta+\Lambda_\chi
+
\sum_{i=1}^{2}
\bigl(
\|\xi_{i,\tt D}\|_{1/2;\Gamma}
+\|f_i\|_{0,\varrho;\O}
\bigr)
\Biggr\}
\le\dfrac12,
\end{aligned}
\end{equation}
then there exists $C_{\rm rel}>0$, independent of $h$,
such that
\begin{equation}\label{eq:posteriori-reliability-2D}
\mathtt E\le C_{\rm rel}\Theta.
\end{equation}
\end{theorem}

\begin{proof}
Substituting \eqref{eq:reliability-first-two-blocks}
into the ionic bound \eqref{eq:abstract-ionic-error-bound},
and using \eqref{def:coupled-residual-constants},
bounds the remaining coupling coefficient by the
left-hand side of
\eqref{eq:reliability-small-data-condition}.
Here, $C_{\rm coup}>0$ depends only on
$\widetilde{C}_{\rm BP}$ and $\widetilde{C}_{\rm N}$,
independently of $h$ and the prescribed data.
Thus, \eqref{eq:reliability-small-data-condition}
allows the $\mathtt E_3$ term to be absorbed,
with a factor bounded below by $1/2$.
Returning to \eqref{eq:reliability-first-two-blocks}
therefore gives
\begin{equation}\label{abstract-residual-global-bound}
\mathtt E
\le
C_{\rm abs}
\bigl(
\mathcal R_{\rm B}
+\mathcal R_{\rm P}
+\mathcal R_{\rm N}
\bigr),
\end{equation}
where $C_{\rm abs}>0$ is independent of $h$.
Combining \eqref{abstract-residual-global-bound}
with the localized residual bounds
\eqref{eq:localized-residual-sums}
and the estimator definition
\eqref{def:posteriori-global-estimator}
proves \eqref{eq:posteriori-reliability-2D}
with $C_{\rm rel}:=C_{\rm abs}C_{\rm loc}$.
\end{proof}

\subsection{Efficiency}

For efficiency, gradient indicators require local bubble arguments,
whereas volume residuals follow directly from the equations.
Assume that $\chi_{\tt D}$ and $\xi_{i,\tt D}$ restrict to polynomials
of uniformly bounded degree on boundary edges. Using the bubble,
inverse, and lifting estimates of \cite[Section~5.2]{ccov2022}, the
following lemma includes the mixed exponents needed for ionic gradients.

\begin{lemma}[Local bounds for gradient indicators]
\label{lem:gradient-indicator-efficiency}
Let $G=\nabla v$ in the distributional sense,
with $G\in\bL^\vartheta(\O)$ and $v\in \rL^\varpi(\O)$,
and let $G_h$ and $v_h$ be piecewise polynomials
of uniformly bounded degree.
For $(\varpi,\vartheta)\in\{(2,2),(r,r),(\rho,2)\}$, there exists $C>0$, independent of $h$, such that
\begin{align}
h_T\|G_h-\nabla v_h\|_{0,\varpi;T}
&\le
C\left(
\|v-v_h\|_{0,\varpi;T}
+
h_T^{1+2/\varpi-2/\vartheta}\|G-G_h\|_{0,\vartheta;T}
\right), \notag
\\[1ex]
h_T\|\operatorname{rot}(G_h)\|_{0,\vartheta;T}
&\le
C\|G-G_h\|_{0,\vartheta;T}, \label{eq:local-gradient-efficiency}
\\[1ex]
h_e^{1/\varpi}\|\jump{v_h}_e\|_{0,\varpi;e}
&\le
C\left(
\|v-v_h\|_{0,\varpi;\omega_e}
+
h_e^{1+2/\varpi-2/\vartheta}
\|G-G_h\|_{0,\vartheta;\omega_e}
\right), \notag
\\[1ex]
h_e^{1/\vartheta}
\|\jump{G_h\boldsymbol{\tau}_e}_e\|_{0,\vartheta;e}
&\le
C\|G-G_h\|_{0,\vartheta;\omega_e},
\qquad e\in\cE_h^{\rm int}.\notag
\end{align}
If $v$ has boundary trace $d_{\tt D}$, whose
restriction to each boundary edge is a polynomial
of uniformly bounded degree, then
\[
\begin{aligned}
h_e^{1/\varpi}\|v_h-d_{\tt D}\|_{0,\varpi;e}
&\le
C\left(
\|v-v_h\|_{0,\varpi;T}
+
h_e^{1+2/\varpi-2/\vartheta}\|G-G_h\|_{0,\vartheta;T}
\right),
\\[1ex]
h_e^{1/\vartheta}
\|G_h\boldsymbol{\tau}_e
-\partial_{\boldsymbol{\tau}_e}d_{\tt D}\|_{0,\vartheta;e}
&\le
C\|G-G_h\|_{0,\vartheta;T},
\end{aligned}
\]
where $e\subset\partial T\cap\Gamma$.
The estimates also apply row-wise to tensor fields.
The constant $C$ depends only on the exponents,
the polynomial degree bounds, and mesh shape regularity.
\end{lemma}

\begin{proof}
Set $R_T:=G_h-\nabla v_h$.
Since
$
R_T=(G_h-G)+\nabla(v-v_h),
$
testing with an element-bubble vector field
$\bw_T$ and integrating by parts gives
\[
\int_T R_T\cdot\bw_T
=
\int_T(G_h-G)\cdot\bw_T
-
\int_T(v-v_h)\div(\bw_T).
\]
Polynomial norm equivalence and bubble-weighted
duality reduce the estimate of
$\|R_T\|_{0,\varpi;T}$ to such pairings.
Moreover, the scaled inverse inequalities yield
\[
\begin{aligned}
\|\div(\bw_T)\|_{0,\varpi';T}
\le Ch_T^{-1}\|\bw_T\|_{0,\varpi';T}
\qan
\|\bw_T\|_{0,\vartheta';T}
\le
Ch_T^{2/\varpi-2/\vartheta}\|\bw_T\|_{0,\varpi';T},
\end{aligned}
\]
where $\varpi'$ and $\vartheta'$ are the conjugate exponents.
H\"older's inequality therefore proves the
first bound in
\eqref{eq:local-gradient-efficiency}.

For the rotation residual, we use
$\operatorname{rot}(G)=0$ and an element-bubble
scalar test function $\psi_T$ to obtain
\[
\int_T\operatorname{rot}(G_h)\,\psi_T
=
\int_T(G_h-G)\cdot\operatorname{curl}(\psi_T).
\]
Bubble-weighted duality and the inverse estimate for $\nabla\psi_T$ prove the rotation bound
in \eqref{eq:local-gradient-efficiency}.

The interior-edge estimates use edge-bubble liftings of the
potential and tangential jumps, integration by parts on adjacent
elements, and the volume bounds, as in
\cite[Section~5.2]{ccov2022}. Changing the volume exponent gives
the factor $h_e^{2/\varpi-2/\vartheta}$ for the potential jump.

On boundary edges, use $v=d_{\tt D}$ and
$G\boldsymbol\tau_e=\partial_{\boldsymbol\tau_e}d_{\tt D}$
in the trace sense. Polynomial data permit the same bubble
duality. Applying the argument row-wise gives the tensor version.
\end{proof}

Applying these local estimates to the choices in
\eqref{eq:residual-localization-choices}, together with the balance
equations, gives the efficiency bound.

\begin{theorem}
\label{thm:posteriori-efficiency-2D}
Assume the domain, mesh, and boundary hypotheses
of Section~\ref{subsec:posteriori-preliminaries},
and the polynomial boundary-data assumption of
Lemma~\ref{lem:gradient-indicator-efficiency}.
Assume also the continuous and discrete well-posedness
hypotheses of
Theorems~\ref{Thm:unique-solution}
and~\ref{Thm:discrete-unique-solution}.
Then there exists $C_{\rm eff}>0$, depending only on
the physical parameters, the domain, the mesh
shape-regularity, the uniform stability constants,
and the prescribed data norms, but independent of $h$,
such that
\begin{equation}\label{eq:posteriori-efficiency-2D}
C_{\rm eff}\Theta
\le
\mathtt E.
\end{equation}
\end{theorem}

\begin{proof}
We apply Lemma~\ref{lem:gradient-indicator-efficiency}
to the gradient identities
\eqref{eq:gradient-reconstruction}, with the field
and potential choices in
\eqref{eq:residual-localization-choices}.
For displacement, pressure, and electrostatic potential,
the gradient and potential exponents coincide.
Summing the local estimates
(cf. \eqref{eq:local-gradient-efficiency})
with the corresponding exponents, and using finite
patch overlap and $h_T\le\operatorname{diam}(\O)$,
bounds these gradient contributions by $C\mathtt E$.

For the ionic potential indicators, the exponent pair
$(\rho,2)$ introduces the factor $h_T^{2/\rho}$.
Since $\rho\ge4$, finite patch overlap yields
\[
\begin{aligned}
\left(
\sum_{T\in\Th}
h_T^2
\|\bt_i-\bt_{i,h}\|_{0;\omega_T}^{\rho}
\right)^{1/\rho}
\le
\operatorname{diam}(\O)^{2/\rho}
\left(
\sum_{T\in\Th}
\|\bt_i-\bt_{i,h}\|_{0;\omega_T}^{2}
\right)^{1/2}
\le C\|\bt_i-\bt_{i,h}\|_{0;\O}.
\end{aligned}
\]
The concentration errors are aggregated in $\rL^\rho$,
while the ionic rotational and tangential-jump
contributions are aggregated in $\rL^2$.
Thus, all gradient indicator contributions are bounded
by $C\mathtt E$.

The volume residuals
(cf. \eqref{def:posteriori-volume-residuals})
are estimated by subtracting the continuous equations
from their expressions at the discrete solution.
H\"older's inequality, with
$1/r+1/\rho=1/2$, bounds the electrostatic and advective
constitutive products. The fixed charge cancels from $r_{{\rm ch},h}$,
and contributes only $\Lambda_\chi\|\bvphi-\bvphi_h\|_{0,r;\O}$
to $r_{{\rm mom},h}$.
For the ionic balance products, we use
$1/\rho+1/2\le1/\varrho$ and the corresponding
Lebesgue embedding.
The remaining linear terms are bounded by the error
norms in \eqref{def:block-errors}.
Using the a priori estimates
(cf. \eqref{bound:a-priori} and \eqref{bound:a-priori-discr})
to bound the coefficient-field norms in terms of
the prescribed data, we obtain
\[
\begin{aligned}
&\|r_{{\rm mom},h}\|_{0;\O}
+\|r_{{\rm mass},h}\|_{0;\O}
+\|r_{{\rm ch},h}\|_{0,r;\O}
+
\sum_{i=1}^{2}
\left(
\|r_{i,{\rm con},h}\|_{0;\O}
+\|r_{i,{\rm bal},h}\|_{0,\varrho;\O}
\right)
\le C\mathtt E.
\end{aligned}
\]
Combining these bounds according to
\eqref{def:posteriori-local-indicators} and
\eqref{def:posteriori-global-estimator}
gives $\Theta\le C\mathtt E$, with $C$ independent of $h$.
Taking $C_{\rm eff}:=C^{-1}$ proves
\eqref{eq:posteriori-efficiency-2D}.
\end{proof}

\begin{remark}[Nonpolynomial boundary data]
\label{rem:posteriori-boundary-oscillation}

If the boundary data are not edgewise polynomials, their
approximation defect appears in the efficiency bound.
For data satisfying \eqref{eq:posteriori-boundary-regularity},
we measure this defect by boundary-data oscillation.
Let $\Pi_e^m$ denote the $L^2(e)$-orthogonal
projection onto $\rP_m(e)$, where
$m\ge\ell+1$ is fixed.
For a boundary datum $d_{\tt D}$ and exponents
$\varpi,\vartheta$, define
\begin{equation}\label{def:boundary-oscillation}
\begin{aligned}
\operatorname{osc}_{\Gamma}(d_{\tt D};\varpi,\vartheta)
:={}&
\left(
\sum_{e\in\cE_h^\Gamma}
h_e
\|d_{\tt D}-\Pi_e^m d_{\tt D}\|_{0,\varpi;e}^{\varpi}
\right)^{1/\varpi}
+
\left(
\sum_{e\in\cE_h^\Gamma}
h_e
\|\partial_{\boldsymbol{\tau}_e}
(d_{\tt D}-\Pi_e^m d_{\tt D})\|_{0,\vartheta;e}^{\vartheta}
\right)^{1/\vartheta},
\end{aligned}
\end{equation}
and set
\begin{equation}\label{def:total-boundary-oscillation}
\operatorname{osc}_{\Gamma}
:=
\operatorname{osc}_{\Gamma}(\chi_{\tt D};r,r)
+
\sum_{i=1}^2
\operatorname{osc}_{\Gamma}(\xi_{i,\tt D};\rho,2).
\end{equation}

Repeating the boundary-bubble estimates of
Lemma~\ref{lem:gradient-indicator-efficiency} after splitting each datum
into its polynomial projection and approximation defect
(cf. \cite[Section~6.2]{cmo-IMANUM-2016} and
\cite[Section~5.2]{ccov2022}) yields
$C_{\rm eff}\Theta\le\mathtt E+\operatorname{osc}_{\Gamma}$
with a possibly different $C_{\rm eff}>0$ independent of $h$.
The oscillation vanishes for representable data and is of higher order
for sufficiently smooth traces; it is used only for efficiency, and
reliability \eqref{eq:posteriori-reliability-2D} still holds for the
original data.
\end{remark}
%%%%%%%%%%%%%%%%%%%%%%%%%%%%%%%%%%%%%%%%%%
\subsection{Nonhomogeneous boundary conditions}
\label{subsec:posteriori-nonhomogeneous-boundary}

We now describe the changes when the
displacement and Darcy normal flux are also prescribed by
nonzero data, replacing the corresponding conditions in
\eqref{Biot-NPN-mixed-bc} by
\begin{equation}\label{eq:posteriori-nonhomogeneous-bc}
\u=\u_{\tt D},
\qquad
\z\cdot\bnu=z_{\tt N}
\qon\Gamma.
\end{equation}
Assume that $\u_{\tt D}\in\mathbf H^{1/2}(\Gamma)$,
$\u_{\tt D}|_e\in\mathbf H^1(e)$ on each boundary edge,
and $z_{\tt N}\in L^r(\Gamma)$.
Integration by parts makes the
displacement datum a load in the constitutive equation:
\begin{equation}\label{eq:mechanical-nonhomogeneous-weak}
\bfa_s(\bro,\bta)+\bfb_s(\bta,\u)
=\langle\bta\bnu,\u_{\tt D}\rangle_\Gamma
\qquad\forall\,\bta\in\bX.
\end{equation}

The Darcy datum prescribes the normal trace of the trial flux; the
test functions still have zero normal trace.
The mean constraints must be compatible with these traces:
the divergence theorem and \eqref{prop:divu:trrho} give
\begin{subequations}
\begin{align}
c_\lambda\int_\O\tr(\bro)
&=\int_\Gamma\u_{\tt D}\cdot\bnu,
\label{eq:nonhomogeneous-trace-mean}
\\
c_0\int_\O p
&=\int_\O g-\alpha\int_\Gamma\u_{\tt D}\cdot\bnu
-\int_\Gamma z_{\tt N}.
\label{eq:nonhomogeneous-mass-compatibility}
\end{align}
\end{subequations}

Identity \eqref{eq:nonhomogeneous-mass-compatibility} replaces
\eqref{eq:pressure-mean}. If $c_0>0$, it fixes the natural mean
$\overline p$ of $p$ in \eqref{def:p-mean}, with the right-hand side divided by
$c_0|\O|$, and imposes no condition on the data. If $c_0=0$, it is the
compatibility condition on $g$, $\u_{\tt D}$, and $z_{\tt N}$, and
\eqref{uniq-p-cond} fixes the mean of $p$. Both pseudostresses are assigned
the mean trace in \eqref{eq:nonhomogeneous-trace-mean}. The trial
spaces are therefore affine; subtracting fixed liftings leaves
increments in the original homogeneous spaces.

By \eqref{eq:mechanical-nonhomogeneous-weak}, the first residual
in \eqref{def:Biot-residual-functionals} becomes
\begin{equation}\label{def:mechanical-nonhomogeneous-residual}
\mathcal R_1^{\rm B}(\bta)
:=\langle\bta\bnu,\u_{\tt D}\rangle_\Gamma
-\bfa_s(\bro_h,\bta)-\bfb_s(\bta,\u_h).
\end{equation}
Its extension to $\bbH(\bdiv;\O)$ satisfies
$\mathcal R_1^{\rm B}(\I)=\int_\Gamma\u_{\tt D}\cdot\bnu
-c_\lambda\int_\O\tr(\bro_h)=0$, so the correction
\eqref{eq:mechanical-extended-orthogonality} still applies, and
Lemma~\ref{lem:gradient-residual-localization} gives the mechanical contribution
\begin{equation}\label{def:mechanical-nonhomogeneous-indicator}
\widetilde\theta_{{\rm B},2,T}^{\,2}
:=\|r_{{\rm mom},h}\|_{0;T}^2
+\|r_{{\rm mass},h}\|_{0;T}^2
+A_T(M_h,\u_h,\u_{\tt D};2)^2
+B_T(M_h,\u_{\tt D};2)^2.
\end{equation}
All other contributions in \eqref{def:posteriori-local-indicators}
are unchanged; in particular, the Darcy contributions remain $A_T^0$
and $B_T^0$, since their test functions have homogeneous normal trace.
Let $\Theta_{\tt D}$ denote the aggregation
\eqref{def:posteriori-global-estimator} with
$\theta_{{\rm B},2,T}$ replaced by
$\widetilde\theta_{{\rm B},2,T}$.

For general normal data, impose $\z_h\cdot\bnu=z_{{\tt N},h}$,
where $z_{{\tt N},h}|_e\in\rP_\ell(e)$ is the moment projection:
$\int_e(z_{\tt N}-z_{{\tt N},h})q_h=0$ for all $q_h\in\rP_\ell(e)$
and $e\in\cE_h^\Gamma$.
Define the normal-data approximation contribution by
\begin{equation}\label{def:normal-boundary-data-indicator}
\operatorname{osc}_{\Gamma,\tt N}
:=\left(
\sum_{e\in\cE_h^\Gamma}
h_e\|z_{\tt N}-z_{{\tt N},h}\|_{0,r;e}^{r}
\right)^{1/r},
\qquad
\Theta_{\rm nh}:=\Theta_{\tt D}+\operatorname{osc}_{\Gamma,\tt N}.
\end{equation}
The corresponding local contribution is obtained by restricting
the edge sum to $\cE_h^\Gamma(T)$.
If the normal datum is represented exactly, then $\operatorname{osc}_{\Gamma,\tt N}=0$.

The following result is an extension, not a complete theorem. We do not
prove well-posedness and uniform stability of the nonhomogeneous
continuous and discrete problems. We assume them, together with the data
bounds and the absorption condition of
Theorem~\ref{thm:posteriori-reliability-2D}, and we only indicate how
the proofs of Theorems~\ref{thm:posteriori-reliability-2D}
and~\ref{thm:posteriori-efficiency-2D} are modified.
\begin{lemma}[Extension to nonhomogeneous data, under assumptions]%\label{prop:posteriori-nonhomogeneous-boundary}
Assume that the continuous and discrete problems with
\eqref{eq:posteriori-nonhomogeneous-bc} are well posed.
Assume also the domain and mesh hypotheses of
Section~\ref{subsec:posteriori-preliminaries}, and the uniform
block stability, data-dependent field bounds, and residual
absorption hypothesis used in
Theorem~\ref{thm:posteriori-reliability-2D},
with the nonhomogeneous boundary data included in those bounds.
Then the estimator with normal-data oscillation
(cf. \eqref{def:normal-boundary-data-indicator}) satisfies
\begin{equation}\label{eq:nonhomogeneous-boundary-reliability}
\mathtt E\le C_{\rm rel}\Theta_{\rm nh},
\end{equation}
where the error is measured in the norms specified in
\eqref{def:block-errors}, and $C_{\rm rel}$ is independent of $h$.
Under the corresponding efficiency hypotheses of
Theorem~\ref{thm:posteriori-efficiency-2D},
\begin{equation}\label{eq:nonhomogeneous-boundary-efficiency}
C_{\rm eff}\Theta_{\rm nh}
\le\mathtt E
+\operatorname{osc}_{\Gamma,\tt D}
+\operatorname{osc}_{\Gamma,\tt N},
\end{equation}
where $\operatorname{osc}_{\Gamma,\tt D}
:=\operatorname{osc}_{\Gamma}
+\operatorname{osc}_{\Gamma}(\u_{\tt D};2,2)$.
The oscillation terms are defined in
\eqref{def:boundary-oscillation}--\eqref{def:total-boundary-oscillation},
with the displacement contribution understood componentwise.
\end{lemma}

\begin{proof}

For exactly imposed normal data, the fluxes have the same normal
trace and the pseudostresses the same mean trace. Their errors
belong to the homogeneous spaces, so the argument of
Theorem~\ref{thm:posteriori-reliability-2D} applies with
\eqref{def:mechanical-nonhomogeneous-residual} and
\eqref{def:mechanical-nonhomogeneous-indicator}.

For projected normal data, set $d_{\tt N}:=z_{\tt N}-z_{{\tt N},h}$,
whose mean vanishes on each boundary edge. Subtracting edge averages of a
test function in $\rW^{1/r,s}(\Gamma)$ and applying the fractional
Poincar\'e and H\"older inequalities
(cf. \cite[Lemma~3.8]{Goffi2021}) bounds its negative trace norm, and
surjectivity of the normal trace
\cite[Appendix~B, Theorem~5]{ld-CAMWA-2021} gives a lifting
$\boldsymbol\ell_h\in\bH^r(\div_r;\O)$ with
$\boldsymbol\ell_h\cdot\bnu=d_{\tt N}$ and
$\|\boldsymbol\ell_h\|_{r,\div;\O}\le C\operatorname{osc}_{\Gamma,\tt N}$
(using also $\rL^r\hookrightarrow\rL^2$ for $r>2$).
Block stability applied to the corrected error
$\z-\z_h-\boldsymbol\ell_h$, which has zero normal trace, together with
continuity of the forms and the assumed field bounds, controls the lifting
terms by $C\operatorname{osc}_{\Gamma,\tt N}$. Localizing the original
residuals by Galerkin orthogonality and absorbing the coupling
proves \eqref{eq:nonhomogeneous-boundary-reliability}.

For efficiency, use the row-wise boundary estimates of
Lemma~\ref{lem:gradient-indicator-efficiency} with
$\u_{\tt D}$ and the projected-data argument of
Remark~\ref{rem:posteriori-boundary-oscillation}.
This bounds $\Theta_{\tt D}$ by the error and the stated
Dirichlet oscillation. Adding $\operatorname{osc}_{\Gamma,\tt N}$ proves
\eqref{eq:nonhomogeneous-boundary-efficiency}.
\end{proof}

\section{Numerical results}\label{sec:numerical-results}

\begingroup
\renewcommand{\topfraction}{0.9}
\renewcommand{\bottomfraction}{0.6}
\renewcommand{\textfraction}{0.07}
\renewcommand{\floatpagefraction}{0.85}
\setcounter{topnumber}{3}
\setcounter{totalnumber}{4}

This section illustrates the accuracy and performance of the
fully mixed finite element scheme
\eqref{Biot-PNP-fully-mixed-formulation-disc} and examines
the reliability and efficiency of the residual-based
{\it a posteriori} error estimator $\Theta$
(cf.~\eqref{def:posteriori-global-estimator}).
The stability, radius, and smallness conditions of
Sections~\ref{section4}--\ref{sec:posteriori-2D}
are sufficient conditions for the theoretical results, and
we do not verify them for the parameters used below.
The experiments therefore also probe the behavior of the scheme
beyond these restrictions.
We consider three examples using the finite element subspaces
defined in \eqref{def-FE}.
The choices $\ell=0$ and $\ell=1$ are denoted by
$\mathbb{RT}_0$--$\mathbf P_0$ and
$\mathbb{RT}_1$--$\mathbf P_1$, respectively.
This notation includes the row-wise Raviart--Thomas
approximation of the pseudostress and the vector-valued
polynomial spaces for the displacement and concentration
gradients.
Throughout the experiments, we take
$r=\rho=4$ and $s=\varrho=4/3$.

The implementation uses {\tt FEniCS} \cite{lmw-FEniCS-2012},
its built-in Newton updates, and an automatically differentiated
Jacobian.
Iteration stops when consecutive coefficient vectors satisfy
\[
\frac{|\mathbf{coeff}^{m+1}-\mathbf{coeff}^{m}|_{\DOF}}
{|\mathbf{coeff}^{m+1}|_{\DOF}}
\le {\tt tol},
\]
where $|\cdot|_{\DOF}$ is the Euclidean norm and $\DOF$ counts
all degrees of freedom in \eqref{def-FE}.
We set ${\tt tol}=10^{-6}$ and allow at most $30$ Newton
updates per nonlinear solve.
The zero-mean conditions on $p_h$ and $\tr(\bro_h)$ are
enforced by two real Lagrange multipliers, which are
included in this count. In all examples the natural pressure mean vanishes,
$\overline p=0$ in \eqref{def:p-mean}, so $p_h$ is also the full discrete
pressure (Lemma~\ref{lem:pressure-mean}).
The quadrature degrees are $16$, $40$, and $20$ in
Examples~1, 2, and~3, respectively.

We report the block errors defined in \eqref{def:block-errors}:
$\mathtt E_1$ corresponds to the Biot variables
$(\bro,\u,\z,p)$,
$\mathtt E_2$ to the Poisson variables $(\bvphi,\chi)$,
and $\mathtt E_3$ to the Nernst--Planck variables
$(\bt_i,\bsi_i,\xi_i)$ for both species.
The global error and effectivity index are
\[
\mathtt E:=\mathtt E_1+\mathtt E_2+\mathtt E_3,
\qquad
\operatorname{eff}(\Theta):=\frac{\mathtt E}{\Theta}.
\]
Values of the effectivity index below or above $1$ indicate
an overestimation or underestimation of the error,
respectively.
The reliability and efficiency estimates
(cf.~\eqref{eq:posteriori-reliability-2D} and
\eqref{eq:posteriori-efficiency-2D}) imply that this index
is bounded above and away from zero under their
corresponding hypotheses.

For two consecutive meshes with $\DOF$ and $\DOF'$
degrees of freedom and corresponding error or estimator
values $e$ and $e'$, we compute
\begin{equation}\label{def:numerical-uniform-rate}
\mathcal R(e)
:=-d\,
\frac{\log(e/e')}{\log(\DOF/\DOF')}.
\end{equation}
For uniformly refined two-dimensional meshes,
$\DOF^{-1/d}\simeq h$, so
\eqref{def:numerical-uniform-rate} measures the order
with respect to the meshsize.
On adaptive meshes, it measures decay with respect to
$\DOF^{-1/2}$.
Each $\mathcal R$ column in the tables refers to the
quantity immediately preceding it.

The adaptive computations follow the
solve--estimate--mark--refine procedure of \cite{Verfurth}.
After solving
\eqref{Biot-PNP-fully-mixed-formulation-disc}
on the current mesh, we compute
\[
\Theta_T
:=
\theta_{{\rm B},2,T}
+\theta_{{\rm B},r,T}
+\theta_{{\rm P},r,T}
+\sum_{i=1}^{2}
\bigl\{
\theta_{i,2,T}
+\theta_{i,\rho,T}
+\theta_{i,\varrho,T}
\bigr\},
\qquad T\in\Th,
\]
from \eqref{def:posteriori-local-indicators}.
Unless the prescribed mesh or DOF limit has been reached,
we mark each $T'\in\Th$ satisfying
\begin{equation}\label{eq:numerical-mean-marking}
\Theta_{T'}
\ge
C_{\tt adm}\,\frac{1}{\#\Th}
\sum_{T\in\Th}\Theta_T,
\qquad C_{\tt adm}=0.8,
\end{equation}
where $\#\Th$ denotes the number of triangles.
The marked elements are refined using the conforming
FEniCS refinement procedure.
This mean-value strategy is a heuristic: we do not claim D\"orfler-type
marking properties or convergence of the adaptive loop, and we use it only
to test the reliability and efficiency of $\Theta$ in practice.

\subsection*{Example 1: A smooth manufactured solution}

We assess convergence under uniform refinement on
$\O=(0,1)^2$, with
\[
\mu=\kappa=\nu=\varepsilon=\kappa_1=1,
\qquad
\lambda=2,
\qquad
\alpha=c_0=\delta=0.1,
\qquad
\kappa_2=0.7.
\]
The sources and boundary data are chosen so that the
exact primary fields are
\begin{align*}
\u(x,y)
&=10^{-3}x(1-x)y(1-y)(1,x)^{\tt t},
\qquad
p(x,y)=10^{-3}\cos(\pi x)\cos(\pi y),
\\[1ex]
\chi(x,y)
&=10^{-3}\sin(\pi x)\sin(\pi y),
\qquad
\xi_i(x,y)=10^{-3}\bigl\{10+i\,x(1-x)y(1-y)\bigr\},
\quad i\in\{1,2\}.
\end{align*}
The remaining exact fields follow from the constitutive
relations in \eqref{eq:strong-complete}.
A prescribed charge source $f_\chi$ allows the potential
and concentrations to be chosen independently.
It enters both the Poisson balance and the Maxwell force
through the total charge $\xi_1-\xi_2+f_\chi$
(cf. \eqref{Biot-NPN-mixed-2} and \eqref{Biot-NPN-mixed-5-2}),
and the estimator through \eqref{def:posteriori-volume-residuals}.
The exact fields satisfy the homogeneous displacement
and Darcy normal-flux conditions, and $\int_\O p=0$
(cf.~\eqref{Biot-NPN-mixed-bc},
\eqref{prop:tr:mean:zero}, and \eqref{eq:pressure-mean}).

The meshes are obtained by dividing the square into
$n\times n$ subsquares and splitting each into two
triangles, giving $2n^2$ elements and $h=\sqrt2/n$.
We use $n=4,8,16,32,64$ for $\ell=0$ and
$n=4,8,16,32$ for $\ell=1$.
Newton's method starts from the zero coefficient vector
on each mesh.

Table~\ref{tab:numerical-smooth} shows first- and
second-order convergence for
$\mathbb{RT}_0$--$\mathbf P_0$ and
$\mathbb{RT}_1$--$\mathbf P_1$, respectively,
for all three block errors and the global error.
The Biot and Poisson contributions dominate the global
error.
The estimator exhibits the same decay, with effectivity
indices approaching $0.94$ and $0.60$, respectively.

\begin{table}[!t]
\centering
\small
\setlength{\tabcolsep}{2pt}
\renewcommand{\arraystretch}{0.95}
\setlength{\aboverulesep}{1.5pt}
\setlength{\belowrulesep}{1.5pt}
\begin{tabular*}{\textwidth}{@{\extracolsep{\fill}}r*{5}{rr}r@{}}
\toprule
& \multicolumn{2}{c}{Biot}
& \multicolumn{2}{c}{Poisson}
& \multicolumn{2}{c}{Nernst--Planck}
& \multicolumn{2}{c}{Global}
& \multicolumn{2}{c}{Estimator}
& \\
\cmidrule(lr){2-3}
\cmidrule(lr){4-5}
\cmidrule(lr){6-7}
\cmidrule(lr){8-9}
\cmidrule(lr){10-11}
$\DOF$
& $\mathtt E_1$ & $\mathcal R$
& $\mathtt E_2$ & $\mathcal R$
& $\mathtt E_3$ & $\mathcal R$
& $\mathtt E$ & $\mathcal R$
& $\Theta$ & $\mathcal R$
& $\operatorname{eff}(\Theta)$
\\
\midrule
\multicolumn{12}{@{}l}{$\mathbb{RT}_0$--$\mathbf P_0$}\\[1pt]
658
& 4.141e-03 & --
& 4.180e-03 & --
& 4.412e-04 & --
& 8.762e-03 & --
& 8.911e-03 & --
& 0.9833
\\
2530
& 2.105e-03 & 1.00
& 2.130e-03 & 1.00
& 2.260e-04 & 0.99
& 4.461e-03 & 1.00
& 4.647e-03 & 0.97
& 0.9600
\\
9922
& 1.057e-03 & 1.01
& 1.070e-03 & 1.01
& 1.137e-04 & 1.01
& 2.241e-03 & 1.01
& 2.358e-03 & 0.99
& 0.9501
\\
39298
& 5.290e-04 & 1.01
& 5.357e-04 & 1.01
& 5.696e-05 & 1.00
& 1.122e-03 & 1.01
& 1.186e-03 & 1.00
& 0.9459
\\
156418
& 2.645e-04 & 1.00
& 2.680e-04 & 1.00
& 2.849e-05 & 1.00
& 5.610e-04 & 1.00
& 5.942e-04 & 1.00
& 0.9441
\\
\addlinespace[3pt]
\multicolumn{12}{@{}l}{$\mathbb{RT}_1$--$\mathbf P_1$}\\[1pt]
2018
& 5.940e-04 & --
& 6.830e-04 & --
& 6.222e-05 & --
& 1.339e-03 & --
& 2.198e-03 & --
& 0.6091
\\
7874
& 1.509e-04 & 2.01
& 1.743e-04 & 2.01
& 1.583e-05 & 2.01
& 3.410e-04 & 2.01
& 5.674e-04 & 1.99
& 0.6011
\\
31106
& 3.789e-05 & 2.01
& 4.381e-05 & 2.01
& 3.982e-06 & 2.01
& 8.568e-05 & 2.01
& 1.434e-04 & 2.00
& 0.5975
\\
123650
& 9.487e-06 & 2.01
& 1.097e-05 & 2.01
& 9.982e-07 & 2.01
& 2.145e-05 & 2.01
& 3.600e-05 & 2.00
& 0.5959
\\
\bottomrule
\end{tabular*}
\caption{Example~1: block errors, global error, estimator,
experimental convergence rates, and effectivity index
under uniform refinement.}
\label{tab:numerical-smooth}
\end{table}

We also evaluate the displacement gradient and Cauchy stress
reconstructed through
\eqref{eq:gradient-reconstruction} and
\eqref{eq:Cauchy-reconstruction}.
For the individual mechanical and reconstructed errors,
we use
\[
\begin{aligned}
\mathrm e(\bro)&:=\|\bro-\bro_h\|_{\bX},
&
\mathrm e(\u)&:=\|\u-\u_h\|_{\bQ},
&
\mathrm e(\nabla\u)&:=\|\nabla\u-(\nabla\u)_h\|_{0;\O},
&
\mathrm e(\bsi_{\rm p})&:=\|\bsi_{\rm p}-\bsi_{{\rm p},h}\|_{0;\O}.
\end{aligned}
\]
Under uniform refinement, the errors of both reconstructed
fields approach orders one and two for $\ell=0$ and
$\ell=1$, respectively.
Their rates on the finest meshes are $0.99$ for $\ell=0$
and $1.98$ for $\ell=1$.

Finally, we examine robustness with respect to $\lambda$
on the fixed mesh with $n=16$, containing $512$ triangles.
For this test, we prescribe
\[
\psi(x,y):=x^2(1-x)^2y^2(1-y)^2,
\qquad
\u_\lambda
:=10^{-3}\left[
\bigl(\partial_y\psi,-\partial_x\psi\bigr)^{\tt t}
+\frac{\mu}{\lambda+\mu}\nabla\psi
\right].
\]
This displacement satisfies $\u_\lambda=\mathbf0$ on $\Gamma$
and has nonzero divergence for every finite $\lambda$, with
\[
\div(\u_\lambda)
=10^{-3}\frac{\mu}{\lambda+\mu}\Delta\psi,
\qquad
(\lambda+\mu)\div(\u_\lambda)
=10^{-3}\mu\Delta\psi.
\]
Thus, the volumetric contribution is retained while the exact
pseudostress and prescribed sources remain uniformly bounded
as $\lambda$ increases.
The exact pressure, electric potential, and concentrations
remain unchanged, and the remaining fields and sources
are recomputed for each $\lambda$.
We retain $\mu=1$ and consider
$\lambda\in\{2,10^2,10^4,10^6,10^8,10^{10}\}$.
The systems have $9{,}922$ and $31{,}106$ DOF
for $\ell=0$ and $\ell=1$, respectively.

Table~\ref{tab:numerical-lambda-robustness} shows that the
mechanical errors stabilize as $\lambda$ increases, with
no observed deterioration in the nearly incompressible
regime for this manufactured family.
Newton's method converges in three updates for every case.
Using the reconstruction
$(\lambda+\mu)c_\lambda\tr(\bro_h)$, the relative $\rL^2$
error in $(\lambda+\mu)\div(\u_\lambda)$ is approximately
$12.0\%$ for $\ell=0$ and $0.55\%$ for $\ell=1$
at $\lambda=10^{10}$.
These results support the robustness of the approximation
with respect to $\lambda$.

\begin{table}[!htbp]
\centering
\small
\setlength{\tabcolsep}{2pt}
\renewcommand{\arraystretch}{0.95}
\begin{tabular*}{\textwidth}{@{\extracolsep{\fill}}rrrrrrrrr@{}}
\toprule
& \multicolumn{4}{c}{$\mathbb{RT}_0$--$\mathbf P_0$}
& \multicolumn{4}{c}{$\mathbb{RT}_1$--$\mathbf P_1$}
\\
\cmidrule(lr){2-5}\cmidrule(lr){6-9}
$\lambda$
& $\mathrm e(\bro)$
& $\mathrm e(\u)$
& $\mathrm e(\nabla\u)$
& $\mathrm e(\bsi_{\rm p})$
& $\mathrm e(\bro)$
& $\mathrm e(\u)$
& $\mathrm e(\nabla\u)$
& $\mathrm e(\bsi_{\rm p})$
\\
\midrule
$2$
& 9.146e-05 & 9.082e-07 & 1.308e-05 & 2.391e-05
& 4.780e-06 & 6.210e-08 & 6.305e-07 & 1.151e-06
\\
$10^2$
& 7.943e-05 & 8.692e-07 & 1.175e-05 & 2.237e-05
& 4.156e-06 & 5.894e-08 & 5.976e-07 & 1.110e-06
\\
$10^4$
& 7.911e-05 & 8.697e-07 & 1.172e-05 & 2.235e-05
& 4.139e-06 & 5.894e-08 & 5.970e-07 & 1.110e-06
\\
$10^6$
& 7.910e-05 & 8.697e-07 & 1.172e-05 & 2.235e-05
& 4.139e-06 & 5.894e-08 & 5.970e-07 & 1.110e-06
\\
$10^8$
& 7.910e-05 & 8.697e-07 & 1.172e-05 & 2.235e-05
& 4.139e-06 & 5.894e-08 & 5.970e-07 & 1.110e-06
\\
$10^{10}$
& 7.910e-05 & 8.697e-07 & 1.172e-05 & 2.235e-05
& 4.139e-06 & 5.894e-08 & 5.970e-07 & 1.110e-06
\\
\bottomrule
\end{tabular*}
\caption{Example~1: mechanical errors as $\lambda$ increases
on a fixed mesh of $512$ triangles, with $9{,}922$ and
$31{,}106$ DOF for $\ell=0$ and $\ell=1$, respectively.
The displacement gradient and Cauchy stress are reconstructed
algebraically.}
\label{tab:numerical-lambda-robustness}
\end{table}

\subsection*{Example 2: Adaptive resolution of an internal layer}

We assess estimator-driven adaptive refinement for a solution
with an internal layer, posed on the interior of the polygon
with successive vertices
$\mathbf0,\boldsymbol v_0,\ldots,\boldsymbol v_{12}$, where
\[
\boldsymbol v_j
=
r_j\bigl(\cos(j\pi/8),\sin(j\pi/8)\bigr)^{\tt t},
\qquad
r_j=
\begin{cases}
1.45,&j\ \text{even},\\
0.95,&j\ \text{odd}.
\end{cases}
\]
The domain has a re-entrant corner of angle $3\pi/2$
at the origin. The manufactured solution below is smooth on
$\overline\O$, so this corner introduces no singularity in the
exact fields, and refinement is driven by the internal layer.
We use the material parameters of the uniform-refinement
test in Example~1, with $\lambda=2$.
Writing
\[
H(x,y):=\tanh\!\left(\frac{x+y}{w}\right),
\qquad
K(x,y):=1-H(x,y)^2,
\qquad w=0.1,
\]
we manufacture the sources and boundary data so that
\[
\begin{gathered}
\u=aK(1,1)^{\tt t},
\qquad
p=aH,
\qquad
\chi=aK,
\\[1ex]
\xi_1=b+AH,
\qquad
\xi_2=b-AH,
\qquad
a=10^{-3},\quad b=2\times10^{-2},\quad A=5\times10^{-3}.
\end{gathered}
\]
These profiles exhibit localized variations near $y=-x$,
across which the charge difference
$\xi_1-\xi_2=2AH$ changes sign.
The remaining exact fields follow from
\eqref{eq:strong-complete}, and the functionals and
residuals are adjusted consistently with the
manufactured sources.
Boundary data are the exact traces, with the Darcy
normal flux imposed through its edgewise polynomial
projection.
We use the extension in
Section~\ref{subsec:posteriori-nonhomogeneous-boundary}
and denote its estimator $\Theta_{\rm nh}$ by $\Theta$,
including the normal-data oscillation term.
Since that extension is stated under assumptions and not proved in full,
the bounds of Theorems~\ref{thm:posteriori-reliability-2D}
and~\ref{thm:posteriori-efficiency-2D} do not cover this example
rigorously; it is a numerical test of the extended estimator.
Reflection across $y=-x$ makes $p$ and $\div(\u)$ odd,
so their means vanish, and the natural pressure mean in
\eqref{eq:nonhomogeneous-mass-compatibility} is $\overline p=0$ for these
boundary data.

Both refinement strategies start from the same mesh
with $12$ triangles.
The adaptive sequences use
\eqref{eq:numerical-mean-marking} with $C_{\tt adm}=0.8$.
Tables~\ref{tab:numerical-front-rt0}
and~\ref{tab:numerical-front-rt1} compare uniform
and adaptive refinement for both polynomial degrees.
The initially slow decay reflects the insufficient
resolution of the transition width $w$.
Adaptive refinement reaches comparable convergence
rates with fewer DOF.
For $\ell=0$, the approximation with $86{,}524$ adaptive
DOF has a global error approximately $34\%$ smaller
than that on the last uniform mesh with $234{,}818$ DOF.
For $\ell=1$, the approximation with $97{,}454$ adaptive
DOF has a global error approximately $38\%$ smaller
than that on the last uniform mesh with $185{,}666$ DOF.
These comparisons require about $37\%$ and $52\%$
of the corresponding uniform DOF counts.

\begin{table}[!t]
\centering
\small
\setlength{\tabcolsep}{2pt}
\renewcommand{\arraystretch}{0.95}
\setlength{\aboverulesep}{1.5pt}
\setlength{\belowrulesep}{1.5pt}
\begin{tabular*}{\textwidth}{@{\extracolsep{\fill}}r*{5}{rr}r@{}}
\toprule
& \multicolumn{2}{c}{Biot}
& \multicolumn{2}{c}{Poisson}
& \multicolumn{2}{c}{Nernst--Planck}
& \multicolumn{2}{c}{Global}
& \multicolumn{2}{c}{Estimator}
& \\
\cmidrule(lr){2-3}
\cmidrule(lr){4-5}
\cmidrule(lr){6-7}
\cmidrule(lr){8-9}
\cmidrule(lr){10-11}
$\DOF$
& $\mathtt E_1$ & $\mathcal R$
& $\mathtt E_2$ & $\mathcal R$
& $\mathtt E_3$ & $\mathcal R$
& $\mathtt E$ & $\mathcal R$
& $\Theta$ & $\mathcal R$
& $\operatorname{eff}(\Theta)$
\\
\midrule
\multicolumn{12}{@{}l}{Uniform refinement}\\[1pt]
272
& 7.200e-01 & --
& 1.849e-01 & --
& 4.478e-01 & --
& 1.353e+00 & --
& 1.382e+00 & --
& 0.9786
\\
998
& 6.760e-01 & 0.10
& 1.692e-01 & 0.14
& 3.530e-01 & 0.37
& 1.198e+00 & 0.19
& 1.235e+00 & 0.17
& 0.9704
\\
3818
& 5.387e-01 & 0.34
& 1.530e-01 & 0.15
& 2.643e-01 & 0.43
& 9.560e-01 & 0.34
& 9.830e-01 & 0.34
& 0.9725
\\
14930
& 3.854e-01 & 0.49
& 1.220e-01 & 0.33
& 1.520e-01 & 0.81
& 6.594e-01 & 0.54
& 6.764e-01 & 0.55
& 0.9749
\\
59042
& 2.216e-01 & 0.80
& 7.436e-02 & 0.72
& 8.616e-02 & 0.83
& 3.822e-01 & 0.79
& 3.894e-01 & 0.80
& 0.9813
\\
234818
& 1.112e-01 & 1.00
& 3.691e-02 & 1.02
& 4.528e-02 & 0.93
& 1.933e-01 & 0.99
& 1.966e-01 & 0.99
& 0.9834
\\
\addlinespace[3pt]
\multicolumn{12}{@{}l}{Adaptive refinement}\\[1pt]
272
& 7.200e-01 & --
& 1.849e-01 & --
& 4.478e-01 & --
& 1.353e+00 & --
& 1.382e+00 & --
& 0.9786
\\
528
& 6.775e-01 & 0.18
& 1.693e-01 & 0.27
& 3.533e-01 & 0.71
& 1.200e+00 & 0.36
& 1.239e+00 & 0.33
& 0.9683
\\
1540
& 5.432e-01 & 0.41
& 1.529e-01 & 0.19
& 2.753e-01 & 0.47
& 9.713e-01 & 0.40
& 9.969e-01 & 0.41
& 0.9743
\\
4028
& 3.947e-01 & 0.66
& 1.209e-01 & 0.49
& 1.795e-01 & 0.89
& 6.951e-01 & 0.70
& 7.078e-01 & 0.71
& 0.9821
\\
10094
& 2.426e-01 & 1.06
& 8.523e-02 & 0.76
& 1.076e-01 & 1.11
& 4.355e-01 & 1.02
& 4.418e-01 & 1.03
& 0.9858
\\
29034
& 1.353e-01 & 1.11
& 5.524e-02 & 0.82
& 6.028e-02 & 1.10
& 2.508e-01 & 1.04
& 2.546e-01 & 1.04
& 0.9851
\\
86524
& 6.752e-02 & 1.27
& 2.516e-02 & 1.44
& 3.553e-02 & 0.97
& 1.282e-01 & 1.23
& 1.306e-01 & 1.22
& 0.9818
\\
\bottomrule
\end{tabular*}
\caption{Example~2: comparison of uniform and adaptive refinement
using $\mathbb{RT}_0$--$\mathbf P_0$.}
\label{tab:numerical-front-rt0}
\end{table}

\begin{table}[!t]
\centering
\small
\setlength{\tabcolsep}{2pt}
\renewcommand{\arraystretch}{0.95}
\setlength{\aboverulesep}{1.5pt}
\setlength{\belowrulesep}{1.5pt}
\begin{tabular*}{\textwidth}{@{\extracolsep{\fill}}r*{5}{rr}r@{}}
\toprule
& \multicolumn{2}{c}{Biot}
& \multicolumn{2}{c}{Poisson}
& \multicolumn{2}{c}{Nernst--Planck}
& \multicolumn{2}{c}{Global}
& \multicolumn{2}{c}{Estimator}
& \\
\cmidrule(lr){2-3}
\cmidrule(lr){4-5}
\cmidrule(lr){6-7}
\cmidrule(lr){8-9}
\cmidrule(lr){10-11}
$\DOF$
& $\mathtt E_1$ & $\mathcal R$
& $\mathtt E_2$ & $\mathcal R$
& $\mathtt E_3$ & $\mathcal R$
& $\mathtt E$ & $\mathcal R$
& $\Theta$ & $\mathcal R$
& $\operatorname{eff}(\Theta)$
\\
\midrule
\multicolumn{12}{@{}l}{Uniform refinement}\\[1pt]
806
& 6.445e-01 & --
& 1.533e-01 & --
& 3.130e-01 & --
& 1.111e+00 & --
& 1.553e+00 & --
& 0.7151
\\
3050
& 4.588e-01 & 0.51
& 1.257e-01 & 0.30
& 2.494e-01 & 0.34
& 8.339e-01 & 0.43
& 1.106e+00 & 0.51
& 0.7541
\\
11858
& 2.831e-01 & 0.71
& 9.801e-02 & 0.37
& 1.160e-01 & 1.13
& 4.972e-01 & 0.76
& 5.843e-01 & 0.94
& 0.8509
\\
46754
& 1.299e-01 & 1.14
& 5.180e-02 & 0.93
& 4.841e-02 & 1.27
& 2.301e-01 & 1.12
& 2.458e-01 & 1.26
& 0.9362
\\
185666
& 3.551e-02 & 1.88
& 1.375e-02 & 1.92
& 1.436e-02 & 1.76
& 6.362e-02 & 1.86
& 6.759e-02 & 1.87
& 0.9412
\\
\addlinespace[3pt]
\multicolumn{12}{@{}l}{Adaptive refinement}\\[1pt]
806
& 6.445e-01 & --
& 1.533e-01 & --
& 3.130e-01 & --
& 1.111e+00 & --
& 1.553e+00 & --
& 0.7151
\\
1670
& 4.639e-01 & 0.90
& 1.250e-01 & 0.56
& 2.647e-01 & 0.46
& 8.536e-01 & 0.72
& 1.130e+00 & 0.87
& 0.7555
\\
5318
& 2.902e-01 & 0.81
& 9.618e-02 & 0.45
& 1.495e-01 & 0.99
& 5.358e-01 & 0.80
& 6.322e-01 & 1.00
& 0.8475
\\
13154
& 1.565e-01 & 1.36
& 7.570e-02 & 0.53
& 6.836e-02 & 1.73
& 3.006e-01 & 1.28
& 3.316e-01 & 1.43
& 0.9064
\\
35306
& 8.367e-02 & 1.27
& 4.720e-02 & 0.96
& 2.259e-02 & 2.24
& 1.535e-01 & 1.36
& 1.633e-01 & 1.44
& 0.9399
\\
97454
& 2.004e-02 & 2.81
& 1.205e-02 & 2.69
& 7.412e-03 & 2.19
& 3.950e-02 & 2.67
& 4.113e-02 & 2.72
& 0.9605
\\
300770
& 5.636e-03 & 2.25
& 3.899e-03 & 2.00
& 2.285e-03 & 2.09
& 1.182e-02 & 2.14
& 1.253e-02 & 2.11
& 0.9436
\\
\bottomrule
\end{tabular*}
\caption{Example~2: comparison of uniform and adaptive refinement
using $\mathbb{RT}_1$--$\mathbf P_1$.}
\label{tab:numerical-front-rt1}
\end{table}

The global adaptive rates on the final refinement steps
are $1.23$ and $2.14$ for $\ell=0$ and $\ell=1$,
respectively.
The estimator follows the error decay, with effectivity
indices between $0.9683$ and $0.9858$ for $\ell=0$
and between $0.7151$ and $0.9605$ for $\ell=1$.
Figure~\ref{fig:numerical-front-fields} shows the final
adaptive mesh and two computed fields for $\ell=1$.
Refinement concentrates along the transition near
$y=-x$.
\begin{figure}[!htbp]
\centering
\setlength{\tabcolsep}{3pt}
\begin{tabular}{@{}ccc@{}}
\raisebox{-\height}{%
\includegraphics[
    width=0.32\textwidth,
    trim=18bp 18bp 18bp 18bp,
    clip
]{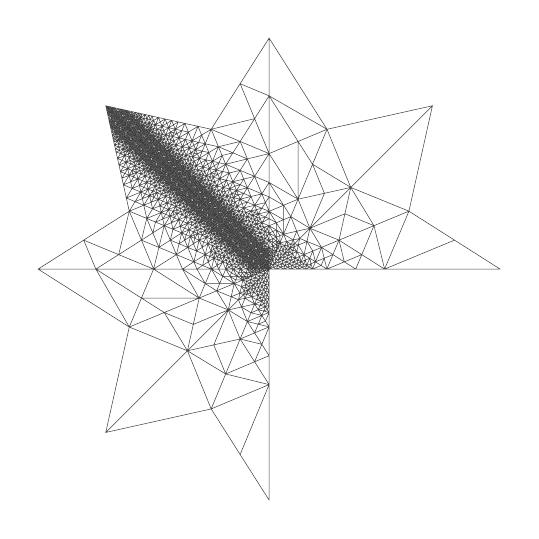}}
&
\raisebox{-\height}{%
\includegraphics[
    width=0.32\textwidth,
    trim=11.28bp 0bp 11.28bp 6.96bp,
    clip
]{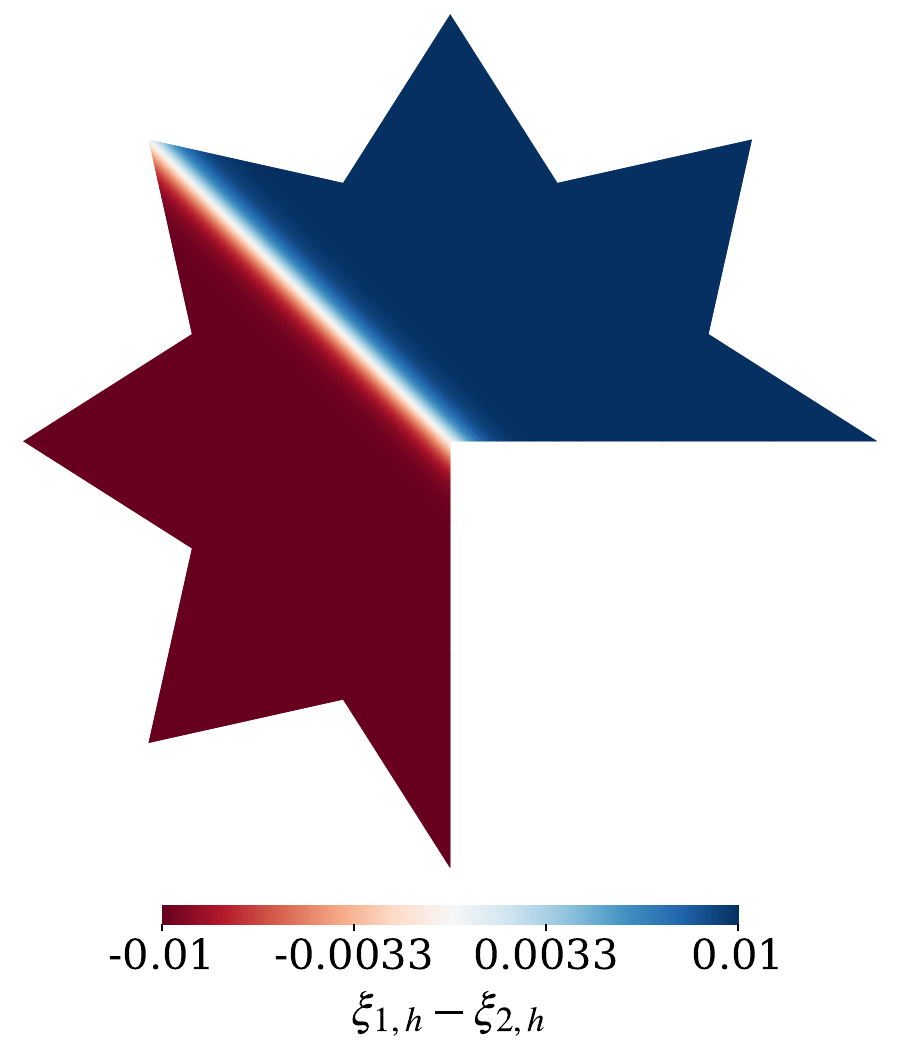}}
&
\raisebox{-\height}{%
\includegraphics[
    width=0.32\textwidth,
    trim=11.28bp 0bp 11.28bp 6.96bp,
    clip
]{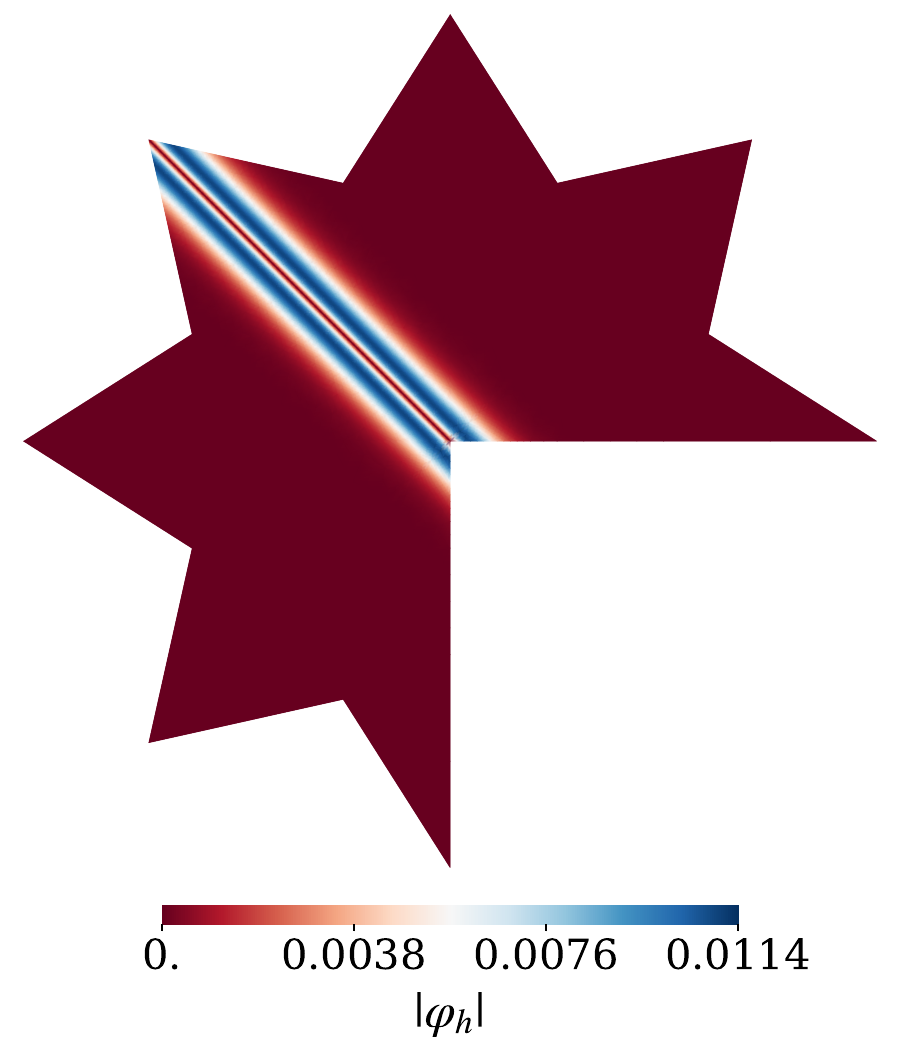}}
\end{tabular}
\caption{Example~2: final adaptive mesh, computed charge difference,
and electric-displacement magnitude (left to right),
using $\mathbb{RT}_1$--$\mathbf P_1$ with $300{,}770$ DOF.}
\label{fig:numerical-front-fields}
\end{figure}

\subsection*{Example 3: Ionic polarization in a charged porous medium}

Inspired by the electrical stimulation of charged hydrogels
studied in \cite{Farooqi2019}, we examine a confined porous
medium with a localized negative fixed charge.
The same poroelastic material occupies the entire square
$\O=(0,1)^2$, and all quantities are dimensionless.
We set
\[
\mu=\kappa=\nu=\kappa_1=1,
\qquad
\lambda=2,
\qquad
\alpha=c_0=\delta=0.1,
\qquad
\kappa_2=0.7,
\qquad
\varepsilon=2\times10^{-4}.
\]
The prescribed fixed-charge source is
\begin{equation}\label{eq:numerical-fixed-charge}
f_\chi(x,y)
:=
-\frac{C_{\rm f}}{2}
\left\{
1+\tanh\left(
\frac{R^2-(x-\tfrac12)^2-(y-\tfrac12)^2}{2R\eta}
\right)
\right\},
\qquad
C_{\rm f}=0.04,\quad R=0.2,\quad\eta=0.02.
\end{equation}
The profile in \eqref{eq:numerical-fixed-charge}
concentrates the charge in a central circular region
with a smooth transition of fixed width.
Both the Poisson load and the electrostatic force
use the total charge $\xi_1-\xi_2+f_\chi$
(cf. \eqref{Biot-NPN-mixed-2} and \eqref{Biot-NPN-mixed-5-2}).
We prescribe
\begin{equation}\label{eq:numerical-charged-boundary}
\u=\mathbf0,
\qquad
\z\cdot\bnu=0,
\qquad
\chi_{\tt D}=V(x-\tfrac12),
\qquad
\xi_{1,\tt D}=\xi_{2,\tt D}=0.02\qon\Gamma.
\end{equation}
We take $\ff=\mathbf0$, $g=0$, and $f_1=f_2=0.02$.
These ionic sources balance the reaction term at
the uncharged rest state.

\begin{figure}[!tb]
\centering
\includegraphics[width=0.3\textwidth]{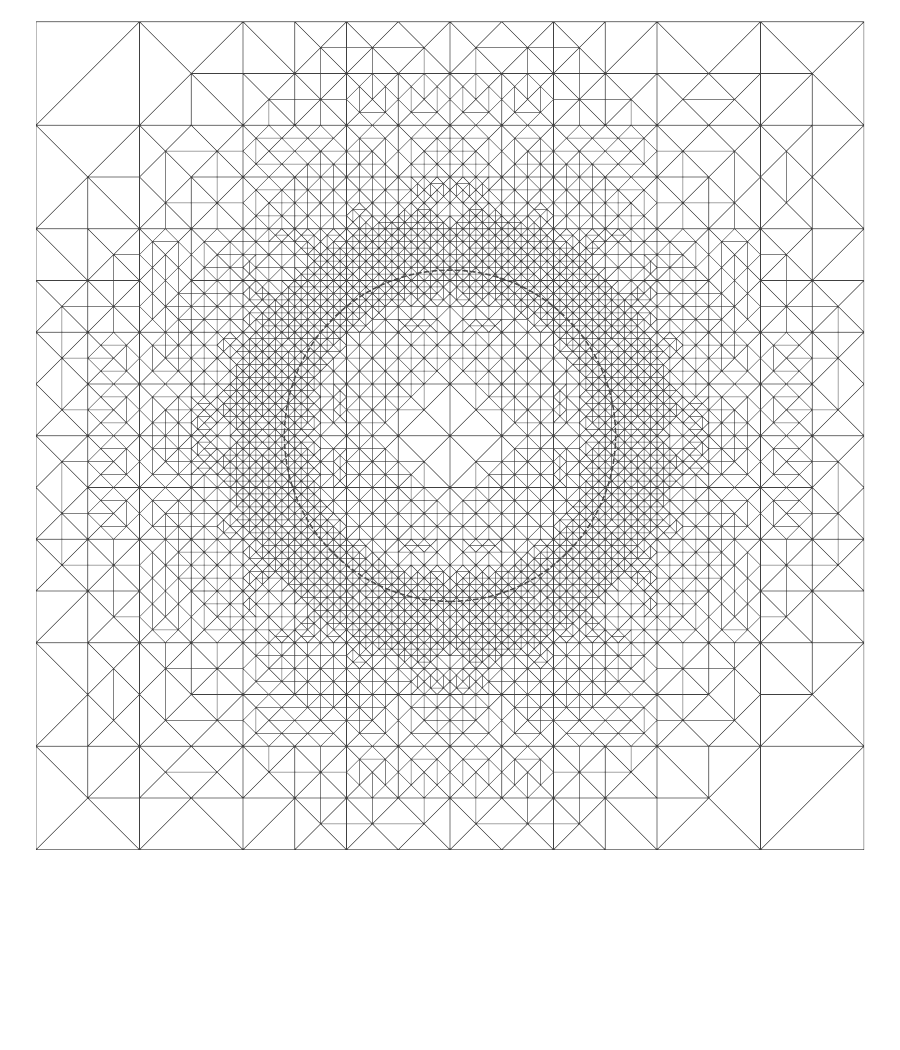}
\hspace{0.02\textwidth}
\includegraphics[width=0.3\textwidth]{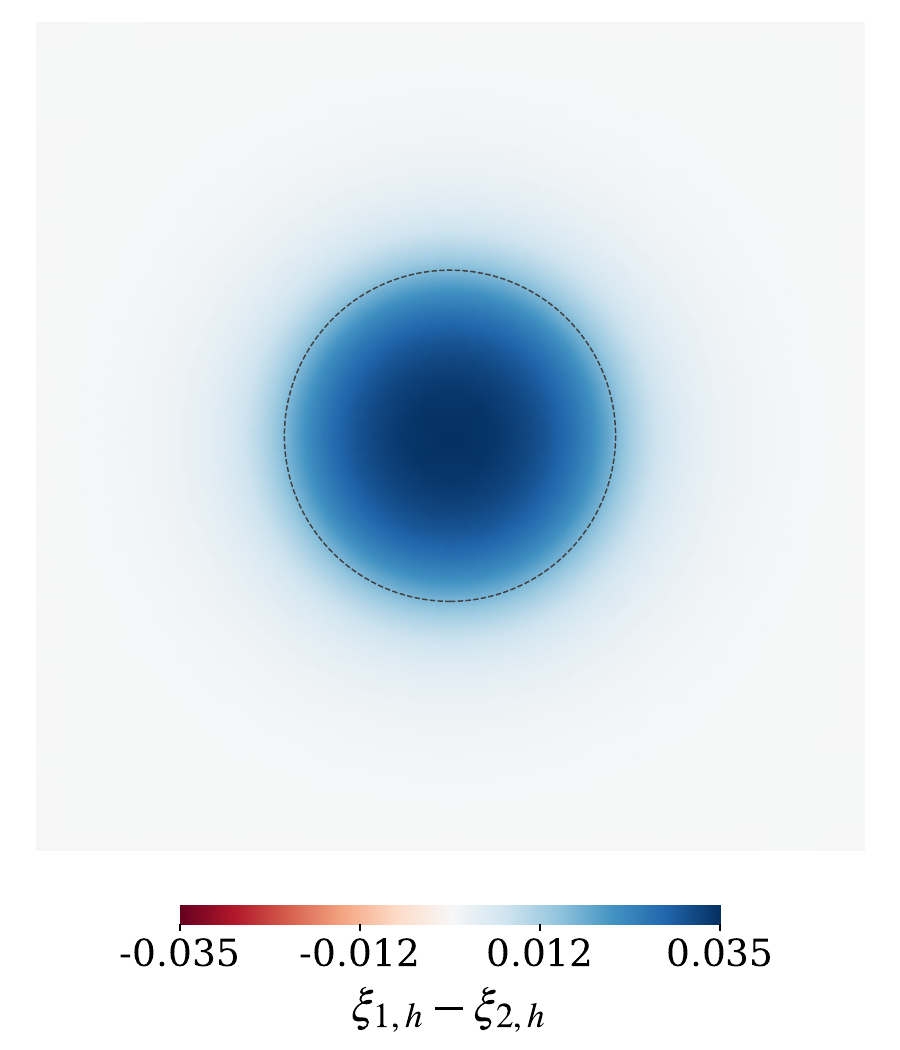}
\caption{Example~3: final adapted mesh after four refinements
(left) and mobile-charge difference at $V=0$ (right).}
\label{fig:numerical-charged-mesh}
\end{figure}
\begin{figure}[!t]
\centering
\includegraphics[width=0.27\textwidth]{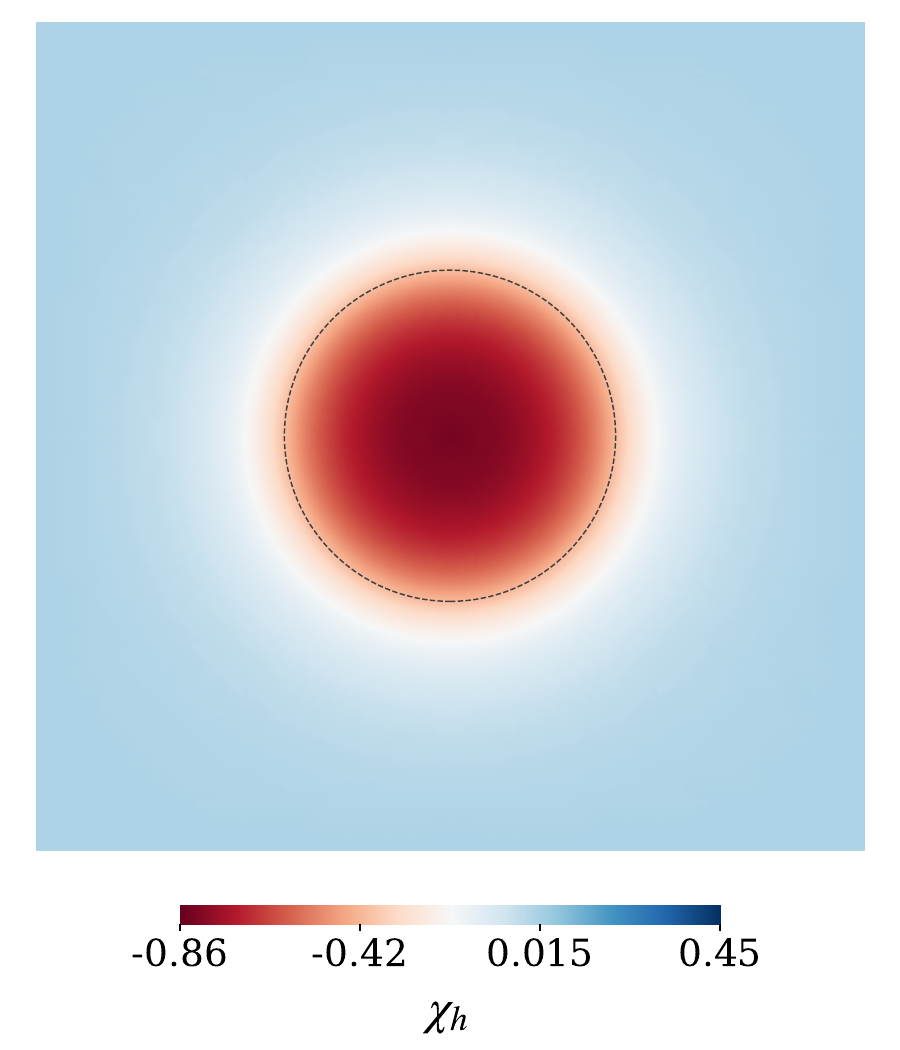}
\hfill
\includegraphics[width=0.27\textwidth]{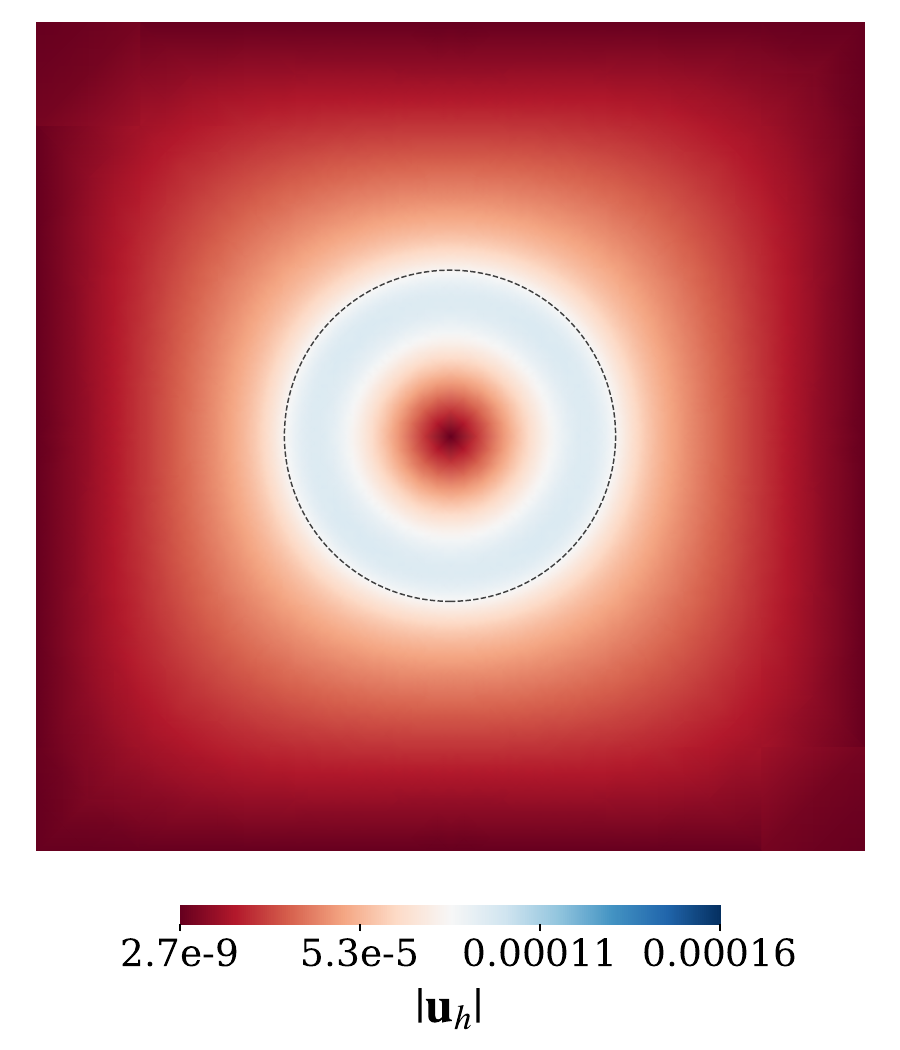}
\hfill
\includegraphics[width=0.27\textwidth]{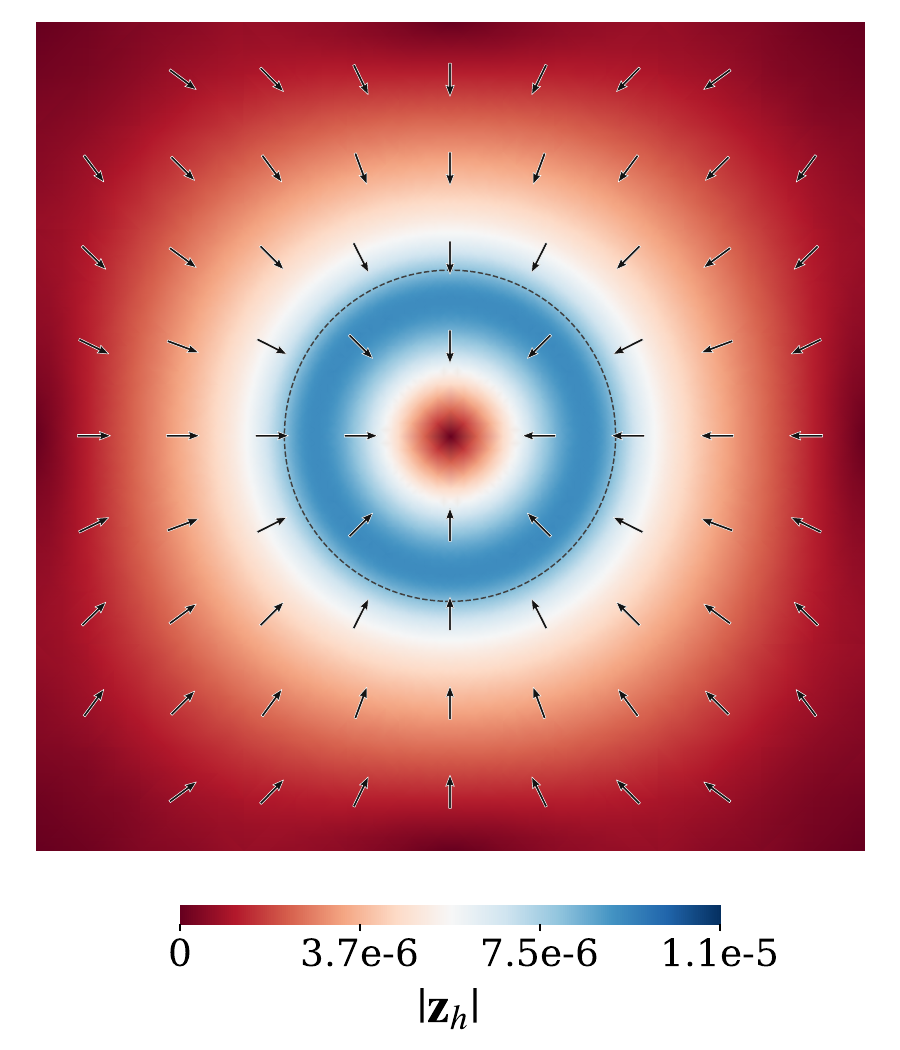}
\\[1.5ex]
\includegraphics[width=0.27\textwidth]{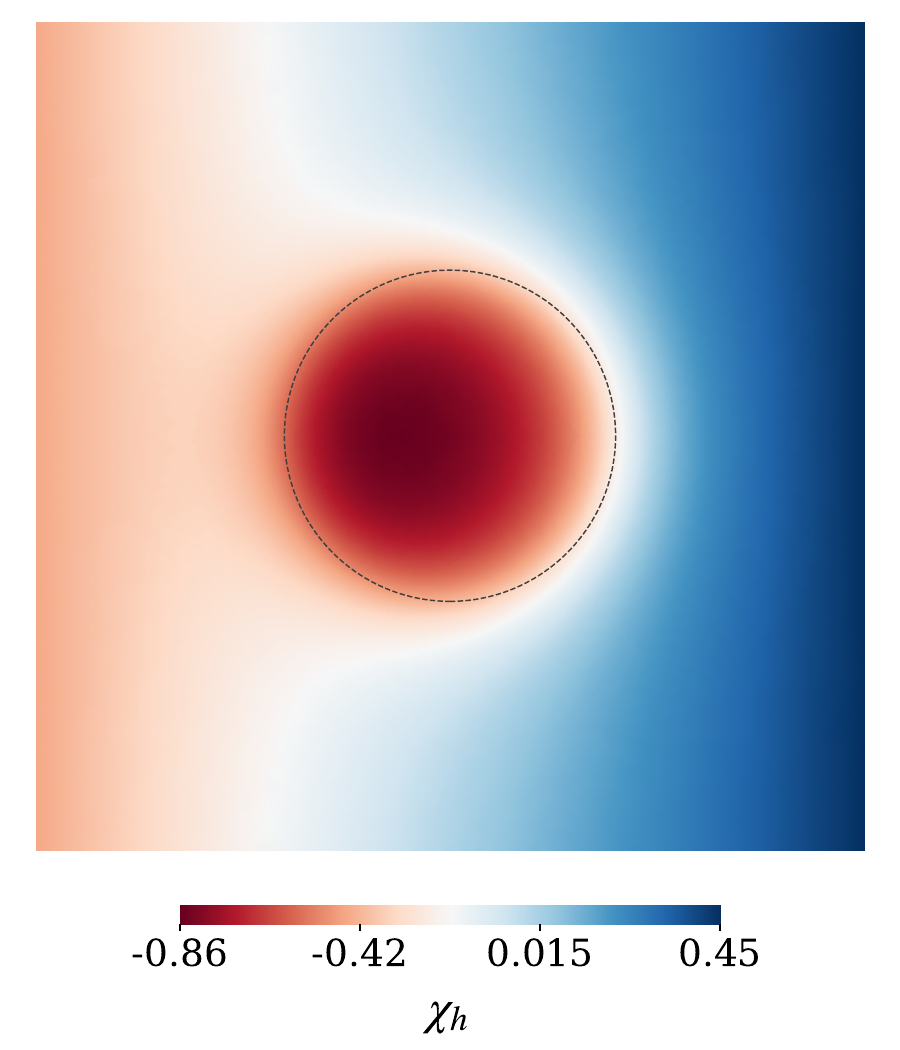}
\hfill
\includegraphics[width=0.27\textwidth]{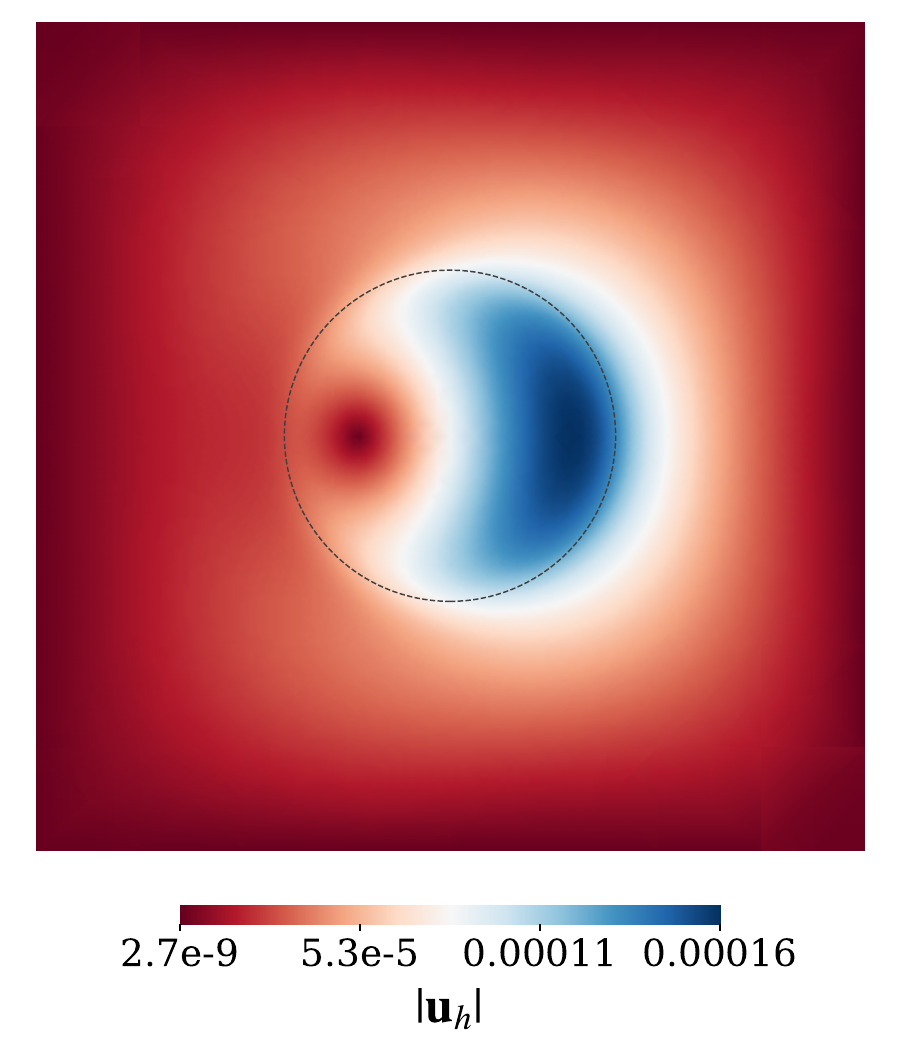}
\hfill
\includegraphics[width=0.27\textwidth]{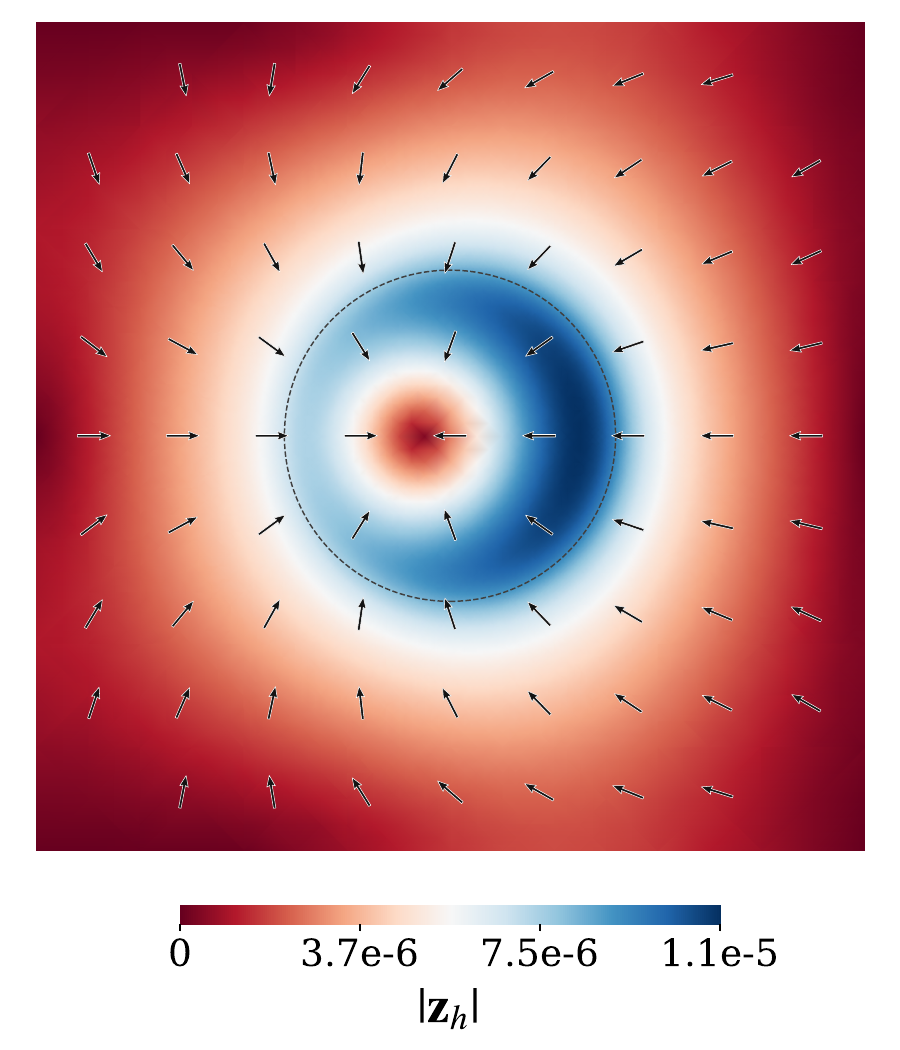}
\caption{Example~3: electric potential, displacement magnitude, and
Darcy-flux magnitude (left to right) for $V=0$ (first row) and $V=0.9$
(second row).}
\label{fig:numerical-charged-flux}
\end{figure}
\begin{figure}[!t]
\centering
\includegraphics[width=0.6\textwidth]{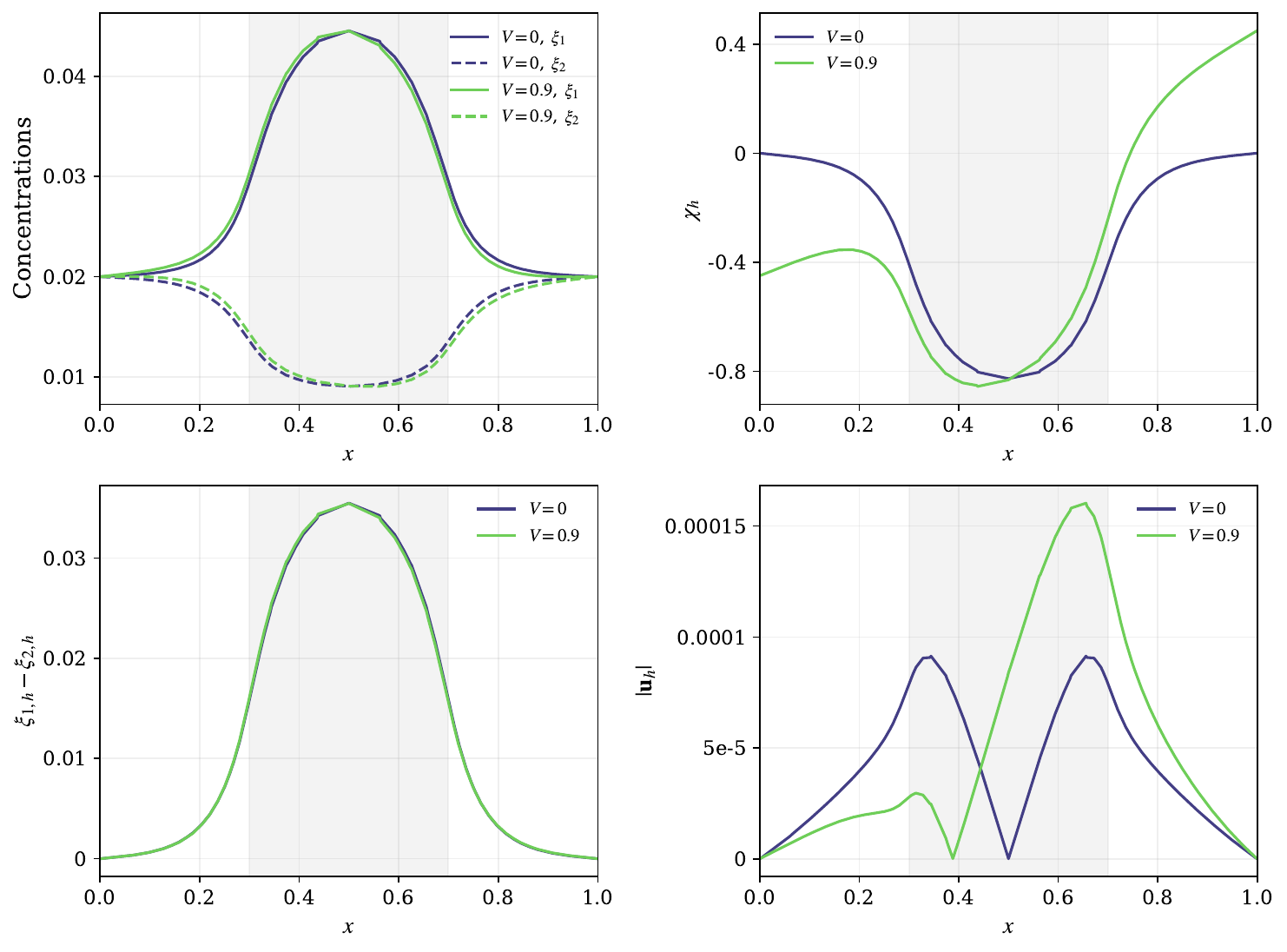}
\caption{Example~3: profiles along $y=1/2$ on the final meshes for
$V=0$ and $V=0.9$: ion concentrations (top left), electric potential
(top right), mobile-charge difference (bottom left), and displacement
magnitude (bottom right). The shaded band marks the charged region.}
\label{fig:numerical-charged-profiles}
\end{figure}
We compare $V=0$ and $V=0.9$ in
\eqref{eq:numerical-charged-boundary}.
In this example the fixed charge and the voltage are applied by load
continuation in six steps, starting from the uncharged rest state
$\xi_1=\xi_2=0.02$ (a failed step is retried with half the increment).
At each step, Newton's method stops when the coefficient criterion above
holds and the algebraic residual is below $10^{-8}$ times its initial value.

We use $\mathbb{RT}_1$--$\mathbf P_1$, starting from
a uniform mesh of $128$ triangles, and perform four
adaptive refinements driven by $\Theta$
(cf.~\eqref{def:posteriori-global-estimator}),
with $C_{\tt adm}=0.8$.
The final meshes have $8{,}020$ and $7{,}938$ triangles, with
$481{,}490$ and $476{,}558$ DOF for $V=0$ and $V=0.9$, respectively.
On these meshes, the six continuation steps require $24$ Newton
updates in total for each voltage, with no retried step, and the
final algebraic residuals are below $4\times10^{-14}$.
The final estimator values are $\Theta\approx2.5\times10^{-2}$ for both
voltages, and the Poisson indicators account for about $96\%$ of them.
Refinement is therefore driven mainly by the electrostatic response
at the fixed-charge interface.

Figure~\ref{fig:numerical-charged-mesh} shows the final
adapted mesh and mobile-charge difference for $V=0$.
Refinement concentrates near the transition of the
prescribed fixed charge
(cf.~\eqref{eq:numerical-fixed-charge}).
The positive mobile-charge difference in the central
region indicates an excess of positive over negative
mobile ions, consistent with the negative prescribed
fixed charge.
Figure~\ref{fig:numerical-charged-flux} compares
the potential, displacement magnitude, and Darcy-flux
magnitude on the final meshes.
For $V=0$, the fields are nearly symmetric about
the centre: the potential is negative in the charged
region, while the displacement and flux magnitudes
are largest in an annular region near the fixed-charge
transition (cf.~\eqref{eq:numerical-fixed-charge}).
The imposed voltage introduces a horizontal variation
in the potential and an asymmetry in the displacement
and flux magnitudes, with larger values on the
right-hand side.
Compared with $V=0$, the computed displacement norm
$\|\u_h\|_{0;\O}$ increases by approximately $20.2\%$
and the flux norm $\|\z_h\|_{0,4;\O}$ by $6.4\%$.

Figure~\ref{fig:numerical-charged-profiles} shows the computed
profiles along the horizontal centreline $y=1/2$.
Inside the charged region, $\xi_{1,h}$ reaches about $0.045$ and
$\xi_{2,h}$ decreases to about $0.009$, from the boundary value
$0.02$. The resulting mobile-charge difference, about $0.035$ at the
centre, screens most of the fixed charge $C_{\rm f}=0.04$.
The concentration profiles are almost insensitive to the imposed
voltage. In contrast, the potential follows the imposed boundary
ramp, and the displacement magnitude loses its symmetry: its maximum
moves to the right-hand transition and increases from about
$0.9\times10^{-4}$ to $1.6\times10^{-4}$.

\FloatBarrier
\endgroup

\section*{Declarations}

\paragraph{Funding.}
This research was partially supported by ANID-Chile through
Fondecyt Postdoctoral project 3260398.

\paragraph{Competing interests.}
The author has no relevant financial or non-financial interests
to disclose.

\paragraph{Ethics approval.}
Not applicable.

\paragraph{Consent to participate.}
Not applicable.

\paragraph{Consent for publication.}
Not applicable.

\paragraph{Data availability.}
Data sharing is not applicable to this article, as no datasets were
generated or analysed beyond the numerical results reported in
Section~\ref{sec:numerical-results}.

\paragraph{Code availability.}
The {\tt FEniCS} codes used to produce the numerical results are
available from the author upon reasonable request.

\paragraph{Author contributions.}
C.~Inzunza is the sole author and carried out the analysis, the
implementation, the numerical experiments, and the writing of the
manuscript.

\end{document}